\documentclass[12pt]{article}
\usepackage{algorithmic}
\usepackage{amsfonts}
\usepackage{amsmath}
\usepackage{amssymb}
\usepackage{bbm}
\usepackage{color}
\usepackage{eurosym}
\usepackage{graphicx}
\usepackage{latexsym}
\usepackage[T1]{fontenc}
\usepackage{upgreek}

\newtheorem{cs}{Case Study}
\newtheorem{problem}{Problem}
\newtheorem{remark}{Remark}

\newcommand{\bfh}{\mbox{$\mbox{\boldmath $h$}$}} 

\newcommand{\bfr}{\mbox{$\mbox{\boldmath $r$}$}} 
\newcommand{\bfs}{\mbox{$\mbox{\boldmath $s$}$}} 
\newcommand{\bft}{\mbox{$\mbox{\boldmath $t$}$}} 
\newcommand{\bfu}{\mbox{$\mbox{\boldmath $u$}$}} 
\newcommand{\bfv}{\mbox{$\mbox{\boldmath $v$}$}} 
\newcommand{\bfx}{\mbox{$\mbox{\boldmath $x$}$}} 
\newcommand{\bfy}{\mbox{$\mbox{\boldmath $y$}$}} 
\newcommand{\bfz}{\mbox{$\mbox{\boldmath $z$}$}} 

\newcommand{\sbfh}{\mbox{$\mbox{\boldmath \scriptsize $h$}$}} 

\newcommand{\bfU}{\mbox{$\mbox{\boldmath $U$}$}} 
\newcommand{\bfV}{\mbox{$\mbox{\boldmath $V$}$}}
\newcommand{\bfW}{\mbox{$\mbox{\boldmath $W$}$}}

\newcommand{\bnab}{\mbox{$\mbox{\boldmath $\nabla$}$}} 
\newcommand{\beps}{\mbox{$\mbox{\boldmath $\epsilon$}$}} 

\newcommand{\bsigma}{\mbox{$\mbox{\boldmath $\sigma$}$}} 

\newcommand{\sbsigma}{\mbox{$\mbox{\boldmath \scriptsize $\sigma$}$}} 

\newcommand{\bfzero}{\mbox{$\mbox{\boldmath $0$}$}} 
\newcommand{\bfone}{\mbox{$\mbox{\boldmath $1$}$}}

\begin{document}

\title{Optimal driving strategies for a fleet of trains on level track with prescribed intermediate signal times and safe separation}

\author{Phil Howlett\footnote{Email: phil.howlett@unisa.edu.au} \and Peter Pudney\footnote{Email: peter.pudney@unisa.edu.au} \and Amie Albrecht\footnote{Email: amie.albrecht@unisa.edu.au}}

\date{11 February 2022}

\maketitle

\begin{abstract}
We propose an analytic solution to the problem of finding optimal driving strategies that minimize total tractive energy consumption for a fleet of trains travelling on the same track in the same direction subject to clearance-time equality constraints that ensure safe separation and compress the line-occupancy timespan.  We assume the track is divided into sections by a set of trackside signals at fixed locations.  For each intermediate signal there is a signal-location segment consisting of the two adjacent sections.  Successive trains are safely separated only if the leading train leaves each signal-location segment before the following train enters.  The fleet can be safely separated by a complete set of clearance times and associated clearance-time inequality constraints.  The problem of finding optimal schedules with safe separation has been solved for two trains but for larger fleets the problem rapidly becomes intractable as the number of trains and signals increases.  The main difficulty is in distinguishing between active equality constraints and inactive inequality constraints. The curse of dimensionality means it is not feasible to check every different combination of active constraints, optimize the corresponding prescribed times and calculate the cost.  Nevertheless we can formulate and solve an alternative problem with active clearance-time equality constraints for successive trains on every signal-location segment.  We show that this problem can be formulated as an unconstrained convex optimization and we propose a viable solution algorithm that finds the optimal schedule and the associated optimal strategies for each train.  Finally we use our solution to find optimal schedules for a busy inter-city shuttle service.
\end{abstract}


\maketitle

\section{Introduction}
\label{s:int}

In modern rail networks train movements are planned to follow strict timetables that allow coordinated operation of the entire network.  Drivers are encouraged to use energy-efficient driving strategies but are expected to reach key locations at predetermined times or within prescribed time windows.  Despite these expectations network safety is paramount and drivers must always comply with safe-operating instructions from the signalling system.  

A safe-operating environment has traditionally been enabled using a system of trackside signals at fixed locations.   These signals divide the track into sections.  A three-aspect signalling system shows a green light if the next two sections of track are clear, a yellow light if the next section is clear but the section after that is occupied, and a red light if the next section is occupied.  A driver will follow the planned schedule when the train passes a green signal but if the train passes a yellow signal the driver must follow a modified speed profile so that the train can stop at the next signal if it remains red.  The train must not pass a red signal.  In normal operation trains will be separated  by at least two fixed signals and one intervening clear section of track.  The signalling system ensures that trains are adequately warned about unscheduled disruptions.

In this paper we wish to find an optimal schedule for a fleet of trains travelling on the same track in the same direction.  The trains must remain safely separated at all times.  In order to clarify our presentation the remainder of the introduction is organised into separate sections\textemdash each one with a specific objective.  Some relate to the problem at hand while others relate to essential background material.  In Section~\ref{s:ppf} we state a preliminary version of the problem.  In Section~\ref{s:sss} we explain how intermediate segment clearance times can be used to ensure that successive trains are safely separated.  Section~\ref{s:em} introduces the technical basis for the paper\textemdash the equations of motion\textemdash and Section~\ref{s:sot} describes the general structural forms for the strategies of optimal type.  These strategies provide the framework for a theory of optimal scheduling. Section~\ref{s:ttsp} describes a typical solution to the two-train separation problem.  In Section~\ref{s:gftsp} the general train separation problem is formulated as two distinct problems.  The first problem is to find an optimal strategy for each train if the stopping patterns are known and the prescribed signal-location times are also known.  This problem has been solved.  The second problem is to find a set of optimal prescribed intermediate clearance times that ensures safe separation and minimizes total energy consumption for the fleet.  This problem has not been solved in any realistic sense because the number of combinations of active constraints which must be checked increases exponentially as the number of trains and signals increases. In Section~\ref{s:mtr} we describe our main theoretical result\textemdash formulation and solution of an alternative train separation problem where a fleet of non-identical trains is separated by a complete set of active clearance-time equality constraints and the task is to find the optimal clearance times.  We explain briefly how our solution can be extended to include additional buffer times between trains and how it can be relaxed by omitting some constraints.  Section~\ref{s:atr} simply notes that we apply our results to a case study of a busy inter-city shuttle service.  Section~\ref{s:rtpf} discusses the role of the train performance functions and Section~\ref{s:term} reviews our terminology.   

\subsection{A preliminary problem formulation}
\label{s:ppf}

We consider a fleet of $m$ successive trains travelling in the same direction on a level track with $n+1$ signals at locations $0 = x_0 < x_1 < \cdots < x_{n-1} < x_n = X$.  For each $i=1,\ldots,m$ the journey for train ${\mathfrak T}_i$ is completely defined by the speed $v_i(x)$ and the elapsed journey time $t_i(x)$ at position $x \in [0,X]$.  We assume that the schedule for ${\mathfrak T}_i$ is defined by a vector $\bfh_i = [h_{i,j}]_{j=0}^n$ of signal location times with $t_i(x_j) = h_{i,j}$ for each $i=1,\ldots,m$ and $j=0,\ldots,n$.  If ${\mathfrak T}_i$ is scheduled to stop at $x_j$ then we assume that $h_{i,j}$ is the departure time and that the stopping time $\sigma_{i,j}$ is known.  We wish to find a schedule $\bfh = [\bfh_1,\ldots,\bfh_m]$ for the entire fleet that minimizes total energy consumption, allows each train to finish the journey on time and ensures that successive trains are safely separated at all times.

\subsection{A safe separation scheme}
\label{s:sss}
  
For each $j=1,\ldots,n$ we define a signal-location segment ${\mathfrak S}_j = [x_{j-1},x_j] \cup [x_j, x_{j+1}] = [x_{j-1}, x_{j+1}]$.  Safe separation for successive trains ${\mathfrak T}_i$ and ${\mathfrak T}_{i+1}$ can be enforced by prescribing segment-clearance times that define the latest allowed exit time for ${\mathfrak T}_i$ and the earliest allowed entry time for ${\mathfrak T}_{i+1}$.  The leading train ${\mathfrak T}_i$ must leave the segment before the following train ${\mathfrak T}_{i+1}$ enters.  Thus we require $h_{i,j+1} = t_i(x_{j+1}) \leq t_{i+1}(x_{j-1}) = h_{i+1,j-1}$ for all $i=1,\ldots,m-1$ and $j=1,\ldots,n-1$.  We could define a safe theoretical schedule by setting $h_{i+1,j-1} = h_{i,j+1}$ but in practice train operators allow for normal stochastic variation in journey times by inserting additional buffer times.  Thus, in practice, we may wish to set $h_{i+1, j-1} = h_{i, j+1} + \delta_i$ where $\delta_i > 0$ for each $i=1,\ldots,m-1$ and each $j=1,\ldots,n-1$.

\subsection{The equations of motion for realistic strategies}
\label{s:em}

We formulate the equations of motion with position $x \in [0,X]$ as the independent variable and with speed $v = v(x) \in [0, \infty)$ and time $t = t(x) \in [0,T]$ as dependent state variables.  The equations are
\begin{eqnarray}
v^{\, \prime} & = & [u - r(v) + g(x)]/v \label{em:vx} \\
t^{\, \prime} & = & 1/v \label{em:tx}
\end{eqnarray}
where $(v,t) = (v(x),t(x))$ is the state variable vector for $x \in [0,X]$ and where $u = u(x) \in {\mathbb R}$ is the known measurable control\textemdash the force per unit mass or acceleration.  We have written $v^{\, \prime} = dv/dx$ and $t^{\, \prime} = dt/dx$.  This formulation decouples the dependent variables and allows us to solve (\ref{em:vx}) for $v = v(x)$ with no knowledge of $t = t(x)$.  We measure distance in metres (m) and time in seconds (s). 

For a journey from $x_{k-1}$ to $x_{\ell}$ where $1 \leq k \leq \ell \leq n$ we assume that $v(x_{k-1}) = v(x_{\ell}) = 0$ with $v(x) > 0$ for all $x \in (x_{k-1},x_{\ell})$ and that $u(x)$ is bounded with $K[v(x)] \leq u(x) \leq H[v(x)]$ for each $x \in (x_{k-1},x_{\ell})$. The bounds $K = K(v) \in(-\infty,0)$ and $H = H(v) \in (0,\infty)$ for $v \in (0, \infty)$ are monotone functions with $K(v) \uparrow 0$ and $H(v) \downarrow 0$ as $v \uparrow \infty$.  The functions $K$ and $H$ define bounds for the maximum braking and driving forces per unit mass in a form that includes\textemdash as special cases\textemdash the specified bounds for a wide range of modern electric and diesel-electric locomotives.  The function $r(v)$ is a general resistance per unit mass with no specific formula assumed.  We define auxiliary functions $\varphi(v) = vr(v)$ and $\psi(v) = v^2 r^{\, \prime}(v)$ and assume only that $\varphi(v)$ is strictly convex with $\varphi(v) \geq 0$ for $v \geq 0$ and $\varphi(v)/v \rightarrow \infty$ as $v \rightarrow \infty$.  See~\cite[Appendix A.3, p 409]{how6} for an explanation.  It follows that both $r(v)$ and $\psi(v)$ are non-negative and strictly increasing for $v \geq 0$.  These properties capture the functional characteristics of the traditional quadratic resistance formula\textemdash the so-called Davis formula \cite{dav1}.  The function $g(x)$ is nominally the component of gravitational acceleration due to track gradient but in practice may also include additional position-dependent resistive forces.  The cost is net mechanical energy usage per unit mass, 
\begin{equation}
\label{em:c}
J_{k-1,\ell} = \int_{x_{k-1}}^{x_{\ell}} (1/2) [u(x) + |u(x)|] \, dx.
\end{equation}
See \cite{alb4, alb5, how6} and the original papers \cite{how7, khm1, liu1} for more information.  The equations (\ref{em:vx}) and (\ref{em:tx}) describe the motion of a point-mass train.  It is known \cite[Section 2.3, pp 23\textendash 24]{how4} that the motion of a train with distributed mass can be modelled as the motion of a point-mass train on a track with modified gradient.  In this paper we assume that energy recovered from regenerative braking is not used to drive the train.  This assumption spawns a more relaxed optimal driving strategy that encourages coasting and discourages braking.  We also assume that the track is level with $g(x) = 0$ for all $x \in [0,X]$.  

\subsection{The strategies of optimal type}
\label{s:sot} 

For a fleet of $m$ trains with scheduled signal times $h_{i,j}$ at each $x_j$ the optimal strategy for train ${\mathfrak T}_i$ will be a sequence of optimal strategies for a succession of smaller journeys between scheduled stops.  Some of these strategies will be simple strategies for journeys on single sections $[x_{k-1},x_k]$ with initial and final speeds $v_i(x_{k-1}) = v_i(x_k) = 0$ and with initial and final time constraints $t_i(x_{k-1}) = h_{i,k-1}$ and $t_i(x_k) = h_{i,k} - \sigma_{i,k}$ where $\sigma_{i,k}$ is the stopping time at $x_k$ but with no intermediate time constraints.  Others will be more complex strategies for journeys on composite segments $[x_{k-1},x_{\ell}]$ where $k < \ell$ with initial and final speeds $v_i(x_{k-1}) = v_i(x_{\ell}) = 0$, initial and final time constraints $t_i(x_{k-1}) = h_{i,k-1}$ and $t_i(x_{\ell}) = h_{i,\ell} - \sigma_{i,\ell}$ where $\sigma_{i,\ell}$ is the stopping time at $x_{\ell}$ and intermediate time constraints $t_i(x_j) = h_{i,j}$ for each $j=k,\ldots,\ell-1$ but no intermediate stops. 

The Pontryagin principle shows that only certain control modes are allowed in an optimal strategy \cite[Section 3, pp 489\textendash 507]{alb4}.  There are three permissible regular modes\textemdash maximum acceleration, coast and maximum brake\textemdash and one singular mode\textemdash speedhold with partial acceleration.  For trains with regenerative braking there is another permissible singular mode\textemdash speedhold with partial brake.  These rules are independent of the initial and final speeds and apply on both level and nonlevel tracks.  A strategy that is generated entirely by a sequence of permissible optimal control modes is called a strategy of optimal type.

For journeys with $v(0) = v(X) = 0$ and no intermediate time constraints the minimum possible journey time $T_{\min}$ has only two phases\textemdash maximum acceleration to some speed $V_{\max}$ and maximum brake.  For $T > T_{\min}$ there is always a rapid-transit strategy of optimal type consisting of a phase of maximum acceleration to speed $V_{\max}$, coast to speed $U$ and maximum brake.  The values of $V_{\max}$ and $U$ are uniquely determined by the distance and time constraints
\begin{equation}
\label{d0Xrt}
X = \int_0^{V_{\max}} dx_a(v) + \int_U^{V_{\max}} |dx_c(v)| + \int_0^U |dx_b(v)|
\end{equation}
and
\begin{equation}
\label{t0Xrt}
T = \int_0^{V_{\max}} dt_a(v) + \int_U^{V_{\max}} |dt_c(v)| + \int_0^U |dt_b(v)|
\end{equation}
where the distance and time differentials for maximum acceleration, coast and maximum brake are defined respectively from (\ref{em:vx}) and (\ref{em:tx}) by
\begin{equation}
\label{dxav}
dx_a(v) = vdt_a(v) = vdv/[H(v) - r(v)],
\end{equation}
\begin{equation}
\label{dxcv}
dx_c(v) = vdt_c(v) = - vdv/r(v)
\end{equation}
and
\begin{equation}
\label{dxbv}
dx_b(v) = vdt_b(v) = -vdv/[K(v)+r(v)]
\end{equation}
for $v \in [0, V_{\sup})$ where $v = V_{\sup}$ is the upper bound on the speed obtained by solving the equation $H(v) - r(v) = 0$.  The rapid-transit strategy is optimal only if $U > U_b(V_{\max})$ where
\begin{equation}
\label{ub}
U_b(v) = \psi(v)/\varphi^{\, \prime}(v).
\end{equation}
If the rapid-transit strategy is not optimal then the optimal strategy is a long-haul strategy of optimal type defined by a phase of maximum acceleration to speed $V$, speedhold at speed $V$, coast to speed $U = U_b(v)$ and maximum brake.  The length of the speedhold segment is defined by
\begin{equation}
\label{dshlh}
\xi(V) = X - \left[ \int_0^V dx_a(v) + \int_{U_b(V)}^V |dx_c(v)| + \int_0^{U_b(V)} |dx_b(v)| \right]
\end{equation}
and $V$ is uniquely determined by the time constraint
\begin{equation}
\label{t0Xlh}
T = \int_0^V dt_a(v) + \xi(V)/V + \int_{U_b(V)}^V |dt_c(v)| + \int_0^{U_b(V)} |dt_b(v)|.
\end{equation}
The Pontryagin principle shows that each strategy of optimal type is uniquely defined by an optimal driving speed.  For a long-haul strategy the optimal driving speed is the speed $V$ on the speedhold segment.  If the rapid-transit strategy is optimal then the optimal driving speed is defined by $\psi(V) = \varphi(V_{\max})U/(V_{\max}-U) \iff U = V_{\max} - V_{\max}\varphi(V_{\max})/[\psi(V) + \varphi(V_{\max})]$.  In this case $V > V_{\max}$ and $U > U_b(V_{\max})$.  See \cite[Section 3.6, pp 494\textendash 497]{alb4} for a detailed discussion.      

Now consider a journey with initial and final speeds $v(0) = v(X) = 0$, initial and final times $t(0)  = 0$ and $t(X) = T$, and prescribed intermediate times $t(x_j) = h_j$ for each $j=1,\ldots,n-1$ but no intermediate stops.  If the section times $h_j - h_{j-1}$ are sufficiently large for each $j=1,\ldots,n$ then an optimal strategy exists and there is an optimal driving speed $V_j$ on each timed section.  On $(x_0, x_1)$ the strategy is maximum acceleration, speedhold at speed $V_1$ on some interval $[a_1,b_1] \subset (x_0,x_1)$ and either maximum acceleration or coast to speed $v(x_1) = U_1 = U_s(V_1,V_2)$ where
\begin{equation}
\label{us}
U_s(v,w) = [\psi(v) - \psi(w)]/[\varphi^{\, \prime}(v) - \varphi^{\, \prime}(w)].
\end{equation}
In general, on $(x_{j-1},x_j)$ for $j=2,\ldots,n-1$ the strategy is either maximum acceleration or coast from speed $v(x_{j-1}) = U_{j-1} = U_s(V_{j-1},V_j)$ to speed $V_j$, speedhold at speed $V_j$ on some interval $[a_j,b_j] \subset (x_{j-1},x_j)$, and either maximum acceleration or coast to speed $v(x_j) = U_j = U_s(V_j,V_{j+1})$.  On $(x_{n-1},x_n)$ the strategy is either maximum acceleration or coast from speed $v(x_{n-1}) = U_{n-1} = U_s(V_{n-1},V_n)$ to speed $V_n$, speedhold at speed $V_n$, coast to speed $U_n = U_b(V_n)$ and maximum brake.  The transition phase on $(b_j,a_{j+1})$ is maximum acceleration if $V_{j+1} > V_j$ and coast if $V_{j+1} < V_j$.  If there is insufficient time for a speedhold phase on the final section $(x_{n-1},x_n)$ of the journey the strategy is maximum acceleration from speed $v(x_{n-1}) = U_{n-1} = U_s^{\dag}(U_n,V_{n-1},V_{\max,n})$ to a maximum speed $V_{\max,n}$, coast to speed $U_n$ and maximum brake where
\begin{equation}
\label{usd}
U_s^{\dag}(u,v,w) = \frac{\varphi(w)u/(w - u) - \psi(v)}{\varphi(w)/(w-u) - \varphi^{\, \prime}(v)}.
\end{equation}
The optimal driving speeds $\{V_j\}_{j=1}^{n-1}$ and the speeds $V_{\max,n}$ and $U_n$ are determined by the time constraints on each section and the distance constraint on the final section.  We still have $v(x_j) = U_j = U_s(V_j,V_{j+1})$ for $j =1,\ldots,n-2$.  The optimal driving speed $V_n$ on $(x_{n-1},x_n)$ is defined by $\psi(V_n) = \varphi(V_{\max,n})U_n/(V_{\max,n} - U_n)$.  We also have $V_n > V_{\max,n}$ and $U_n > U_b(V_{\max,n})$.

\subsection{The two-train separation problem}
\label{s:ttsp}

\begin{figure}[htb]
\begin{center}
\includegraphics[width=5.4cm]{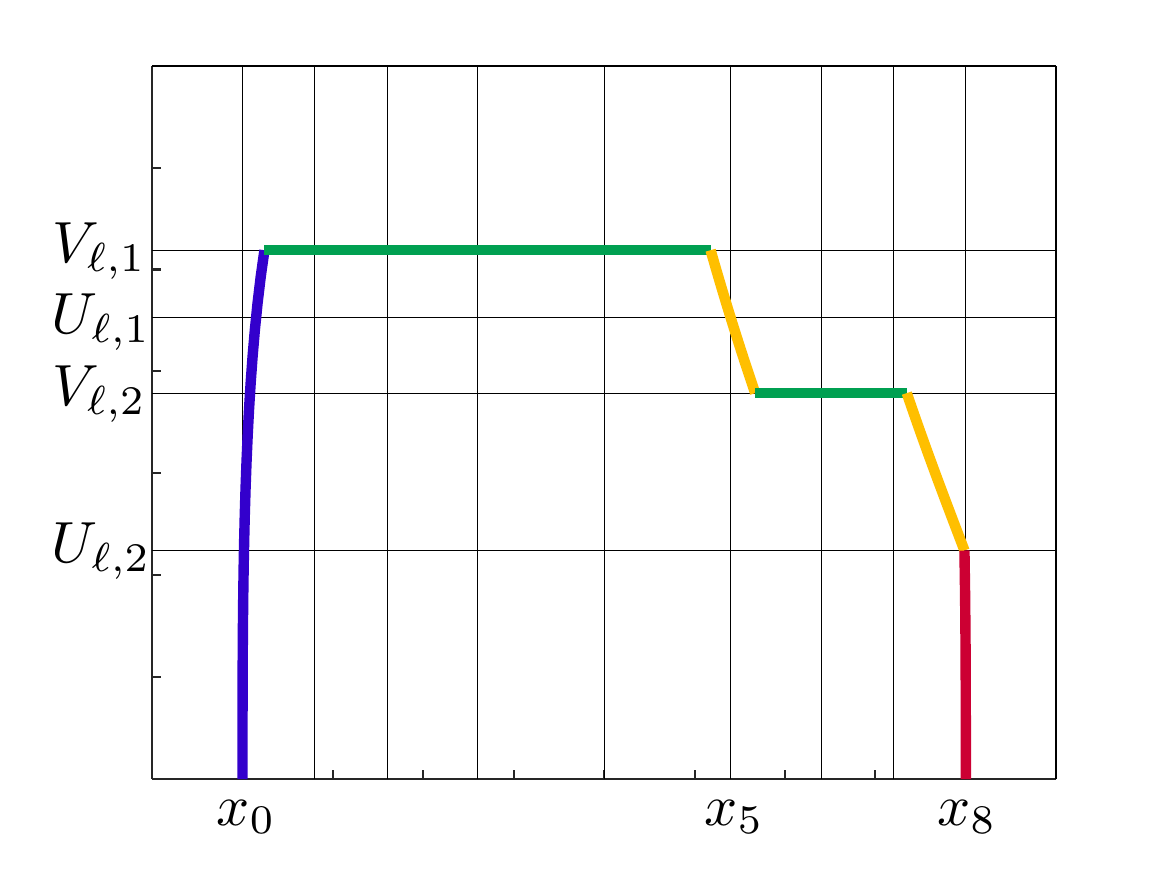}
\includegraphics[width=5.4cm]{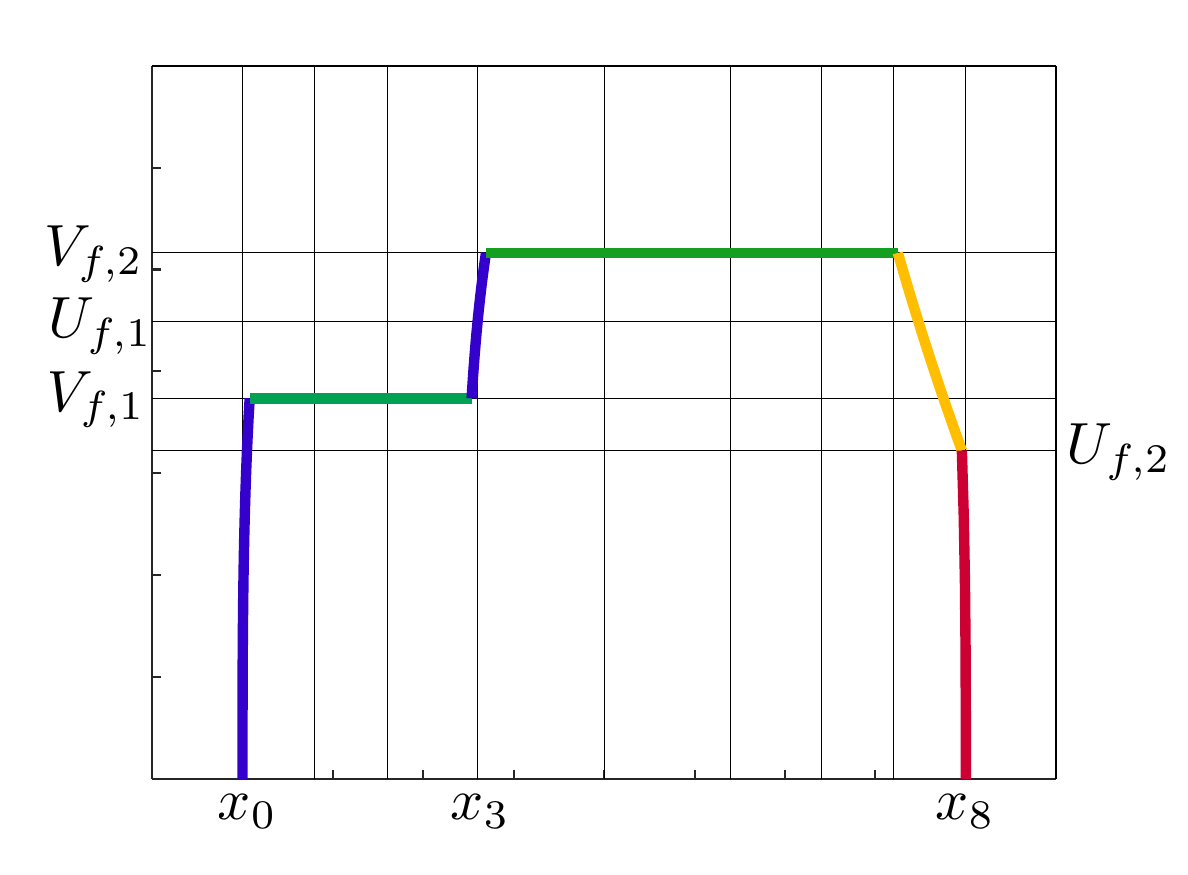}
\includegraphics[width=5.4cm]{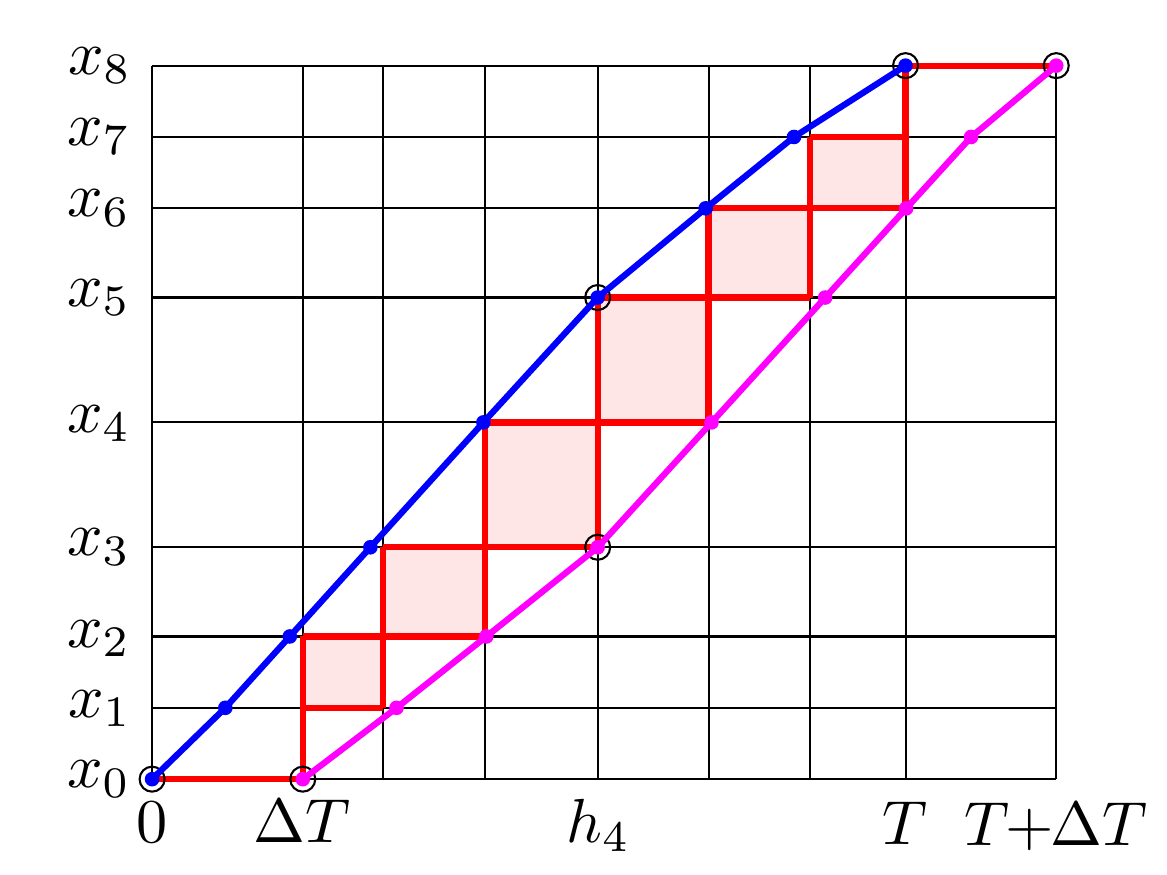}
\end{center}
\caption{\small \boldmath \bf A typical two train separation problem with distance $X$ and time $T$ and signal locations $\bfx = (x_0=0, x_1,\ldots, x_7, x_8=X)$ showing optimal speed profiles for the leading train~${\mathfrak T}_{\ell}$ (left) and following train~${\mathfrak T}_f$ (centre).  The train graph (right) shows one clear section between ${\mathfrak T}_{\ell}$ and ${\mathfrak T}_f$ at all times.  There is one active time constraint on ${\mathfrak S}_4$ with $t_{\ell}(x_5) = t_f(x_3) = h_4$.}
\label{fig1}
\end{figure}

Much of the train separation work has been done for a fleet of two trains.  For each train the optimal strategy is defined by a sequence of optimal driving speeds.  For the leading train $\{V_{\ell,j}\}_{j=1}^n$ is a decreasing sequence and for the following train $\{V_{f,j}\}_{j=1}^n$ is an increasing sequence.  See \cite{alb9}.  The graphs in Figure~\ref{fig1} depict the optimal solution to a two-train separation problem on level track \cite[Section~VII, Example 3, Figure~4]{alb9}.  The speed profiles are shown on the left for the leading train~${\mathfrak T}_{\ell}$ and in the centre for the following train~${\mathfrak T}_f$.  Only one active time constraint is required to ensure safe separation.  The optimal clearance time $t_{\ell}(x_5) = t_f(x_3) = h_4$ on segment ${\mathfrak S}_4 = [x_3,x_4] \cup [x_4, x_5] = [x_3,x_5]$ minimizes total energy consumption.  The train graph shows that the trains are safely separated.  

When speed limits are imposed the optimal strategies are not changed substantially but the speed limits must be enforced~\cite{how3, pud1}.  We will not consider speed limits here but we note that the very nature of optimal strategies ensures that speeds must be kept as low as possible.  Thus virtual speed limits are imposed.

\subsection{A general formulation of the train separation problem}
\label{s:gftsp}

The task of designing optimal strategies for a fleet of trains travelling on the same track in the same direction subject to conditions of safe separation can be formulated as two distinct problems.

\begin{problem}
\label{tsp:p1}
Let $\bfh = [\bfh_1,\ldots,\bfh_m]$, where $\bfh_i = [h_{i,j}]_{j \in {\mathcal C}_i}$ is a given set of prescribed signal times for train ${\mathfrak T}_i$ on some subset ${\mathcal C}_i \subset \{0,\ldots,n\}$ of signal-location indices which includes the indices of all scheduled stops.  Suppose the set of clearance-time constraints $t_i(x_j) = h_{i,j}$ for each $i =1,\ldots,m$ and each $j \in {\mathcal C}_i$ is sufficient to ensure safe separation.  Let ${\mathcal U}(\bfh)$ be the set of feasible controls.  Find a control vector $\bfu = \bfu(x)$ where  $\bfu = [u_1,\ldots,u_m] \in {\mathcal U}(\bfh)$ and $u_i = u_i(x)$ is the control function for ${\mathfrak T}_i$ and the associated speed profile vector $\bfv = \bfv(x)$ where $\bfv = [ v_1,\ldots,v_m]$ and $v_i = v_i(x)$ is the speed profile for ${\mathfrak T}_i$, such that the constraints are satisfied and tractive energy consumption $J = J(\bfu, \bfh)$ is minimized.  That is find $\bfu = \bfu_{\sbfh}$ such that $J(\bfh) = J(\bfu_{\sbfh}, \bfh) = \min \{ J(\bfu, \bfh) \mid \bfu \in {\mathcal U}(\bfh)\}$. $\hfill \Box$
\end{problem}

\begin{problem}
\label{tsp:p2}
Find necessary conditions for a vector $\bfh_0$ of optimal prescribed clearance times that minimizes the total tractive energy $J(\bfh)$ consumed by the fleet over all $\bfh \in {\mathcal H}$, where ${\mathcal H}$ denotes the collection of all feasible sets of prescribed clearance-time constraints that ensure safe separation for the given stopping plan.  Thus $J_0 = J(\bfh_0) = \min \{ J(\bfh) \mid \bfh \in {\mathcal H} \}$.  $\hfill \Box$
\end{problem}

Problem~\ref{tsp:p1} has been solved.  If train ${\mathfrak T}$ must pass certain signal locations at known times and there are no other active time constraints the precise form of the optimal strategy is known \cite{alb2, alb3} and successful computational algorithms can be used on level track to calculate the optimal speed profile~\cite{alb9}.  When the set of active clearance times is known for each train in the fleet the individual optimal strategies can be computed separately~\cite{alb9}.  

Problem~\ref{tsp:p2} has been partially solved.  Necessary conditions are known in various forms for the two-train separation problem~\cite{alb2, alb3, alb8} and more recently in \cite{alb9}.  Similar necessary conditions are also known for the three-train separation problem~\cite[Section 16, pp 162\textendash 165]{alb8} and it is clear that the same arguments and conclusions will remain valid for a fleet of more than three trains.  However computation of the optimal solution\textemdash an optimal speed profile for each train in the fleet\textemdash depends on knowing the locations of the {\em active} constraints.  Safe separation requires $t_i(x_{j+1}) \leq t_{i+1}(x_{j-1})$ for each $i=1,\ldots m-1$ and each $j=1,\ldots n-1$.  This means $(m-1)(n-1)$ inequality constraints and $2^{(m-1)(n-1)}$ possible combinations of active equality constraints.  For each combination it is necessary to compute the optimal times.  There is currently no known way to discard a given combination without checking for safe separation, computing the optimal times and calculating the cost.

\subsection{The main theoretical result}
\label{s:mtr}

In this paper we formulate and solve a train-separation problem for $m$ trains and $n+1$ signals where we aim to minimize total cost subject to a complete set of active clearance time equality constraints that are designed to ensure safe separation and compress the total line-occupancy timespan.  We compress the total line-occupancy timespan by assuming that ${\mathfrak T}_{i+1}$ enters segment ${\mathfrak S}_j = [x_{j-1}, x_{j+1}]$ at the same time as ${\mathfrak T}_i$ leaves.  Thus we assume that $t_{i,j+1} = t_i(x_{j+1}) = t_{i+1}(x_{j-1}) = t_{i+1,j-1}$.  Therefore
\begin{equation}
\label{activect}
t_{i,\, j+1} = t_{i+1,\, j-1} = t_{i+2,\, j-3} = \cdots = t_{i + \lfloor (j+1)/2 \rfloor,\, j + 1 - 2\lfloor (j+1)/2 \rfloor}
\end{equation}
for each $i = 1, \ldots, m-1$ and each $j = 1,\ldots,n$ where $\lfloor (j+1)/2 \rfloor$ denotes the integer part of $(j+1)/2$.  Thus we have $t_{1,2} = t_{2,0}$, $t_{1,3} = t_{2,1}$, $t_{1,4} = t_{2,2} = t_{3,0}$, $t_{1,5} = t_{2,3} = t_{3,1}$, $t_{1,6} = t_{2,4} = t_{3,2} = t_{4,0}$ and so on. 

We will solve the following so-called {\em alternative} train-separation problem.

\begin{problem}
\label{mtr:p3}
Find necessary conditions for a vector $\bfh_0$ of optimal prescribed segment clearance times that minimizes the total tractive energy $J(\bfh)$ required by the fleet over all $\bfh = [h_{i,j}] \in {\mathcal H}_{\mbox{\scriptsize f}}$, where ${\mathcal H}_{\mbox{\scriptsize f}}$ denotes the collection of all feasible sets of prescribed section clearance times with active constraints $t_i(x_j) = h_{i,j}$ and with $h_{i+1,\, j-1} = t_{i+1,\, j-1} = t_{i,\, j+1} = h_{i,\, j+1}$ for each $i=1,\ldots,m-1$ and $j=1,\ldots,n-1$.  The trains may have different performance functions.  $\hfill \Box$
\end{problem}

In practice there are two important modifications that we can make to Problem~\ref{mtr:p3}.
\begin{itemize}
\item We can allow an additional {\em buffer-time} delay $\delta_i \geq 0$ between ${\mathfrak T}_i$ and ${\mathfrak T}_{i+1}$ so that $ t_{i+1,\, j-1} = t_{i,\, j+1} + \delta_i$ for each $i=1,\ldots,m-1$.  Thus we have $t_{2,0} =  t_{1,2} + \delta_1$, $t_{2,1} =  t_{1,3} + \delta_1$, $t_{3,0} =  t_{2,2} + \delta_2 = t_{1,4} + (\delta_1 + \delta_2)$, $t_{3,1} =  t_{2,3} + \delta_2 = t_{1,5} + (\delta_1 + \delta_2)$, $t_{4,0} = t_{3,2} + \delta_3 = t_{2,4} + (\delta_2 + \delta_3) = t_{1,6} + (\delta_1 + \delta_2 + \delta_3)$ and so on.  The optimal speed profiles are not changed by this additional delay.  On a typical train graph this simply means that the planned graph for ${\mathfrak T}_{i+1}$ is shifted to the right by $\delta_1 + \cdots + \delta_i$ time units.  For instance in Figure~\ref{fig1} the train graph for ${\mathfrak T}_2$ is simply moved to the right by $\delta_1$ time units. 
\item We can assume that active constraints are only required for successive trains at certain key locations.  For instance, in the problem depicted in Figure~\ref{fig1}, we could set signal times ${\mathcal H}_1 = [h_{1,j}]_{j=0}^8 = [ 0, h_1, h_2, \ldots, h_7, T ]$ for ${\mathfrak T}_1$ and ${\mathcal H}_2 = [h_{2,j}]_{j=0}^8 = [ \Delta T, h_7, h_8, h_4, h_9, h_{10}, h_{11}, h_{12}, \Delta T + T]$ for ${\mathfrak T}_2$.  The only active intermediate time constraint is $t_{2,3} = h_{2,3} = h_4 = h_{1,5} = t_{1,5}$.  We now optimize over all $\bfh = [h_1,\ldots,h_{12}]$ but must check retrospectively that $h_{1,j+1} \leq h_{2,j-1}$ for all $j = 1,\ldots,7$.
\end{itemize} 

We solve Problem~\ref{mtr:p3} using an extended form of a little-known formula 
\begin{equation}
\label{djdt}
\frac{dJ}{dT} = - \psi(V)
\end{equation}
that expresses the rate of change of cost with respect to journey time for an optimal strategy with no intermediate time constraints as a function of the optimal driving speed.  An embryonic form of the extended formula \cite[Section IV A]{alb9} for the two-train separation problem, shows that the partial rate of change of journey cost with respect to the prescribed signal times is a known function of the preceding and succeeding optimal section driving speeds.  Direct derivations of (\ref{djdt}) have been given for level track \cite{alb6} and more generally for non-level track \cite{how10}.  We note that (\ref{djdt}) was known much earlier in the Russian literature using a high-level derivation based on a deep understanding of the Pontryagin Principle.  See, for instance, a brief reference to the formula in \cite[Section 3, p 922, following equation (12)]{liu1}.

\subsection{Application of the theoretical results}
\label{s:atr}

The theoretical results are applied to a case study of the regular $15$-minute daily passenger shuttle service from Glasgow to Edinburgh in Scotland.

\subsection{The role of the train performance functions}
\label{s:rtpf}

For each train ${\mathfrak T}_i$ the precise details of the optimal strategy\textemdash the optimal driving speeds, the tractive energy consumption, and the optimal switching locations\textemdash depend only on the stopping pattern, the prescribed signal times and the performance functions $(H_i(v), K_i(v), r_i(v))$.  Hence the optimal strategy for ${\mathfrak T}_i$ can be computed separately with no knowledge of the other optimal driving strategies.  Consequently we can use individual performance functions without changing the overall structure of the solution.  This is confirmed by our theoretical results and illustrated in Case Study~\ref{cs5}.

\subsection{Terminology}
\label{s:term}

It is convenient to summarize our terminology. 

\begin{itemize}
  \item In the case studies we assume each train ${\mathfrak T}$ has bounded acceleration with $-K(v) \leq u \leq H(v)$ where $H(v) = \max\{ P_0, P_1/v\}$ and $K(v) = \max\{Q_0, Q_1/v\}$ for $v > 0$ where $P_0, P_1, Q_0, Q_1 \in {\mathbb R}_+$ are positive constants.  The resistive acceleration is given by $r(v) = r_0 + r_1 v + r_2 v^2$ for $v \geq 0$ where $r_0, r_1, r_2 \in {\mathbb R}_+$ are positive constants.  We make repeated use of the functions $\varphi(v) = vr(v)$ and $\psi(v) = v^2r^{\, \prime}(v)$.
  \item  The distance and time differentials $dx_a(v) = v dt_a(v)$, $dx_c(v) = vdt_c(v)$ and $dx_b(v) = vdt_b(v)$ for phases of maximum acceleration, coast and maximum brake are defined as functions of the speed $v$ by (\ref{dxav}), (\ref{dxcv}) and (\ref{dxbv}) respectively.  
  \item The known speed functions $U_b(v)$, $U_s(v,w)$ and $U_s^{\dag}(u,v,w)$ which depend only on the resistance function for the train are defined by (\ref{ub}), (\ref{us}) and (\ref{usd}) respectively.  
  \item For our theoretical derivations it is convenient to assume there is no {\em buffer time} in the form of an additional delay between ${\mathfrak T}_i$ and ${\mathfrak T}_{i+1}$.  That is, in general, we assume $t_{i+1,\, j-1} = t_{i,\, j+1}$ for each $i=1,\ldots,m-1$ and $j=1,\ldots,n-1$.  If a buffer time $\delta_i > 0$ is subsequently imposed between successive trains ${\mathfrak T}_i$ and ${\mathfrak T}_{i+1}$ with $t_{i+1,\, j-1} = t_{i,\, j+1} + \delta_i$ the optimal speed profiles do not change and one simply moves the optimal train graph for ${\mathfrak T}_{i+1}$ an additional $\Delta_i = \delta_1+ \cdots + \delta_i$ units to the right.  Buffer times are essential in practice.   
  \item The term {\em clearance times} for sections $[x_{j-1},x_j]$ or segments ${\mathfrak S}_j = [x_{j-1},x_{j+1}]$ is a collective term for both the {\em earliest allowed entry times} and the {\em latest allowed exit times}.   The times $t_{i,\, j} = t_i(x_j)$ are also known as {\em signal times}.
  \item An inequality constraint, either $t_{i,\, j} \geq h_{i,\, j}$ or $t_{i,\, j} \leq h_{i,\, j}$, is described as {\em active} if $t_{i,\, j} = h_{i,\, j}$ and {\em inactive} otherwise.
  \item A {\em strategy of optimal type} satisfies the necessary conditions for optimality but may not satisfy the imposed time constraints.
  \item A {\em feasible} strategy satisfies the imposed time constraints.
  \item For $s \in {\mathbb R}$ with $s > 0$ we write $\lfloor s \rfloor$ to denote the integer part of $s$.
  \item For vectors $\bfs = [s_j]_{j=1}^p, \bft = [t_j] _{j=1}^p \in {\mathbb R}^p$ the Hadamard product is the vector $\bfr = \bfs \circ \bft$ defined by $\bfr = [r_j]_{j=1}^p$ where $r_j = s_j t_j$ for each $j=1,\ldots,p$.
\end{itemize}

\section{Structure of the paper}
\label{s:sp}

In Section~\ref{s:pw} we review the relevant literature.  In Section~\ref{s:con} we outline the contribution of this paper to the theory and practice of optimal train control.  We propose a systematic scheme for definition of the clearance-time constraints in Section~\ref{s:ctitc}.  The constraints are displayed in a convenient table format.  The next two sections contain the main theoretical results.  In Section~\ref{s:opcss} we solve the optimal scheduling problem for a fleet of trains using constant-speed strategies subject to safe separation of successive trains enforced by active equality constraints.  The problem is formulated as an unconstrained convex optimization.  We find analytic formul{\ae} for the cost gradient with respect to the prescribed clearance times and for the associated Hessian  matrix.  An analogous optimal scheduling  problem for realistic strategies is formulated and solved in Section~\ref{s:oprods}.  Once again the problem is formulated as an unconstrained convex optimization and we establish an analogous analytic formula for the cost gradient.  A general procedure for calculation of realistic optimal strategies for journeys with intermediate time constraints is described in Section~\ref{s:cositc}.  In Section~\ref{s:gsa} we outline a proposed algorithm to find an optimal schedule for a fleet of trains subject to a typical set of safe separation conditions enforced by active equality constraints.  Section~\ref{s:cs} contains a suite of case studies where the theoretical techniques are used to find an optimal schedule for a busy intercity passenger shuttle service subject to safe-separation constraints for a set of nominal sections defined by the key locations.  In Section~\ref{s:tssje} we formulate a rudimentary model for the stochastic evolution of a typical train journey where the driver follows advice provided by an on-board driver advisory system. These ideas are then applied to the optimal schedule developed in the earlier case studies.  Selected background mathematical material is explained and extended in the Appendix.  

\section{Previous work}
\label{s:pw}

There are several widely cited papers which propose analytic solutions to the classic single train control problem.  When continuous control is allowed we refer to \cite{alb1, alb4, alb5, how9} and to earlier papers, \cite{how7, khm1, liu1}.  These authors all use the commonly accepted model of a point-mass train and find solutions by applying classical methods of constrained optimization.   See also \cite{bar2}.  The strategies have been implemented in real time on very fast trains where updated strategies are routinely calculated on-board in a matter of a few seconds~\cite[Section 9, pp 534\textendash 535]{alb5}.  The point-mass model does not consider in-train forces which may be significant in heavy-haul trains \cite{zhu1}.  Other authors use pseudo-spectral methods \cite{gov1, wanp2} but these methods are too slow for on-board use \cite{sch2}.  The train control problem has also been solved for the discrete control systems used on diesel-electric locomotives.  See \cite{che1, how2, pud1} for solutions on level track and subsequent papers \cite{how3, how5, how6} for solutions on tracks with nonzero gradients. 

There are numerous papers that use standard methods of operations research\textemdash mathematical programming, job-shop scheduling, graph theory and various heuristic search procedures\textemdash to find efficient train schedules on complex rail networks.  The work is often broadly-based with objectives that include improved efficiencies obtained by adjusting allowed journey times~\cite{su2}, and best ordering of scheduled services \cite {bur1, bur2, cap1, lius1} and by coordinating arrival and departure times to capture energy from regenerative braking \cite{lilo1, lilo2}.   Other authors considered selecting the best meeting locations~\cite{hig1, wany2}, and providing improved service to customers~\cite{wanp1, yans}.  The development of efficient timetables on complex networks has also spawned work on conflict detection and resolution~\cite{cor1, dar1} and schedule recovery from disruption~\cite{wanp2}.

See \cite{sch1,yanx3, yin1} for recent broadly-based and detailed reviews. 

A recent paper~\cite{wanp3} considers multi-train trajectory optimization for the development of energy-efficient schedules.  The main focus is on single track corridors where conflicts between opposing trains must be resolved and on double track corridors where the focus is on the development of efficient timetables for successive trains and\textemdash to a lesser extent\textemdash on provision for occasional overtaking.  The first step is to replace the existing arrival and departure times with time-window constraints in order to relax the timetable.  Energy-efficient speed profiles are then used either individually (if there are no conflicts) or collectively (if conflicting train paths are involved) to find optimal arrival and departure times within the relaxed time windows. The two problems are reformulated as a multiple-phase optimal control problem and solved by a pseudo-spectral method.  There is no attempt to develop a cohesive theory of optimal scheduling.

\section{Contribution}
\label{s:con}

Our focus is on developing a theory of optimal scheduling and using it to find analytic solutions to unsolved timetable problems on dedicated double track corridors where all trains have similar performance functions.  This includes busy metropolitan transport systems in many major cities, regular intercity passenger shuttle services between neighbouring cities and dedicated heavy-haul freight corridors.  The main theoretical contributions are complementary solutions to the {\em alternative} train-separation problem where the trains are separated by active clearance-time equality constraints  using, firstly, constant-speed strategies and, secondly, realistic strategies.  We show that this problem can be reformulated as an unconstrained convex optimization and that an analogous cost gradient formula defines the optimal solution in each instance. We apply our results to a set of case studies for the daily shuttle service between Glasgow and Edinburgh in Scotland.   Our aim is to minimize total energy consumption subject to active intermediate clearance-time equality constraints that ensure safe separation of successive trains and preserve overall journey times.  For the constant-speed strategies we show that the optimal schedule can be found using a single multi-dimensional Newton iteration.  For the realistic strategies we use the method of steepest descent to find the optimal schedules with a sequence of Newton iterations required at each step to determine the individual optimal strategies.  We also propose a theoretical model to quantify the normal stochastic variation in section running times and hence find appropriate additional buffer times for the separation constraints.  The effectiveness of these buffer times is tested in our final case study.

\section{The time constraints for individual trains in the fleet}
\label{s:ctitc}

We will use a single vector  $\bfh = [ h_k ] \in {\mathbb R}^{n+m-2}$ of unknown prescribed times to define the separation constraints.  The prescribed intermediate times for ${\mathfrak T}_1$ are all unknown with $t_{1,j} = h_j$ for each $j=1,\ldots,n-1$.  For $1 \leq i \leq m-1$ the prescribed times for train ${\mathfrak T}_{i+1}$ are defined by the equality constraint $t_{i+1,j-1} = t_{i,j+1}$ for each $j=1,\ldots,n-1$.  The time $t_{i+1,n-1}$ at $x_{n-1}$ is not constrained because ${\mathfrak T}_i$ has already completed the journey.  Thus we set $t_{i+1,n-1} = h_{n+i-1}$ where $h_{n+i-1}$ is unknown.   For each additional train there is one additional unknown.  The easiest way to visualize the scheduled times for the entire fleet is to construct a table of constraints and unknowns. 

\begin{table}[htb]
\begin{center}
\begin{tabular}{|c|c|c|c|c|c|c|c|c|c|c|c|c|} \hline
${\mathfrak T}_i$ & $t_{i,0}$ & $t_{i,1}$ & $t_{i,2}$ & $t_{i,3}$ & $t_{i,4}$ & $\cdot$ & $t_{i,2k-4}$ & $t_{i,2k-3}$ & $t_{i,2k-2}$ & $t_{i,2k-1}$ & $ t_{i,2k}$ \\ \hline
${\mathfrak T}_1$ & $0$ & $h_1$ & $h_2$ & $h_3$ & $h_4$ & $\cdot$ & $h_{2k-4}$ & $h_{2k-3}$ & $h_{2k-2}$ & $h_{2k-1}$ & T \\ \hline
${\mathfrak T}_2$ & $h_2$ & $h_3$ & $h_4$ & $h_5$ & $h_6$ & $\cdot$ & $h_{2k-2}$ & $h_{2k-1}$ & $T$ & $h_{2k}$ & $h_2+T$ \\ \hline
${\mathfrak T}_3$ & $h_4$ & $h_5$ & $h_6$ & $h_7$ & $h_8$ & $\cdot$ & $T$ & $h_{2k}$ & $h_2+T$ & $h_{2k+1}$ & $h_4+T$ \\ \hline
$\vdots$ & $\vdots$ & $\vdots$ & $\vdots$ & $\vdots$ & $\vdots$ & $\vdots$ & $\vdots$ & $\vdots$ & $\vdots$ & $\vdots$ & $\vdots$ \\ \hline
${\mathfrak T}_{k-1}$ & $h_{2k-4}$ & $h_{2k-3}$ & $h_{2k-2}$ & $h_{2k-1}$ & $T$ & $\cdot$ & $h_{2k-8}+T$ & $h_{3k-3}$ & $h_{2k-6}+T$ & $h_{3k-2}$ & $h_{2k-4}+T$ \\ \hline
${\mathfrak T}_k$ & $h_{2k-2}$ & $h_{2k-1}$ & $T$ & $h_{2k}$ & $h_2+T$ & $\cdot$ & $h_{2k-6}+T$ & $h_{3k-2}$ & $h_{2k-4}+T$ & $h_{3k-1}$ & $h_{2k-2}+T$ \\ \hline \hline
${\mathfrak T}_{k+1}$ & $T$ & $h_{2k}$ & $h_2+T$ & $h_{2k+1}$ & $h_4+T$ & $\cdot$ & $h_{2k-4}+T$ & $h_{3k-1}$ & $h_{2k-2}+T$ & $h_{3k}$ & $h_{2k}+T$ \\ \hline
${\mathfrak T}_{k+2}$ & $h_2+T$ & $h_{2k+1}$ & $h_4+T$ & $h_{2k+2}$ & $h_6+T$ & $\cdot$ & $h_{2k-2}+T$ & $h_{3k}$ & $h_{2k}+T$ & $h_{3k+1}$ & $h_{2k+2}+T$ \\ \hline
$\vdots$ & $\vdots$ & $\vdots$ & $\vdots$ & $\vdots$ & $\vdots$ & $\vdots$ & $\vdots$ & $\vdots$ & $\vdots$ & $\vdots$ & $\vdots$ \\ 
\end{tabular}
\end{center}
\vspace{0.4cm}
\caption{\boldmath Table of constraints and unknown times when $n = 2k \in {\mathbb N}$.}
\label{tab1}
\end{table}

\begin{table}[htb]
\begin{center}
\begin{tabular}{|c|c|c|c|c|c|c|c|c|c|c|c|c|} \hline
${\mathfrak T}_i$ & $t_{i,0}$ & $t_{i,1}$ & $t_{i,2}$ & $t_{i,3}$ & $t_{i,4}$ & $\cdot$ & $t_{i,2k-3}$ & $t_{i,2k-2}$ & $t_{i,2k-1}$ & $t_{i,2k}$ & $ t_{i,2k+1}$ \\ \hline
${\mathfrak T}_1$ & $0$ & $h_1$ & $h_2$ & $h_3$ & $h_4$ & $\cdot$ & $h_{2k-3}$ & $h_{2k-2}$ & $h_{2k-1}$ & $h_{2k}$ & T \\ \hline
${\mathfrak T}_2$ & $h_2$ & $h_3$ & $h_4$ & $h_5$ & $h_6$ & $\cdot$ & $h_{2k-1}$ & $h_{2k}$ & $T$ & $h_{2k+1}$ & $h_2+T$ \\ \hline
${\mathfrak T}_3$ & $h_4$ & $h_5$ & $h_6$ & $h_7$ & $h_8$ & $\cdot$ & $T$ & $h_{2k+1}$ & $h_2+T$ & $h_{2k+2}$ & $h_4+T$ \\ \hline
$\vdots$ & $\vdots$ & $\vdots$ & $\vdots$ & $\vdots$ & $\vdots$ & $\vdots$ & $\vdots$ & $\vdots$ & $\vdots$ & $\vdots$ & $\vdots$ \\ \hline
${\mathfrak T}_{k-1}$ & $h_{2k-4}$ & $h_{2k-3}$ & $h_{2k-2}$ & $h_{2k-1}$ & $h_{2k}$ & $\cdot$ & $h_{2k-8}+T$ & $h_{3k-3}$ & $h_{2k-6}+T$ & $h_{3k-2}$ & $h_{2k-4}+T$ \\ \hline
${\mathfrak T}_k$ & $h_{2k-2}$ & $h_{2k-1}$ & $h_{2k}$ & $T$ & $h_{2k+1}$ & $\cdot$ & $h_{2k-6}+T$ & $h_{3k-2}$ & $h_{2k-4}+T$ & $h_{3k-1}$ & $h_{2k-2}+T$ \\ \hline \hline
${\mathfrak T}_{k+1}$ & $h_{2k}$ & $T$ & $h_{2k+1}$ & $h_2+T$ & $h_{2k+2}$ & $\cdot$ & $h_{2k-4}+T$ & $h_{3k-2}$ & $h_{2k-2}+T$ & $h_{3k}$ & $h_{2k}+T$ \\ \hline
${\mathfrak T}_{k+2}$ & $h_{2k+1}$ & $h_2+T$ & $h_{2k+2}$ & $h_4+T $ & $h_{2k+3}$ & $\cdot$ & $h_{2k-2}+T$ & $h_{3k-2}$ & $h_{2k}+T$ & $h_{3k+1}$ & $h_{2k+2}+T$ \\ \hline
$\vdots$ & $\vdots$ & $\vdots$ & $\vdots$ & $\vdots$ & $\vdots$ & $\vdots$ & $\vdots$ & $\vdots$ & $\vdots$ & $\vdots$ & $\vdots$ \\ 
\end{tabular}
\end{center}
\vspace{0.4cm}
\caption{\boldmath Table of constraints and unknown times when $n = 2k+1 \in {\mathbb N}$.}
\label{tab2}
\end{table}

Although the patterns in Tables~\ref{tab1} and \ref{tab2} are ultimately different they are essentially the same above the double line, where $i \leq \lfloor n/2 \rfloor$.  Thus, a general algebraic description is much easier when $m \leq \lfloor n/2 \rfloor$.

\section{The optimization problem with constant-speed strategies}
\label{s:opcss}

For each feasible $\bfh$ assume that ${\mathfrak T}_i$ travels at constant speed
$$
W_{i,j}(\bfh) = (x_j - x_{j-1})/[t_{i,j}(\bfh) - t_{i,j-1}(\bfh)]
$$
on section $(x_{j-1},x_j)$ for each $i = 1,\ldots,m$ and $j=1,\ldots,n$ subject to the constraints listed in Table~\ref{tab1} or Table~\ref{tab2}.  Let $r_i(v)$ be the resistive acceleration for ${\mathfrak T}_i$ and let $\varphi_i(v) = vr_i(v)$ and $\psi_i(v) = v^2r_i^{\, \prime}(v)$ be the associated functions.  Problem~\ref{mtr:p3} can now be posed as an unconstrained optimization.

\begin{problem}
\label{css:p4}
Let
\begin{equation}
L(\bfh)  = \sum_{i=1}^m \sum_{j=1}^n r_i[W_{i,j}(\bfh)] (x_j - x_{j-1})
\end{equation}
denote the total cost for the fleet to complete the journey.   Find $\bfh = [h_k]_{k=1}^{n+m-2}$ to minimize $L$ subject to the time constraints defined in Table~\ref{tab1} or Table~\ref{tab2}.  In the special case where $m \leq \lfloor n/2 \rfloor$ the time constraints for ${\mathfrak T}_i$ for all $i=1,\ldots,m$ are given by $t_{i,\, j} = h_{2i-2+j}$ for each $j=0,\ldots,n-2i+1$; $t_{i,n-2i+2} = T$; $t_{i,n-2i+2\ell+1} = h_{n+\ell-1}$ for each $\ell = 1,\ldots,i-1$; and $t_{i,n-2i+2\ell+2} = h_{2\ell} + T$ for each $\ell = 1,\ldots,i-1$.  Note that the trains are not assumed to be identical.  $\hfill \Box$
\end{problem} 

To minimize $L$ we must solve the equation
\begin{equation}
\label{delk}
\bnab L(\bfh) = \bfzero \iff \left[ \partial L/\partial h_k \right]_{k=1}^{n+m-2} = \bfzero.
\end{equation}   
The total cost for train ${\mathfrak T}_i$ is $L_i(\bfh) = \sum_{j=1}^n r_i[W_{i,j}(\bfh)](x_j - x_{j-1})$.  For each fixed $k \in \{1,\ldots,n\}$ it is useful to begin by identifying the terms which depend on the variable $h_k$.  Let ${\mathcal F}_{i,k} = \{ j \mid t_{i,j} = h_k + c_{i,j}$ for some constant $c_{i,j} \geq 0$ \}.  It is easy to see from Table~\ref{tab1} and Table~\ref{tab2} that $j \in {\mathcal F}_{i,k} \Rightarrow j-1, j+1 \notin {\mathcal F}_{i,k}$.  Therefore
\begin{eqnarray}
\label{dkidhk}
\frac{\partial L_i}{\partial h_k} & = & \sum_{j \in {\mathcal F}_{i,k}} \left[ r_i^{\, \prime}(W_{i,j+1})(x_{j+1}- x_j) \cdot (x_{j+1} - x_j)/(t_{i,j+1} - t_{i,j})^2 \right. \nonumber \\
& & \hspace{4cm} \left.- r_i^{\, \prime}(W_{i,j})(x_j- x_{j-1}) \cdot (x_j - x_{j-1})/(t_{i,j} - t_{i,j-1})^2 \right] \nonumber \\
& = & \sum_{j \in {\mathcal F}_{i,k}} \left[ \psi_i(W_{i,j+1}) - \psi_i(W_{i,j}) \right].
\end{eqnarray}
In order to apply a Newton iteration to solve (\ref{delk}) it is necessary to calculate the Hessian matrix of second derivatives. In this regard we find
\begin{equation}
\label{d2Lidk2}
\frac{\partial^2 L_i}{\partial h_k^2} = \sum_{j \in {\mathcal F}_{i,k}} \left[ \rule{0cm}{0.4cm} W_{i,j+1} \psi_i^{\, \prime}(W_{i,j+1})/(t_{i,j+1}- t_{i,j}) + W_{i,j}\psi_i^{\, \prime}(W_{i,j})/(t_{i,j}- t_{i,j-1}) \right]
\end{equation}
and since $j \in {\mathcal F}_{i,k} \Rightarrow j \notin {\mathcal F}_{i,\ell}$ when $\ell \neq k$ we  have
\begin{equation}
\label{d2Lidelldk}
\frac{\partial^2 L_i}{\partial h_{\ell} \partial h_k} = \hspace{2mm} - \hspace{-3.5mm}  \sum_{j \in {\mathcal F}_{i,k}\, \mbox{\scriptsize \&}\, j+1 \in {\mathcal F}_{i,\ell}} \hspace{-3.5mm} W_{i,j+1} \psi_i^{\, \prime}(W_{i,j+1})/(t_{i,j+1}- t_{i,j})\hspace{2mm} - \hspace{-5mm} \sum_{j \in {\mathcal F}_{i,k}\, \mbox{\scriptsize \&}\, j-1 \in {\mathcal F}_{i,\ell}} \hspace{-5mm}W_{i,j}\psi_i^{\, \prime}(W_{i,j})/(t_{i,j}- t_{i,j-1})
\end{equation}
when $\ell \neq k$.  We will show that each $L_i(\bfh)$ is convex.  We may assume without loss of generality that ${\mathfrak T}_i$ satisfies a set of intermediate time constraints $t_i(x_j) = g_{i,j}$ for each $j=1,\ldots,n$.  It follows from (\ref{d2Lidk2}) and (\ref{d2Lidelldk}) that the Hessian matrix $P_i = [p_{i,k,\ell}] = [\partial^2 L_i /(\partial g_{i,k} \partial g_{i,\ell})] \in {\mathbb R}^{n \times n}$ is given by
$$
p_{i,k,\ell} = \left\{ \begin{array}{ll}
- \beta_{i,k} & \mbox{when}\ \ell = k-1 \\
\beta_{i,k+1} + \beta_{i,k} & \mbox{when}\ \ell = k \\
- \beta_{i,k+1} & \mbox{when}\ \ell = k+1 \\
0 & \mbox{otherwise} \end{array} \right.
$$
where $W_{i,k} = (x_k - x_{k-1})/(g_{i,k} - g_{i,k-1})$ and $\beta_{i,k} = W_{i,k} \psi_i^{\, \prime}(W_{i,k})/( g_{i,k} - g_{i,k-1}) > 0$.  Let $P_{i,r} = [p_{i,k,\ell}]_{k,\ell = 1}^r$ denote the principal minor for each $r=1,\ldots,n$.  Now it can be shown that
$$
\det P_{i,r} = \mbox{$\sum_{\sbsigma \in {\mathcal P_r}} \left[ \prod_{j=1}^r \beta_{i,\sigma(j)} \right]$} > 0
$$
where ${\mathcal P}_r$ is the set of all strictly increasing sequences $\bsigma = \{ \sigma(j)\}_{j=1}^r$ with $1 \leq \sigma(1) < \cdots < \sigma(r) \leq n$ for each $r=1,\ldots,n$.  It follows from Sylvester's criterion \cite{gil1} that $P_i$ is positive definite.  Hence $L_i(\bfh)$ is convex  for each $i=1,\ldots,m$.  It follows that $L(\bfh)$ is also convex.   Thus the solution to (\ref{delk}) is unique.  It is conceptually straightforward but algebraically complicated to write down general formul{\ae} for the Newton iteration.  We prefer to elaborate those details in particular case studies where (\ref{delk}) will solved directly using a Newton iteration.
  
\section{The optimization problem with realistic optimal driving strategies}
\label{s:oprods}

For each feasible $\bfh = [h_k]_{k=1}^{n+m-2}$ assume that all trains use realistic strategies subject to the time constraints listed in Table~\ref{tab1} or Table~\ref{tab2}.  For each $i=1,\ldots,m$ let $\bfV_i(\bfh) = \{ V_{i,j}(\bfh)\}_{j=1}^n$ denote the vector of optimal driving speeds for ${\mathfrak T}_i$ and let $J_i(\bfh) = J_i[\bfV_i(\bfh)]$ denote the total cost for ${\mathfrak T}_i$.  The constrained optimization posed in Problem~\ref{mtr:p3} is now an unconstrained optimization problem in $\bfh$.

\begin{problem}
\label{rs:p5}
Let
\begin{equation}
J(\bfh)  = \sum_{i=1}^m J_i[\bfV_i(\bfh)]
\end{equation}
denote the total cost for the fleet to complete the journey.   Find $\bfh = [h_k]_{k=1}^{n+m-2}$ to minimize $J$ subject to the time constraints defined in Table~\ref{tab1} or Table~\ref{tab2}.  If $m \leq \lfloor n/2 \rfloor$ the time constraints take the algebraic form described in Problem~\ref{css:p4}.  The trains are not assumed to be identical. $\hfill \Box$
\end{problem}

To minimize $J$ we must solve the equations
\begin{equation}
\label{delj}
\bnab J(\bfh) = \bfzero \iff \left[ \partial J/\partial h_k \right]_{k=1}^{n+m-2} = \bfzero.
\end{equation}   
For train ${\mathfrak T}_i$ we will argue that
\begin{equation}
\label{djidhk}
\frac{\partial J_i}{\partial h_k} = \sum_{j \in {\mathcal F}_{i,k}} \left[ \psi_i(V_{i,j+1}) - \psi_i(V_{i,j}) \right]
\end{equation}
where ${\mathcal F}_{i,k} = \{ j \mid t_{i,j} = h_k + c_{i,j}$ for some constant $c_{i,j} \geq 0$ \}.  In the following argument it is convenient to drop the subscript $i$ and argue for a generic train ${\mathfrak T}$.  We may assume without loss of generality that $v_0 = v(x_0) = 0$ and $v_n = v(x_n) = 0$ with $v(x) > 0$ for all $x \in (0,X)$ and that $t(x_j) = g_j$ for each $j=0,\ldots,n$ where $g_0 < g_1 < \cdots < g_{n-1} < g_n$ and where $g_0$ and $g_n = g_0 +T$ are fixed.  The optimal strategy, described in Section~\ref{s:sot}, is a sequence of speedhold phases on segments $[a_j,b_j] \subset (x_{j-1},x_j)$ separated by transition phases of either maximum acceleration or coast when passing through the signal locations $\{x_j\}_{j=1}^{n-1}$.  Because we do not know whether a particular transition phase will be maximum acceleration or coast we will use an elementary notational subterfuge.  For the transition phase at $x_j$ the acceleration will be denoted by $A_j(v)$.   Thus $v = v(x)$ is defined by $v v^{\, \prime} = A_j(v) - r_i(v)$ for $x \in (b_j, a_{j+1})$ for each $j=1,\ldots,n-1$ where $A_j(v) = H(v)$ for a phase of maximum acceleration and $A_j(v) = 0$ for a phase of coast.  The time taken on $(x_0,x_1)$ is
\begin{equation}
\label{time1}
g_1 - g_0 = \int_0^{V_1} \frac{dv}{H(v) - r(v)} + \int_{V_1}^{U_1} \frac{dv}{A_1(v) - r(v)} + \frac{\xi_1}{V_1},
\end{equation}
where
\begin{equation}
\label{delx1}
\xi_1 = x_1 - x_0 - \int_0^{V_1} \frac{vdv}{H(v) - r(v)} - \int_{V_1}^{U_1} \frac{vdv}{A_1(v) - r(v)}
\end{equation}
is the length of the speedhold segment $[a_1,b_1]$, the time taken on $(x_{j-1},x_j)$ is
\begin{equation}
\label{timej}
g_j - g_{j-1} = \int_{U_{j-1}}^{V_j} \frac{dv}{A_{j-1}(v)-r(v)} + \int_{V_j}^{U_j} \frac{dv}{A_j(v) - r(v)} + \frac{\xi_j}{V_j}
\end{equation}
where
\begin{equation}
\label{delxj}
\xi_j =  x_j - x_{j-1} - \int_{U_{j-1}}^{V_j} \frac{vdv}{A_{j-1}(v)-r(v)} - \int_{V_j}^{U_j} \frac{vdv}{A_j(v) - r(v)}
\end{equation}
is the length of the speedhold segment $[a_j,b_j]$ for each $j = 1,\ldots,n-1$ and the time taken on $(x_{n-1},x_n)$ is
\begin{equation}
\label{timen}
T - g_{n-1} =  \int_{U_{n-1}}^{V_n} \frac{dv}{A_{n-1}(v) - r(v)} + \int_{U_n}^{V_n} \frac{dv}{r(v)} + \int_0^{U_n} \frac{dv}{K(v) + r(v)} + \frac{\xi_n}{V_n}
\end{equation}
where
\begin{equation}
\label{delxn}
\xi_n = x_n - x_{n-1} - \int_{U_{n-1}}^{V_n} \frac{vdv}{A_{n-1}(v) - r(v)} - \int_{U_n}^{V_n} \frac{vdv}{r(v)} - \int_0^{U_n} \frac{vdv}{K(v) + r(v)}
\end{equation}
is the length of the speedhold segment $[a_n,b_n]$.  If we differentiate (\ref{time1}), (\ref{timej}) and (\ref{timen}) with respect to $g_k$ and rearrange the terms we get
\begin{equation}
\label{dv1dgk}
\frac{\xi_1}{V_1^2} \frac{\partial V_1}{\partial g_k} = - \delta_{1,k} + \frac{1 - U_1/V_1}{A_1(U_1) - r(U_1)} \frac{\partial U_1}{\partial g_k},
\end{equation}
\begin{equation}
\label{dvjdgk}
\frac{\xi_j}{V_j^2} \frac{\partial V_j}{\partial g_k} = \delta_{j-1,k} - \delta_{j,k} - \frac{1 - U_{j-1}/V_j}{A_{j-1}(U_{j-1}) - r(U_{j-1})}  \frac{\partial U_{j-1}}{\partial g_k} + \frac{1 - U_j/V_j}{A_j(U_j) - r(U_j)} \frac{\partial U_j}{\partial g_k}
\end{equation}
for each $j=2,\ldots,n-1$ and
\begin{equation}
\label{dvndgk}
\frac{\xi_n}{V_n^2} \frac{\partial V_n}{\partial g_k} = \delta_{n-1,k} - \frac{1 - U_{n-1}/V_n}{A_{n-1}(U_{n-1}) - r(U_{n-1})} \frac{\partial U_{n-1}}{\partial g_k} - \frac{(1 - U_n/V_n)K(U_n)}{r(U_n)(K(U_n)+r(U_n))} \frac{\partial U_n}{\partial g_k}
\end{equation}
where $\delta_{r,s} = 1$ if $s=r$ and $\delta_{r,s} = 0$ otherwise.  The cost of the strategy is given by
\begin{equation}
\label{jtc}
J = \int_0^{V_1} \frac{v H(v) dv}{H(v)-r(v)} + \sum_{j=1}^n r(V_j) \xi_j + \sum_{j=1}^{n-1} \int_{V_j}^{V_{j+1}} \frac{vA_j(v) dv}{A_j(v)-r(v)}. 
\end{equation}
Differentiation with respect to $g_k$ and simplification gives
\begin{equation}
\label{djtcdgk}
\frac{\partial J}{\partial g_k} = \sum_{j=1}^{n-1} \frac{[r(V_{j+1}) - r(V_j)]U_j}{A_j(U_j) - r(U_j)} \frac{\partial U_j}{\partial g_k} + \frac{r(V_n)U_n K(U_n)}{r(U_n)(K(U_n)+r(U_n))} \frac{\partial U_n}{\partial g_k} + \sum_{j=1}^n r^{\, \prime}(V_j) \xi_j \frac{\partial V_j}{\partial g_k}
\end{equation}
and if we eliminate terms in $\partial V_j/\partial g_k$ using (\ref{dv1dgk}), (\ref{dvjdgk}) and (\ref{dvndgk}) we obtain
\begin{eqnarray}
\label{djtcdgks}
\frac{\partial J}{\partial g_k} & = & \sum_{j=1}^{n-1} \frac{[r(V_{j+1}) - r(V_j)]U_j}{A_j(U_j) - r(U_j)} \frac{\partial U_j}{\partial g_k} + \frac{r(V_n)U_n K(U_n)}{r(U_n)(K(U_n)+r(U_n))} \frac{\partial U_n}{\partial g_k} \nonumber \\
& & +\psi(V_1)\left[ -\delta_{1,k} + \frac{1 - U_1/V_1}{A_1(U_1)-r(U_1)}\frac{\partial U_1}{\partial g_k}\right]  \nonumber \\
& & + \sum_{j=2}^{n-1} \psi(V_j) \left[ \delta_{j-1,k} - \delta_{j,k} - \frac{1 - U_{j-1}/V_j}{A_{j-1}(U_{j-1}) - r(U_{j-1})}  \frac{\partial U_{j-1}}{\partial g_k} + \frac{1 - U_j/V_j}{A_j(U_j) - r(U_j)} \frac{\partial U_j}{\partial g_k} \right] \hspace{1.25cm} \nonumber \\
& & + \psi(V_n) \left[ \delta_{n-1,k} - \frac{1 - U_{n-1}/V_n}{A_{n-1}(U_{n-1}) - r(U_{n-1})} \frac{\partial U_{n-1}}{\partial g_k} - \frac{(1 - U_n/V_n)K(U_n)}{r(U_n)(K(U_n)+r(U_n))} \frac{\partial U_n}{\partial g_k} \right].
\end{eqnarray}
If we collect like terms and use (\ref{ub}) and (\ref{us}) we see that the coefficient of
$$
(K(U_n)/[r(U_n)(K(U_n)+r(U_n))])/(\partial U_n/\partial g_k)
$$
is $\varphi^{\, \prime}(V_n) U_n - \psi(V_n) = 0$ and that the coefficient of
$$
(1 /[A_j(U_j) - r(U_j)])(\partial U_j/\partial g_k)
$$
is $[\varphi^{\, \prime}(V_{j+1}) - \varphi^{\, \prime}(V_j)] U_j + \psi(V_j) - \psi(V_{j+1}) = 0$ for each $j=1,\ldots,n-1$.  Therefore 
$$
\frac{\partial J}{\partial g_k} = - \psi(V_1)\delta_{1,k} + \sum_{j=2}^{n-1} \psi(V_j) [ \delta_{j-1,k} - \delta_{j,k}] + \psi(V_n) \delta_{n-1,k} = \psi(V_{k+1}) - \psi(V_k).
$$

To determine an optimal strategy for a journey with prescribed intermediate times but no intermediate stops it may be necessary to use a progressive sequence of Newton iterations.  We will describe a general procedure for that purpose in the next section.  When we find the optimal strategy we also find the optimal driving speeds $V_{i,j}$.  Thus we can calculate $\psi_i(V_{i,j})$.  If there is no speedhold phase on the final section $(x_{n-1},x_n)$ then we must use the formula $\psi_i(V_{i,n}) = \varphi(V_{\max,i,n})U_{i,n}/(V_{\max,i,n}-U_{i,n})$.  See Appendix~\ref{s:rccitc} for a detailed justification in a special case.  The general argument is similar.

The question still remains as to how we solve (\ref{delj}).   We cannot use a Newton iteration in the way that we could for Problem~\ref{css:p4} because, aside from the need for a separate Newton iteration to find the optimal strategies and the associated optimal driving speeds for each train on the segments between scheduled stops, there is no known analytic expression for $\partial V_{i,j}/ \partial h_{\ell}$ and hence no analytic expression for the Hessian matrix.  However the cost function (\ref{djidhk}) is a sum of terms in essentially the same form as the locally convex functions considered in \cite[Appendices C, D and F]{alb9}.  We can use the same algebraic arguments here to estimate the derivatives $\partial V_{i,j}/ \partial h_{\ell}$ and hence show that each $J_i(\bfh)$ is locally convex.  Thus $J(\bfh) = \sum_{i=1}^m J_i(\bfh)$ is also locally convex.  Consequently we can use the method of steepest descent to find a uniquely defined set of clearance times $\bfh_0$ and a corresponding locally optimal cost $J(\bfh_0)$.

To apply the method of steepest descent we need an initial set of feasible time constraints. We will use the optimal times obtained from the solution to Problem~\ref{css:p4} for that purpose.  Let us suppose then that we have an initial set of times defined by a particular value of $\bfh$.  We calculate optimal realistic strategies for each train in the fleet subject to the given constraints.  This gives us a complete set of optimal driving speeds $\{V_{i,j}\}$ for each $i=1,\ldots,m$ and each $j=1,\ldots,n$.  Thus we can calculate $\bnab J$ and find a new value of $\bfh$ according to the formula $\bfh_{\mbox{\scriptsize \rm new}} = \bfh_{\mbox{\scriptsize old}} - r \bnab J(\bfh_{\mbox{\scriptsize old}})$ for an appropriate value of $r > 0$.  Because $J$ is convex $r$ can be determined using a golden search \cite[Appendix ~F]{alb9}.  We now find optimal realistic strategies for the new value of $\bfh$ and calculate a reduced value of $J$.  The procedure can then be repeated using the new optimal driving speeds to further reduce $J$ until we find an approximate solution to (\ref{delj}).

{\bf A remarkable analogy.} The formul{\ae} (\ref{dkidhk}) and (\ref{djidhk}) are remarkable because they are essentially identical.  Thus the necessary conditions for optimality are basically the same for constant-speed strategies as they are for realistic strategies.  This does not mean the optimal speeds computed for the constant-speed strategies are identical to the optimal driving speeds computed for the realistic strategies.  There are however important structural similarities to the optimal train graphs.  In each case the trains must travel faster on the longer sections and in each case the optimal solutions are a compromise between minimizing the higher speeds and equalising the section traversal times.

\section{Calculating an optimal strategy for a journey with intermediate time constraints but no intermediate stops}
\label{s:cositc}

We wish to calculate a realistic optimal strategy for a typical train ${\mathfrak T}$ subject to a known set of prescribed intermediate times $t(x_j) = t_j$ for each $j=1,\ldots,n$ with $0 = t_0 < t_1 < \cdots < t_{n-1} < t_n = T$ where $T$ is the total journey time.  The optimal strategy is uniquely defined by a sequence of optimal driving speeds $\{V_j\}_{j=1}^n$.  For each set of prescribed times it may be difficult to decide {\em a priori} whether the transition from speed $V_j$ on $[a_j,b_j] \subset (x_{j-1},x_j)$ to speed $V_{j+1}$ on $[a_{j+1},b_{j+1}] \subset (x_j, x_{j+1})$ will require a phase of coast or a phase of maximum acceleration.  In such cases a progressive procedure could be used.  At each stage the relevant equations can be solved by a Newton iteration.  Every optimal journey begins with a phase of maximum acceleration followed by a phase of speedhold at the optimal driving speed $V_1$.  At Stage~1 we choose $V_1 = {V_1}^{(1)}$ to satisfy
$$
\mbox{$\int_0^{V_1}$} dt_a(v) + (1/V_1)\left[x_1 - x_0 - \mbox{$\int_0^{V_1}$} dx_a(v)\right] = t_1
$$
where $t_1$ is the prescribed arrival time at $x_1$.  Suppose now that the speedhold phase at speed $V_1 = {V_1}^{(1)}$ is continued until $x_2$ and that $t_1 + [x_2 - x_1]/V_1 = t_2 + \epsilon$ where $\epsilon > 0$ is a small positive number.  Hence we must go a little faster on $(x_1, x_2)$.  Thus we insert a phase of maximum acceleration to change from the driving speed $V_1 = {V_1}^{(1)}$ on $[a_1,b_1] \subset (x_0, x_1)$ to a slightly higher driving speed $V_2$ on $[a_2,b_2] \subset (x_1,x_2)$.  It is known that $v(x_1) = U_1 = U_s(V_1,V_2)$ always lies between $V_1$ and $V_2$.  Thus the train must begin increasing speed before it reaches $x_1$.  This means the true optimal driving speed $V_1$ must be slightly lower than ${V_1}^{(1)}$.  Thus, at Stage~2 we assume a transition phase of maximum acceleration and choose $V_1= {V_1}^{(2)}$ and $V_2 = {V_2}^{(2)}$ to satisfy
\begin{eqnarray*}
\mbox{$\int_0^{U_1}$} dt_a(v) + (1/V_1) \left[x_1 - x_0 - \mbox{$\int_0^{U_1}$} dx_a(v) \right] & = & t_1 \\
\mbox{$\int_{U_1}^{V_2}$} dt_a(v) + (1/V_2) \left[ x_2 - x_1 - \mbox{$\int_{U_1}^{V_2}$} dx_a(v) \right] & = & t_2 - t_1.
\end{eqnarray*}
This means ${V_2}^{(2)}$ will be slightly larger than ${V_1}^{(2)}$.  Now suppose we continue at speed $V_2 = {V_2}^{(2)}$ through until $x_3$ and find that $t_2 + [x_3 - x_2]/V_2 = t_3 + M$ where $M > 0$ is large.  Thus we need to go much faster on $(x_2,x_3)$. Hence we insert a phase of maximum acceleration at Stage~3 to change from speed $V_2 = {V_2}^{(2)}$ on $[a_2,b_2] \subset (x_1, x_2)$ to a much higher speed $V_3$ on $[a_3,b_3] \subset (x_2,x_3)$.  Therefore $U_2 = U_s(V_2,V_3)$ is also much larger than $V_2 = {V_2}^{(2)}$ and so the train must start going much faster as it approaches $x_2$.  If we wish to solve the system
\begin{eqnarray*}
\mbox{$\int_0^{U_1}$} dt_a(v) + (1/V_1)\left[x_1 - x_0 - \mbox{$\int_0^{U_1}$} dx_a(v) \right] & = & t_1 \\
\mbox{$\int_{U_1}^{U_2}$} dt_a(v) + (1/V_2) \left[ x_2 - x_1 - \mbox{$\int_{U_1}^{U_2}$} dx_a(v) \right] & = & t_2 - t_1 \\
\mbox{$\int_{U_2}^{V_3}$} dt_a(v) +(1/V_3) \left[ x_3 - x_2 - \mbox{$\int_{U_2}^{V_3}$} dx_a(v) \right]  & = & t_3 - t_2
\end{eqnarray*} 
at Stage~3 we must choose $V_2 = {V_2}^{(3)}$ to be substantially less than ${V_2}^{(2)}$.  This in turn means we must choose ${V_1}^{(3)} > {V_1}^{(2)}$.  If $\epsilon$ is sufficiently small and $M$ is sufficiently large these changes will mean that ${V_1}^{(3)} > {V_2}^{(3)}$.  Therefore we will now need to use a phase of coast to change from driving speed $V_1={V_1}^{(3)}$ on $[a_1,b_1] \subset (x_0,x_1)$ to a lower driving speed $V_2={V_2}^{(3)}$ on $[a_2,b_2] \subset (x_2,x_2)$.  As we progress further and further with the calculation the impacts on our earlier decisions become less and less pronounced.  It is most likely, in the case described above, that the revised decision will ultimately be the correct decision.  The new optimal driving speeds will indicate the correct choices for the revised calculation.  The final stage on the interval $(x_{n-1},x_n)$ is different.  There are two possibilities.  If there is a speedhold phase at speed $V_n$ on some interval $[a_n,b_n] \subset (x_{n-1},x_n)$ then the speed at $x_{n-1}$ is given by $U_{n-1} = U_s(V_{n-1},V_n)$ and the strategy takes the form of a transition phase\textemdash either maximum acceleration or coast\textemdash followed by speedhold at the optimal driving speed $V_n$, coast to the optimal braking speed $U_n = U_b(V_n)$ and maximum brake.  If there is insufficient time for a speedhold phase then the transition phase through $x_{n-1}$ must be maximum acceleration to speed $V_{\max,n}$, coast to speed $U_n$ and maximum brake.  In this case the speed at $x_{n-1}$ must be $U_{n-1} = U_{n-1}^{\dag}(U_n,V_{n-1},V_{\max,n})$ where the variables $V_{\max,n}$ and $U_n$ must be chosen to satisfy the distance and time constraints on $(x_{n-1},x_n)$.  See Appendix~\ref{s:itc} for details and an application to a journey with one intermediate time constraint.

\begin{remark}
\label{cositc:r1}
It is often possible to decide {\em a priori} on the correct choices for the transitional phases in which case a single Newton iteration is sufficient to determine the optimal times.  In most cases $W_{i,j} < W_{i,j+1} \Rightarrow V_{i,j} < V_{i,j+1}$ and $W_{i,j} > W_{i,j+1} \Rightarrow V_{i,j} > V_{i,j+1}$.  If the journey for ${\mathfrak T}_i$ consists of a sequence of shorter journeys between scheduled stops then the separate journeys depend only on the prescribed times and can be computed independently. $\hfill \Box$
\end{remark}

\section{A general solution algorithm}
\label{s:gsa}

The following algorithm finds a near optimal schedule for a fleet of trains.

\begin{algorithmic}

\STATE{\bf Step \mbox{\boldmath $1$}.} Define the table of active time constraints.  Choose an error margin $\epsilon > 0$ for the cost.
\STATE{\bf Step \mbox{\boldmath $2$}.} Define $0 = j(i,0)< j(i,1) < \cdots < j(i, p_i) = n$ for each $i=1,\ldots,m$ to denote the stopping pattern $0 = x_{j(i,0)}  < x_{j(i,1)} < \cdots < x_{j(i,p_i)} = x_n$ for ${\mathfrak T}_i$.
\STATE{\bf Step \mbox{\boldmath $3$}.} Use a Newton iteration to find optimal times and speeds for the constant-speed strategies.  Repeat with appropriate speed penalties if the timetable is unsuitable for realistic strategies. 
\STATE{\bf Step \mbox{\boldmath $4$}.} For each $i=1,\ldots,m$ find the vector of optimal driving speeds $\bfV_i$ for ${\mathfrak T}_i$ using realistic strategies of optimal type on each segment of the journey subject to the times generated at Step $3$.  For ${\mathfrak T}_i$ the separate journey segments will be $(x_{j(i,0)}, x_{j(i,1)}), (x_{j(i,1)},\, x_{j(i,2)}),\, \ldots,(x_{j(i, p_i-1)}, x_{j(i,p_i)})$.  If $j(i,r+1) = j(i,r)+1$ the strategy on segment $(x_{j(i,r)}.x_{j(i,r+1)})$ will be a long-haul strategy or a rapid-transit strategy.  If $j(i,r+1) > j(i,r)+1$ the strategy on segment $(x_{j(i,r)}, x_{j(i,r+1)})$ will be an optimal strategy with one or more intermediate time constraints.  Calculate $J(\bfV_1,\ldots,\bfV_m)$.
\STATE{\bf Step \mbox{\boldmath $5$}.}  Calculate the gradient $\bnab J = \partial J/ \partial \bfh$ of the cost with respect to the prescribed intermediate times defined by $\bfh_{\mbox{\scriptsize old}}$ for the current optimal driving speed vectors $\bfV_1,\ldots,\bfV_m$.  Move an appropriate distance $r$ in the direction of $(-1) \bnab J$ to find improved intermediate times $\bfh_{\mbox{\scriptsize new}} = \bfh_{\mbox{\scriptsize old}} - r \bnab J$.  An approximate optimal  value for $r$ can be found using a golden search.  This requires running Step $4$ for each selected value of $r$. 
\STATE{\bf Step \mbox{\boldmath $6$}.} Repeat Step $4$ using the new value $\bfh_{\mbox{\scriptsize new}} = \bfh_{\mbox{\scriptsize old}} - r_{\mbox{\scriptsize opt}} \bnab J(\bfh_{\mbox{\scriptsize old}})$.
\STATE{\bf Step \mbox{\boldmath $7$}.} If $J(\bfh_{\mbox{\scriptsize old}}) - J(\bfh_{\mbox{\scriptsize new}}) >  \epsilon$ go to Step $5$.  Else End.
\end{algorithmic} 

\section{Case studies}
\label{s:cs}

We will apply our results to the premium shuttle service operated by Abellio ScotRail on the Glasgow Queen Street to Edinburgh Waverley via Falkirk line.  The trains in this service are a British Rail Class 385~AT~200 Hitachi electric multiple unit operated as three or four cars driven respectively by three or four 250 kW (335 hp) tandem motors with an IGBT converter/inverter.  The multiple unit can be effectively modelled in each case as a point mass train using the equations of motion (\ref{em:vx}) and (\ref{em:tx}) with constraints on the acceleration given by $H(v) \leq \min\{ P_0, P_1/v \}$ ms$^{-2}$ where $P_0 = 0.84$ ms$^{-2}$ and $P_1 = 9.10$ m$^2$s$^{-3}$ and $K(v) \leq \max\{ -Q_0, -Q_1/v \}$ ms$^{-2}$ where $Q_0 = 1.00$ ms$^{-2}$ and $Q_1 = 10.83$ m$^2$s$^{-3}$ and with resistance $r(v) = r_0 + r_1v + r_2v^2$ ms$^{-2}$ where $r_0 = 0.12$ ms$^{-2}$, $r_1 = 0.0001$ s$^{-1}$ and $r_2 = 0.00004$ m$^{-1}$.  We assume that the track is level.  The key locations and a typical weekday timetable for departures from Glasgow Queen Street at $15$-minute intervals over a $60$-minute period are shown in Table~\ref{tab3}.   The station codes are Glasgow Queen Street (GLQ), Bishopbriggs (BBG), Lenzie (LNZ), Croy (CRO), Falkirk High (FKK), Polmont (PMT), Linlithgow (LIN), Winchburgh Junction (WGJ), Haymarket (HYM) and Edinburgh Waverley (EDB).

\begin{table}[htb]
\begin{center}
\begin{tabular}{|l|c|c|c|c|c|} \hline
Station & Position & ${\mathfrak T}_1$ & ${\mathfrak T}_2$ & ${\mathfrak T}_3$ & ${\mathfrak T}_4$ \\ \hline
Glasgow Queen Street (GLQ) & $x_0 = 00000$ & $t_{1,0} = 0000$ & $t_{2,0} = 0900$ & $t_{3,0} = 1800$ & $t_{4,0} = 2700$ \\ \hline
Bishopbriggs (BBG)& $x_1 = 05140$ & \textemdash & \textemdash & \textemdash & \textemdash \\ \hline
Lenzie (LNZ) & $x_2 = 09980$ & \textemdash & \textemdash & \textemdash & \textemdash \\ \hline
Croy (CRO) & $x_3 = 18350$ & $t_{1,3} = 780$ & \textemdash & $t_{3,3} = 2520$ & \textemdash \\ \hline
Falkirk High (FKK) & $x_4 = 34820$ & $t_{1,4} = 1380$ & $t_{2,4} = 2040$ & $t_{3,4} = 3120$ & $t_{4,4} = 3820$ \\ \hline
Polmont (PMT) & $x_5 = 40250$ & \textemdash & $t_{2,5} = 2340$ & \textemdash & $t_{4,5} = 4120$ \\ \hline
Linlithgow (LIN) & $x_6 = 47590$ & \textemdash & $t_{2,6} = 2700$ & \textemdash & $t_{4,6} = 4780$ \\ \hline
Winchburgh Junction (WGJ) & $x_7 = 55610$ & \textemdash & \textemdash & \textemdash & \textemdash \\ \hline
Haymarket (HYM) & $x_8 = 73010$ & $t_{1,8} = 2760$ & $t_{2,8} = 3840$ & $t_{3,8} = 4560$ & $t_{4,8} = 5460$ \\ \hline
Edinburgh Waverley (EDB) & $x_9 = 75700$ & $t_{1,9} = 3180$ & $t_{2,9} = 4140$ & $t_{3,9} = 4860$ & $t_{4,9} = 5820$ \\ \hline
\end{tabular}
\end{center}
\vspace{0.4cm}
\caption{\boldmath Published departure times for a $60$-minute period on a typical weekday morning.}
\label{tab3}
\end{table}

In practice safe separation is assured using a three-aspect signalling system and a set of designated sections defined by the trackside signals.  The precise details depend on the number and location of the designated signals.  We will illustrate the underlying issues using nominal track sections defined by the key locations shown in Table~\ref{tab3}.  We will assume that in normal operation we require one clear section of track between successive trains at all times.  The mathematical optimization does not depend on the number of sections, the lengths of the sections or on a common journey time.  Nor does it depend on the precise train performance functions although each train must have the capacity to complete the entire journey  in the scheduled time using a realistic strategy.

We begin by designing energy-efficient schedules using (unrealistic) constant-speed strategies for each train on each section.  To ensure realistic outcomes we may need to impose speed penalties on selected sections that recognise the realistic constraints on maximum speed imposed by the train performance functions and\textemdash if the train has to stop at particular stations\textemdash to allow additional time for the train to speed up and slow down.  For the Class 385~AT~200 Hitachi electric multiple-unit trains the upper bound on speed $v = V_{\sup}$ can be calculated by solving the equation $H(v) - r(v) = 0$.  This gives $V_{\sup} \approx 44.6011$ ms$^{-1}$.  For a typical speed $V_t = 30$ ms$^{-1}$ we have
$$
\tau_{a,t} = \int_0^{V_t} dt_a(v) \approx 82.91 \quad \mbox{and} \quad \xi_{a,t} = \int_0^{V_t} dx_a(v) \approx 1626.11.
$$ 
Thus the train takes $83$ s to reach speed $30$ ms$^{-1}$ and travels $1626$ m.  We also have
$$
\tau_{c,t} = \int_{U_b(V_t)}^{V_t} |dt_c(v)| \approx 147.33 \quad \mbox{and} \quad \xi_{c,t} = \int_{U_b(V_t)}^{V_t} |dx_c(v)| \approx 2857.28.
$$
Hence the train takes $147$ s to coast from a typical optimal driving speed $V_t = 30$ ms$^{-1}$ to the associated optimal braking speed $U_b(V_t) = \psi(V_t)/\varphi^{\, \prime}(V_t) \approx 9.62$ ms$^{-1}$ and travels $2857$ m.   The advantage of the constant-speed model\textemdash with or without speed penalties\textemdash  is that the optimal schedules can be conveniently computed using a single multi-dimensional Newton iteration.

Once we have an initial feasible schedule for the realistic strategies\textemdash using the weighted optimal constant-speed schedule\textemdash we can find an optimal schedule using the method of steepest descent.  At each step the optimal strategies between stops for each train can be computed by a Newton iteration\textemdash an elegant but much more protracted procedure than that required for the constant-speed schedules.  The case studies show that the realistic schedule is very similar in structure to the constant-speed schedule.  In each case the trains must go faster on the longer sections. 

\begin{figure}[htb]
\begin{center}
\includegraphics[width=8.1cm]{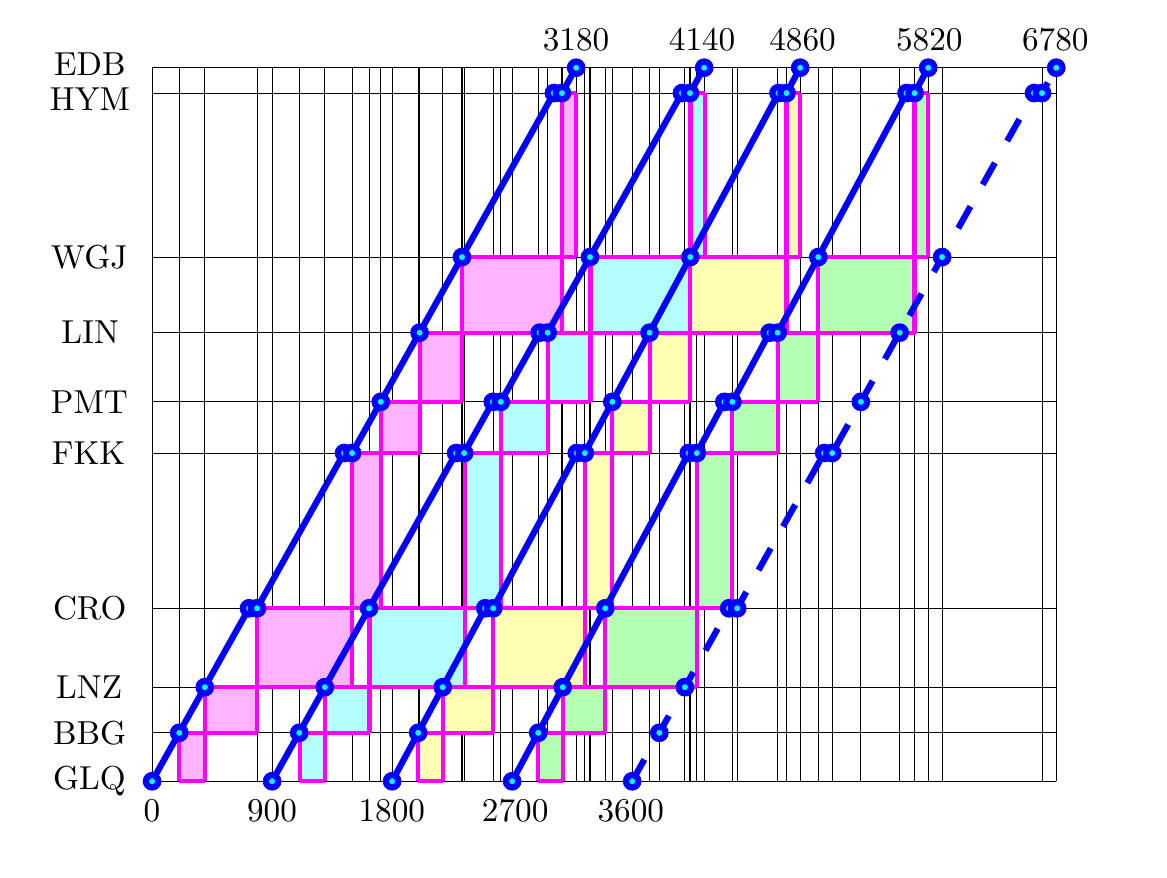}
\includegraphics[width=8.1cm]{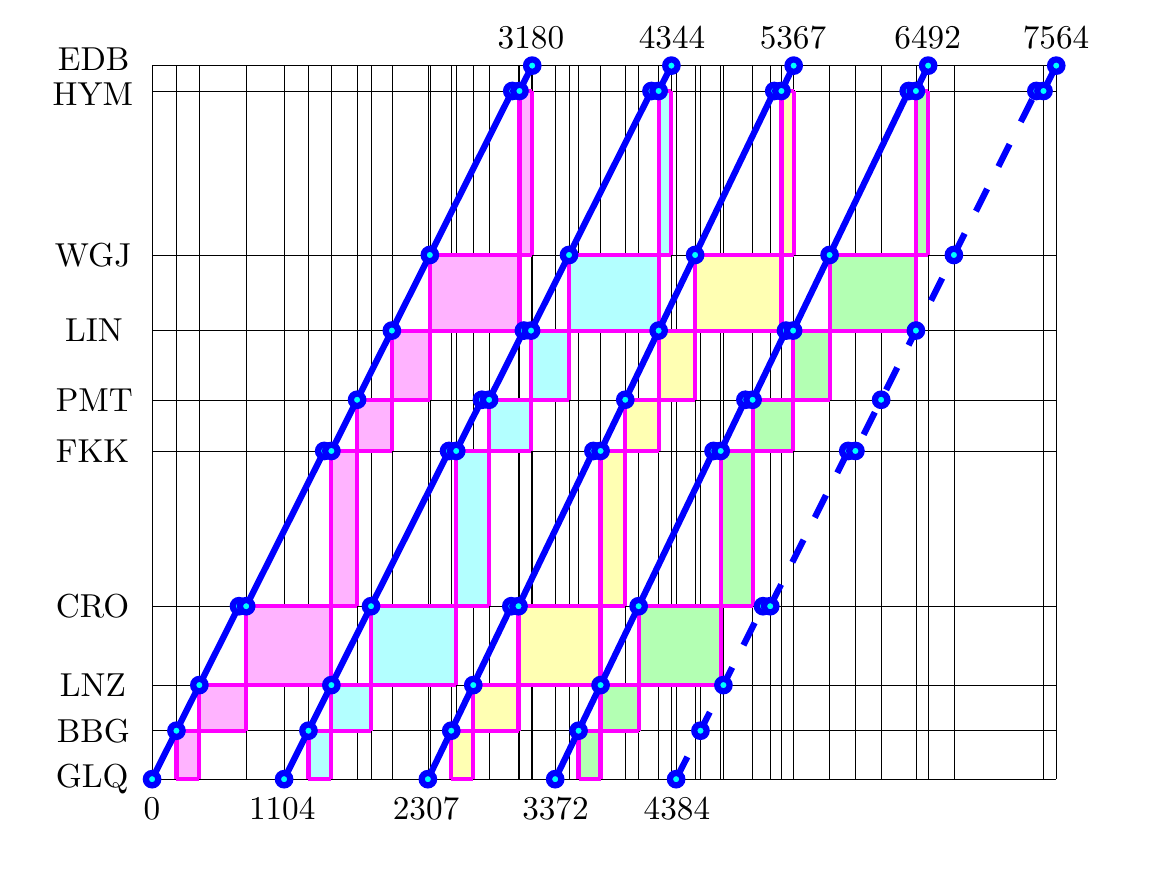}
\end{center}
\caption{\boldmath \small \bf Case Study~\ref{cs1}.  Train graphs for the regular shuttle service from GLQ to EDB using optimal constant-speed strategies with no separation constraints (left) and with delayed departures to satisfy the separation conditions (right).  The shaded rectangles show that the unrestricted optimal strategies do not satisfy the separation rules but that the delayed strategies do.  The cost is $L \approx 45156$ {J kg}$^{-1}$ in each case.  The line-occupancy timespan increases from $5820$ to $6492$~s when the separation conditions are imposed.}
\label{fig2}
\end{figure}

\begin{cs}
\label{cs1}

We use the published departure times for GLQ and EDB shown in Table~\ref{tab3}.  We assume each train travels at the same constant speed on every section.  The journey times are $3180$~s for ${\mathfrak T}_1$, $3240$~s for ${\mathfrak T}_2$, $3060$~s for ${\mathfrak T}_3$ and $3120$~s for ${\mathfrak T}_4$.  Trains ${\mathfrak T}_1$ and ${\mathfrak T}_3$ stop at CRO, FKK and HYM.  Trains ${\mathfrak T}_2$ and ${\mathfrak T}_4$ stop at FKK, PMT, LIN and HYM.  Where stops are scheduled the stopping time is $\sigma = 60$~s.  The scheduled time is the departure time except at EDB.

The optimal speeds are $W_1 = X/(T_1 - 3 \sigma) = 75700/3000 \approx 25.23$~ms$^{-1}$ for ${\mathfrak T}_1$, $W_2 = X/(T_2 - 4 \sigma) = 75700/3000 \approx 25.23$~ms$^{-1}$ for ${\mathfrak T}_2$, $W_3 = X/(T_3 - 3 \sigma) = 75700/2880 \approx 26.28$~ms$^{-1}$ for ${\mathfrak T}_3$ and $W_4= X/(T_4 - 4\sigma) = 75700/2880  \approx 26.28$~ms$^{-1}$ for ${\mathfrak T}_4$.  The total cost is
$$
L = [r(W_1) +r(W_2)+r(W_3)+r(W_4)]X \approx 45156\ \mbox{{J kg}}^{-1}.
$$
For costing purposes we assume the trains have equal masses.  The optimal schedules are shown in Table~\ref{tab4}.  The brackets indicate that the train does not stop.  The train graphs are shown on the left in Figure \ref{fig2}.  The shaded rectangles show that the separation conditions are not always satisfied.  To satisfy the separation conditions using these strategies the departure times must be delayed.  The changes are $[0900, 4140] \rightarrow [1104, 4340]$~s for ${\mathfrak T}_2$, $[1800, 4860] \rightarrow [2307, 5367]$~s for ${\mathfrak T}_3$, $[2700, 5820] \rightarrow [3372, 6492]$~s for ${\mathfrak T}_4$ and $[3600, 6780] \rightarrow [4384, 7564]$~s for the repeat service ${\mathfrak T}_{1,\star}$ in the next cycle.  The line-occupancy timespan for the fleet of four trains increases from $5820$ to $6492$~s.  The train graphs for the separated schedule are shown on the right in Figure~\ref{fig2}. $\hfill \Box$
\end{cs}

\begin{table}[htb]
\begin{center}
\begin{tabular}{|c|c|c|c|c|c|} \hline
Station & ${\mathfrak T}_1$ & ${\mathfrak T}_2$ & ${\mathfrak T}_3$ & ${\mathfrak T}_4$ & ${\mathfrak T}_{1,\star}$ \\ \hline
GLQ & $t_{1,0} = 0000$ & $t_{2,0} = 0900$ & $t_{3,0} = 1800$ & $t_{4,0} = 2700$ & $t_{1,\star,0} = 3600$ \\ \hline
BBG & $(t_{1,1} = 0204)$ & $(t_{2,1} = 1104)$ & $(t_{3,1} = 1996)$ & $(t_{4,1} = 2896)$ & $(t_{1,\star,1} = 3804)$ \\ \hline
LNZ & ($t_{1,2} = 0396)$ & {\color{red}\dag} $(t_{2,2} = 1296)$ & {\color{red}\dag} $(t_{3,2} = 2180)$ & {\color{red}\dag} $(t_{4,2} = 3080)$ & {\color{red}\dag} $(t_{1,\star,2} = 3996)$ \\ \hline
CRO & $t_{1,3} = 0787$ & {\color{red}\dag} $(t_{2,3} = 1627)$ & {\color{red}\dag} $t_{3,3} = 2558$ & {\color{red}\dag} $(t_{4,3} = 3398)$ & $(t_{1,\star,3} = 4387)$ \\ \hline
FKK & $t_{1,4} = 1500$ & $t_{2,4} = 2340$ & $t_{3,4} = 3245$ & $t_{4,4} = 4085$ & $t_{1,\star,4} = 5100$ \\ \hline
PMT & $(t_{1,5} = 1715)$ & $t_{2,5} = 2615$ & $(t_{3,5} = 3451)$ & $t_{4,5} = 4351$ & $(t_{1,\star,5} = 5315)$ \\ \hline
LIN & $(t_{1,6} = 2006)$ & {\color{red}\dag} $t_{2,6} = 2966$ & {\color{red}\dag} $(t_{3,6} = 3731)$ & {\color{red}\dag} $t_{4,6} = 4691$ & {\color{red}\dag} $(t_{1,\star,6} = 5605)$ \\ \hline
WGJ & $(t_{1,7} = 2324)$ & $(t_{2,7} = 3284)$ & {\color{red}\dag} $(t_{3,7} = 4036)$ & $(t_{4,7} = 4996)$ & $(t_{1,\star,7} = 5924)$ \\ \hline
HYM & $t_{1,8} = 3073$ & $t_{2,8} = 4033$ & $t_{3,8} = 4758$ & $t_{4,8} = 5718$ & $t_{1,\star,8} = 6673$ \\ \hline
EDB & $t_{1,9} = 3180$ & $t_{2,9} = 4140$ & $t_{3,9} = 4860$ & $t_{4,9} = 5820$ & $t_{1,\star,9} = 6780$ \\ \hline
\end{tabular}
\end{center}
\vspace{0.4cm}
\caption{ \small \bf Case Study~\ref{cs1}.  Optimal constant-speed strategies for each train using the scheduled starting times.  The daggers {\color{red}\dag} show instances where a train is scheduled to reach a signal location too soon.}
\label{tab4}
\end{table}

\begin{remark}
\label{cs1r}
Case Study~\ref{cs1} shows that if each train uses an individual optimal strategy the safe separation conditions will be violated.   If the departure times for subsequent trains are delayed to ensure safe separation the cycle time increases from 60 minutes to $74$ minutes.  $\hfill \Box$
\end{remark}

\begin{figure}[htb]
\begin{center}
\includegraphics[width=8.1cm]{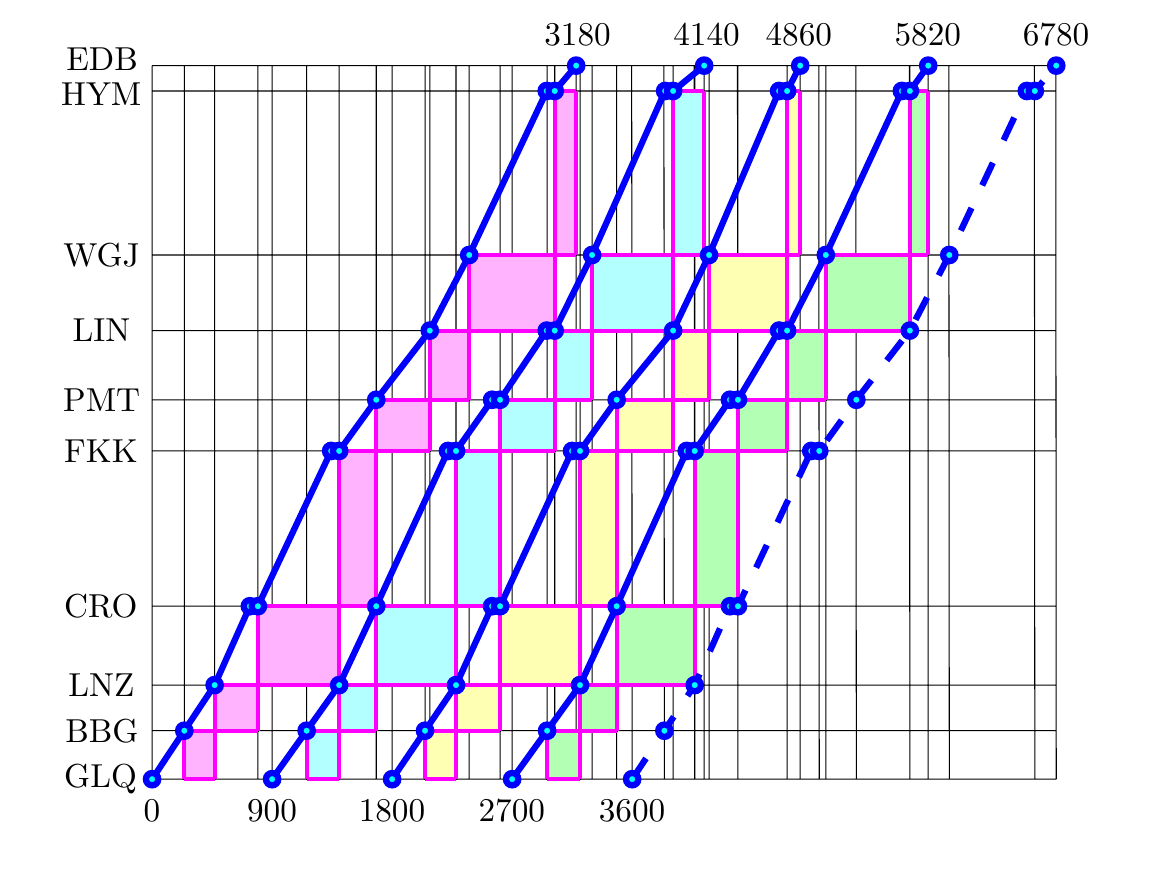}
\includegraphics[width=8.1cm]{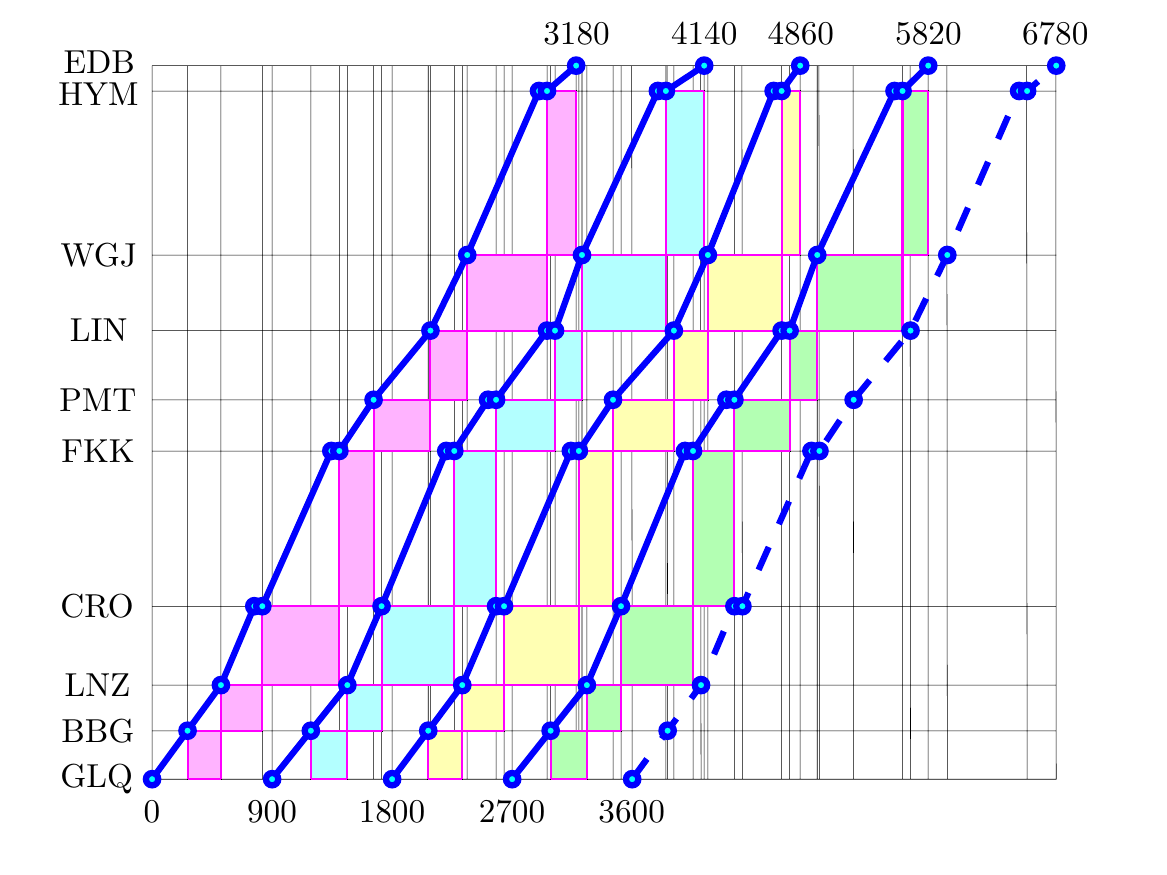}
\end{center}
\caption{\boldmath \small \bf Case Study~\ref{cs2}.  Train graphs for the optimal constant-speed strategies with minimal required  separation (left) and optimal weighted constant-speed strategies with a $60$ s separation buffer (right).}
\label{fig3}
\end{figure}

\begin{cs}
\label{cs2}

We find optimal constant-speed strategies using the scheduled times for GLQ and EDB with selective constraints to eliminate violations of the separation conditions.  In Case Study~\ref{cs1} we showed $t_{i+1,2} < t_{i,4}$, $t_{i+1,3} < t_{i,5}$ and $t_{i+1,6} < t_{i,8}$ for each $i=1,2,3$ and $t_{1,\star,2} < t_{4,4}$ and $t_{1,\star,6} < t_{4,8}$ when there are no separation constraints.  To ensure minimal safe separation we set $t_{i+1,2} = t_{i,4}$, $t_{i+1,3} = t_{i,5}$ and $t_{i+1,6} = t_{i,8}$ for each $i=1,2,3$ and $t_{1,\star,j-1} = t_{4,j+1}$ for $j=1,2,5$.  The variables and constraints are shown in Table~\ref{tab5}.  The vector of unknown times is denoted by $\bfh = (h_1,\ldots,h_{12})$.  We assume that $W_{i,2} = W_{i,1}$ and $W_{i,8} = W_{i,7}$ for each $i=1,\ldots,4$.

\begin{table}[htb]
\begin{center}
\begin{tabular}{|c|c|c|c|c|c|} \hline
Station & ${\mathfrak T}_1$ & ${\mathfrak T}_2$ & ${\mathfrak T}_3$ & ${\mathfrak T}_4$ & ${\mathfrak T}_{1,\star}$ \\ \hline
GLQ & $t_{1,0} = 0$ & $t_{2,0} = 900$ & $t_{3,0} = 1800$ & $t_{4,0} = 2700$ & $t_{1,\star,0} = 3600$ \\ \hline
BBR & \textemdash & \textemdash & \textemdash & \textemdash & \textemdash \\ \hline
LNZ & $t_{1,2} = h_1$ & $t_{2,2} = h_3$ & $t_{3,2} = h_7$ & $t_{4,2} = h_{10}$ & $t_{1,\star,2} = h_1+3600$ \\ \hline
CRO & $t_{1,3} = h_2$ & $t_{2,3} = h_4$ & $t_{3,3} = h_8$ & $t_{4,3} = h_{11}$ & $t_{1,\star,3} = h_2+3600$ \\ \hline
FKK & $t_{1,4} = h_3$ & $t_{2,4} = h_7$ & $t_{3,4} = h_{10}$ & $t_{4,4} = h_1+3600$ & $t_{1,\star,4} = h_3+3600$ \\ \hline
PMT & $t_{1,5} = h_4$ & $t_{2,5} = h_8$ & $t_{3,5} = h_{11}$ & $t_{4,5} = h_2+3600$ & $t_{1,\star,5} = h_4+3600$ \\ \hline
LIN & $t_{1,6} = h_5$ & $t_{2,6} = h_6$ & $t_{3,6} = h_9$ & $t_{4,6} = h_{12}$ & $t_{1,\star,6} = h_5+3600$ \\ \hline
WGJ & \textemdash & \textemdash & \textemdash & \textemdash & \textemdash \\ \hline
HYM & $t_{1,8} = h_6$ & $t_{2,8} = h_9$ & $t_{3,8} = h_{12}$ & $t_{4,8} = h_5+3600$ & $t_{1,\star,8} = h_6+3600$ \\ \hline
EDB & $t_{1,9} = T_1$ & $t_{2,9} = T_2+900$ & $t_{3,9} = T_3+1800$ & $t_{4,9} = T_4+2700$ & $t_{5,9} = T_1+3600$ \\ \hline
\end{tabular}
\end{center}
\vspace{0.4cm}
\caption{\boldmath \small \bf Case Study~\ref{cs2}.  Departure time variables for the schedule optimization subject to minimal required separation.}
\label{tab5}
\end{table}

The speeds for train ${\mathfrak T}_i$ for each $i=1,\ldots,4$ are given by the vector $\bfW_i = \bfW_i(\bfh)$ with components
$$
W_{i,j} = \left\{ \begin{array}{ll}
(x_2 - x_0)/(t_{i,2} - \sigma_{i,2} - t_{i,0}) & \mbox{for}\ j=1,2 \\
(x_j - x_{j-1})/(t_{i,j} - \sigma_{i,j} - t_{i,j-1}) & \mbox{for}\ j=3,\ldots,6 \\
(x_8 - x_6)/(t_{i,8} - \sigma_{i,8} - t_{i,6}) & \mbox{for}\ j=7,8 \\
(x_9 - x_8)/(t_{i,9} - t_{i,8}) & \mbox{for}\ j=9. \end{array} \right.
$$
where $\sigma_{i,j} = \sigma$ if ${\mathfrak T}_i$ stops at $x_j$ for each $j=1,\ldots,8$ and $\sigma_{i,j}  = 0$ otherwise.  The total cost is
$$
L = \sum_{i=1}^4 \left[ \sum_{j=1}^9 r(W_{i,j})(x_j - x_{j-1}) \right]
$$
and the cost gradient vector is
$$
\bnab L = \left[ \begin{array}{l}
\psi_{13} - \psi_{11} + \psi_{45} - \psi_{44} \\
\psi_{14} - \psi_{13} + \psi_{46} - \psi_{45} \\
\psi_{15} - \psi_{14} + \psi_{23} - \psi_{21} \\
\psi_{16} - \psi_{15} + \psi_{24} - \psi_{23} \\
\psi_{17} - \psi_{16} + \psi_{49} - \psi_{47} \\
\psi_{19} - \psi_{17} + \psi_{27} - \psi_{26} \\
\psi_{25} - \psi_{24} + \psi_{33} - \psi_{31} \\
\psi_{26} - \psi_{25} + \psi_{34} - \psi_{33} \\
\psi_{29} - \psi_{27} + \psi_{37} - \psi_{36} \\
\psi_{35} - \psi_{34} + \psi_{43} - \psi_{41} \\
\psi_{36} - \psi_{35} + \psi_{44} - \psi_{43} \\
\psi_{39} - \psi_{37} + \psi_{47} - \psi_{46} \end{array} \right]
$$
where we have written $\psi_{i,j} = \psi(W_{i,j})$ for convenience.  We use a Newton iteration to solve the equation $\bnab L = \bfzero$.  The iteration converges rapidly.  The optimal speeds in ms$^{-1}$ are
$$
\bfW = \left[ \begin{array}{c}
\bfW_1 \\
\bfW_2 \\
\bfW_3 \\
\bfW_4 \end{array} \right] \approx \left[ \begin{array}{ccccccccc}
21.27 & 21.27 & 31.68 & 29.97 & 19.53 & 18.28 & 28.99 & 28.99 & 16.75 \\
19.84 & 19.84 & 30.11 & 30.60 & 20.07 & 20.98 & 30.71 & 30.71 & 11.55 \\
20.83 & 20.83 & 30.94 & 30.53 & 19.76 & 17.34 & 31.99 & 31.99 & 27.35 \\
19.61 & 19.61 & 30.46 & 31.35 & 20.55 & 23.81 & 29.53 & 29.53 & 19.56 \end{array} \right].
$$
The optimal departure times are shown in Table~\ref{tab6} and the corresponding train graphs are displayed on the left in Figure~\ref{fig3}.  The total cost is $L \approx 46295$ {J kg}$^{-1}$.  

\begin{table}[htb]
\begin{center}
\begin{tabular}{|c|c|c|c|c|} \hline
Station & ${\mathfrak T}_1$ & ${\mathfrak T}_2$ & ${\mathfrak T}_3$ & ${\mathfrak T}_4$ \\ \hline
GLQ & $t_{1,0} = 0000$ & $t_{2,0} = 0900$ & $t_{3,0} = 1800$ & $t_{4,0} = 2700$ \\ \hline
BBG & $(t_{1,1} = 0242)$ & $(t_{2,1} = 1159)$ & $(t_{3,1} = 2047)$ & $(t_{4,1} = 2962)$ \\ \hline
LNZ & ($t_{1,2} = 0469)$ & $(t_{2,2} = 1403)$ & $(t_{3,2} = 2279)$ & $(t_{4,2} = 3209)$ \\ \hline
CRO & $t_{1,3} = 0793$ & $(t_{2,3} = 1681)$ & $t_{3,3} = 2610$ & $(t_{4,3} = 3484)$ \\ \hline
FKK & $t_{1,4} = 1403$ & $t_{2,4} = 2279$ & $t_{3,4} = 3209$ & $t_{4,4} = 4069$ \\ \hline
PMT & $(t_{1,5} = 1681)$ & $t_{2,5} = 2610$ & $(t_{3,5} = 3484)$ & $t_{4,5} = 4393$ \\ \hline
LIN & $(t_{1,6} = 2083)$ & $t_{2,6} = 3019$ & $(t_{3,6} = 3907)$ & $t_{4,6} = 4762$ \\ \hline
WGJ & $(t_{1,7} = 2378)$ & $(t_{2,7} = 3300)$ & $(t_{3,7} = 4177)$ & $(t_{4,7} = 5052)$ \\ \hline
HYM & $t_{1,8} = 3019$ & $t_{2,8} = 3907$ & $t_{3,8} = 4762$ & $t_{4,8} = 5683$ \\ \hline
EDB & $t_{1,9} = 3180$ & $t_{2,9} = 4140$ & $t_{3,9} = 4860$ & $t_{4,9} = 5820$ \\ \hline
\end{tabular}
\end{center}
\vspace{0.4cm}
\caption{\boldmath \bf \small Case Study~\ref{cs2}.  Optimal constant-speed schedule with minimal separation.}
\label{tab6}
\end{table}
Ultimately we wish to find realistic strategies to implement the optimized schedule.  An effective schedule needs to allow for normal stochastic variation in journey times by providing additional buffer times between successive trains.  If we set  $t_{i+1,2} = t_{i,4}+60$, $t_{i+1,3} = t_{i,5}+60$ and $t_{i+1,6} = t_{i,8}+60$ for each $i=1,2,3$ we will provide at least a $60$~s buffer on each segment.  There may be some segments on the constant-speed schedule where the traversal time is less than the minimum time required by a realistic strategy.  If we wish to allow more time on $(x_{j-1},x_j)$ for ${\mathfrak T}_i$ we can repeat the optimization with $r(W_{i,j})$ replaced in the cost term by $r(W_{i,j} + p_{i,j})$ where $p_{i,j} > 0$ is a fixed penalty.  In this case we set $p_{13} = 2$, $p_{33} = 1$ and $p_{39} = 7$.  The variables and constraints are shown in Table~\ref{tab7}.

\begin{table}[htb]
\begin{center}
\begin{tabular}{|c|c|c|c|c|c|} \hline
Station & ${\mathfrak T}_1$ & ${\mathfrak T}_2$ & ${\mathfrak T}_3$ & ${\mathfrak T}_4$ & ${\mathfrak T}_{1,\star}$ \\ \hline
GLQ & $t_{1,0} = 0$ & $t_{2,0} = 900$ & $t_{3,0} = 1800$ & $t_{4,0} = 2700$ & $t_{1,\star,0} = 3600$ \\ \hline
BBR & \textemdash & \textemdash & \textemdash & \textemdash & \textemdash \\ \hline
LNZ & $t_{1,2} = h_1$ & $t_{2,2} = h_3+60$ & $t_{3,2} = h_7+60$ & $t_{4,2} = h_{10}+60$ & $t_{1,\star,2} = h_1+3600$ \\ \hline
CRO & $t_{1,3} = h_2$ & $t_{2,3} = h_4+60$ & $t_{3,3} = h_8+60$ & $t_{4,3} = h_{11}+60$ & $t_{1,\star,3} = h_2+3600$ \\ \hline
FKK & $t_{1,4} = h_3$ & $t_{2,4} = h_7$ & $t_{3,4} = h_{10}$ & $t_{4,4} = h_1+3540$ & $t_{1,\star,4} = h_3+3600$ \\ \hline
PMT & $t_{1,5} = h_4$ & $t_{2,5} = h_8$ & $t_{3,5} = h_{11}$ & $t_{4,5} = h_2+3540$ & $t_{1,\star,5} = h_4+3600$ \\ \hline
LIN & $t_{1,6} = h_5$ & $t_{2,6} = h_6+60$ & $t_{3,6} = h_9+60$ & $t_{4,6} = h_{12}+60$ & $t_{1,\star,6} = h_5+3600$ \\ \hline
WGJ & \textemdash & \textemdash & \textemdash & \textemdash & \textemdash \\ \hline
HYM & $t_{1,8} = h_6$ & $t_{2,8} = h_9$ & $t_{3,8} = h_{12}$ & $t_{4,8} = h_5+3540$ & $t_{1,\star,8} = h_6+3600$ \\ \hline
EDB & $t_{1,9} = T_1$ & $t_{2,9} = T_2+900$ & $t_{3,9} = T_3+1800$ & $t_{4,9} = T_4+2700$ & $t_{1,\star,9} = T_1+3600$ \\ \hline
\end{tabular}
\end{center}
\vspace{0.4cm}
\caption{\boldmath \small \bf Case Study~\ref{cs2}.  Departure time variables for a weighted schedule optimization with a $60$ s separation buffer.}
\label{tab7}
\end{table}

Once again the Newton iteration converges rapidly.  The optimal speeds in ms$^{-1}$ are
$$
\bfW = \left[ \begin{array}{c}
\bfW_1 \\
\bfW_2 \\
\bfW_3 \\
\bfW_4 \end{array} \right] \approx \left[ \begin{array}{ccccccccc}
19.21 & 19.21 & 33.43 & 31.74 & 21.09 & 16.58 & 31.68 & 31.68 & 12.89 \\
17.55 & 17.55 & 32.51& 33.99 & 21.43 & 18.96 & 33.20 & 33.20 & 9.51 \\
18.80 & 18.80 & 33.03 & 32.94 & 21.22 & 16.06 & 34.20 & 34.20 & 19.26 \\
17.69 & 17.69 & 32.71 & 34.35 & 21.69 & 20.94 & 31.44 & 31.44 & 15.72 \end{array} \right].
$$
The optimal departure times are shown in Table~\ref{tab8} and the corresponding train graphs are displayed on the right in Figure~\ref{fig3}.  The total cost is $L \approx 47307$ {J kg}$^{-1}$. $\hfill \Box$  
\end{cs}

\begin{table}[htb]
\begin{center}
\begin{tabular}{|c|c|c|c|c|} \hline
Station & ${\mathfrak T}_1$ & ${\mathfrak T}_2$ & ${\mathfrak T}_3$ & ${\mathfrak T}_4$ \\ \hline
GLQ & $t_{1,0} = 0000$ & $t_{2,0} = 0900$ & $t_{3,0} = 1800$ & $t_{4,0} = 2700$ \\ \hline
BBG & $(t_{1,1} = 0268)$ & $(t_{2,1} = 1193)$ & $(t_{3,1} = 2073)$ & $(t_{4,1} = 2991)$ \\ \hline
LNZ & ($t_{1,2} = 0519)$ & $(t_{2,2} = 1469)$ & $(t_{3,2} = 2331)$ & $(t_{4,2} = 3264)$ \\ \hline
CRO & $t_{1,3} = 0830$ & $(t_{2,3} = 1726)$ & $t_{3,3} = 2644$ & $(t_{4,3} = 3520)$ \\ \hline
FKK & $t_{1,4} = 1409$ & $t_{2,4} = 2271$ & $t_{3,4} = 3204$ & $t_{4,4} = 4059$ \\ \hline
PMT & $(t_{1,5} = 1666)$ & $t_{2,5} = 2584$ & $(t_{3,5} = 3460)$ & $t_{4,5} = 4370$ \\ \hline
LIN & $(t_{1,6} = 2109)$ & $t_{2,6} = 3031$ & $(t_{3,6} = 3917)$ & $t_{4,6} = 4780$ \\ \hline
WGJ & $(t_{1,7} = 2381)$ & $(t_{2,7} = 3232)$ & $(t_{3,7} = 4170)$ & $(t_{4,7} = 4994)$ \\ \hline
HYM & $t_{1,8} = 2971$ & $t_{2,8} = 3857$ & $t_{3,8} = 4720$ & $t_{4,8} = 5649$ \\ \hline
EDB & $t_{1,9} = 3180$ & $t_{2,9} = 4140$ & $t_{3,9} = 4860$ & $t_{4,9} = 5820$ \\ \hline
\end{tabular}
\end{center}
\vspace{0.4cm}
\caption{\boldmath \small \bf Case Study~\ref{cs2}.  Optimal weighted constant-speed schedule with a $60$ s separation buffer.}
\label{tab8}
\end{table}

\begin{remark}
\label{cs2r}
Case Study~\ref{cs2} shows that by inserting intermediate time constraints we can find a schedule with no violations and additional buffer times between trains that preserves the scheduled times at GLQ and EDB with only a marginal increase in cost.  The key observation is that the trains go faster on the longer sections.  The optimal schedule is a compromise between reducing the highest speeds and equalising the section traversal times. $\hfill \Box$
\end{remark}

\begin{cs}
\label{cs3}

We construct optimal realistic strategies to fit the optimal weighted constant-speed schedule with $60$ s buffer times obtained in Case Study~\ref{cs2}. 

{\bf Train \boldmath ${\mathfrak T}_1$:}  The first segment is from GLQ to CRO.  The time allowed from GLQ to LNZ is $t_{1,2} - t_{1,0} = 519$~s and from LNZ to CRO is $t_{1,3} - \sigma - t_{1,2} = 251$~s.  The strategies in Case Study~\ref{cs2} suggest that ${\mathfrak T}_1$ must travel faster on the section LNZ to CRO.  If so there are two possible forms for the optimal realistic strategy.  The two forms are described and justified in Appendix~\ref{s:itc}. 

If there is no speedhold phase on the segment from LNZ to CRO then ${\mathfrak T}_1$ will use a strategy of maximum acceleration, speedhold at speed $V_{11} = V_{12}$ on some segment $(a_{12}, b_{12}) \subset (x_0, x_2)$, maximum acceleration through $x_2$ with speed $U_{12} = U_s^{\dag}(U_{13},V_{12},V_{13})$ at $x = x_2$ to reach a maximum speed $V_{13}$ at some point $a_{13} \in (x_2, x_3)$, coast to speed $U_{13}$ at some point $\gamma_{13} \in (x_2, x_3)$ and maximum brake.   The unknown parameters are $V_{12}$, $V_{13}$ and $U_{13}$.  The distance travelled on $(x_2,x_3)$ is
$$
\xi_{1,3}(U_{13},V_{12},V_{13}) = \int_{U_{12}}^{V_{13}} dt_a(v) + \int_{U_{13}}^{V_{13}} |dt_c(v)| + \int_0^{U_{13}} |dt_b(v)|.
$$
The times taken on $(x_0,x_2)$ and $(x_2,x_3)$ are
$$
\tau_{1,2}(U_{13},V_{12}, V_{13}) = \int_0^{U_{12}} dt_a(v) + \frac{1}{V_{12}} \left[ x_2 - \int_0^{U_{12}} dx_a(v) \right],
$$  
\begin{eqnarray*}
\tau_{1,3}(U_{13},V_{12},V_{13}) & = & \int_{U_{12}}^{V_{13}} dt_a(v) + \int_{U_{13}}^{V_{13}} |dt_c(v)| + \int_0^{U_{13}} |dt_b(v)|
\end{eqnarray*}
respectively.  We apply a Newton iteration to solve the system of equations
\begin{eqnarray*}
\xi_{1,3}(U_{13},V_{12},V_{13}) - x_3 + x_2 & = & 0 \\
\tau_{1,2}(U_{13},V_{12},V_{13}) - t_{1,2} & = & 0 \\
\tau_{1,3}(U_{13},V_{12},V_{13}) - t_{1,3} + \sigma + t_{1,2} & = & 0
\end{eqnarray*}
with initial values $(U_{13},V_{12}, V_{13}) = (15, 25, 40)$.  We obtain $V_{11} = V_{12} \approx 17.4819$, $V_{13} \approx 41.2937$ and $U_{13} \approx 23.8935$~ms$^{-1}$.  We calculate $U_{12} \approx 35.2248$~ms$^{-1}$.  The cost is $J_{1,03} \approx 3113.2266$ {J kg}$^{-1}$.  The speed profile is shown in Figure~\ref{fig4}.   Since $U_{13} > U_b(V_{13}) = \psi(V_{13})/\varphi^{\, \prime}(V_{13}) \approx 17.4344$~ms$^{-1}$ this is the optimal strategy.  A strategy of optimal type with a speedhold phase on $(x_2,x_3)$ will not be feasible. 

\begin{figure}[htb]
\begin{center}
\includegraphics[width=5.4cm]{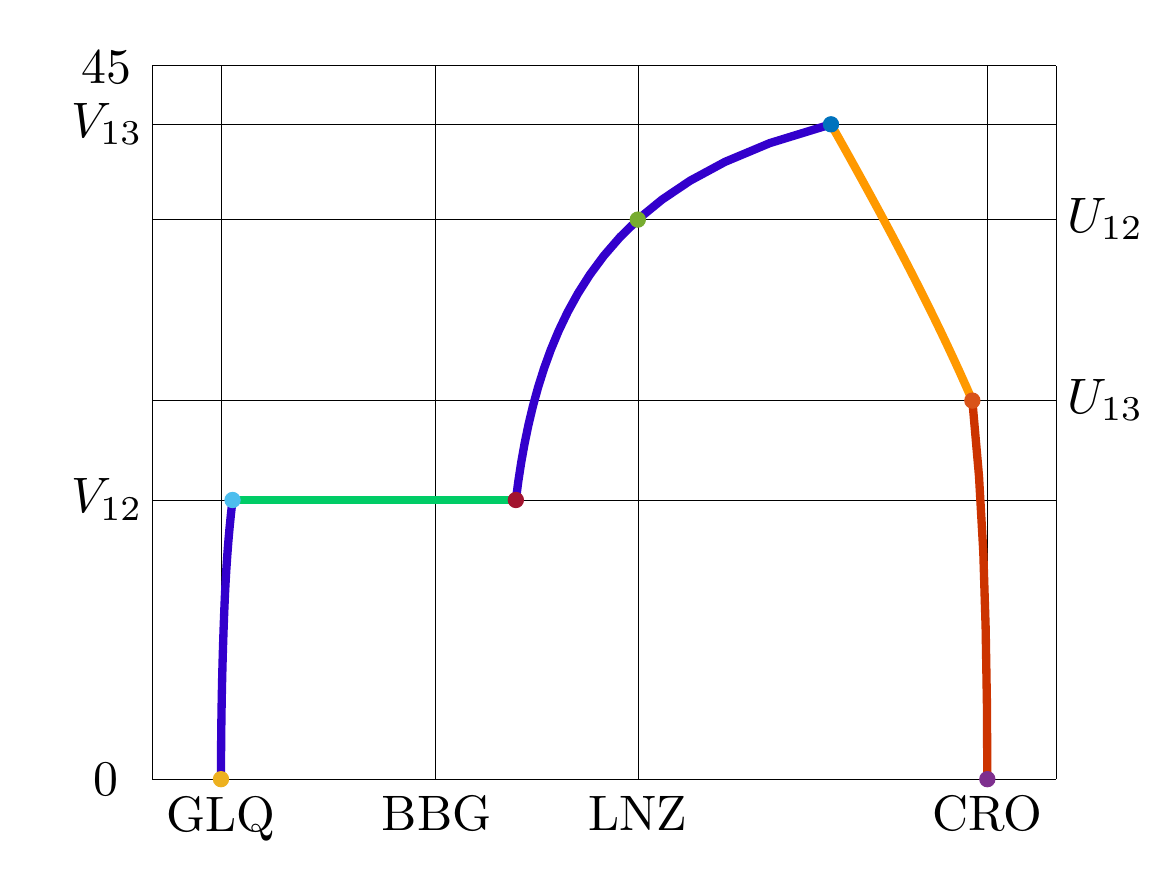}
\includegraphics[width=5.4cm]{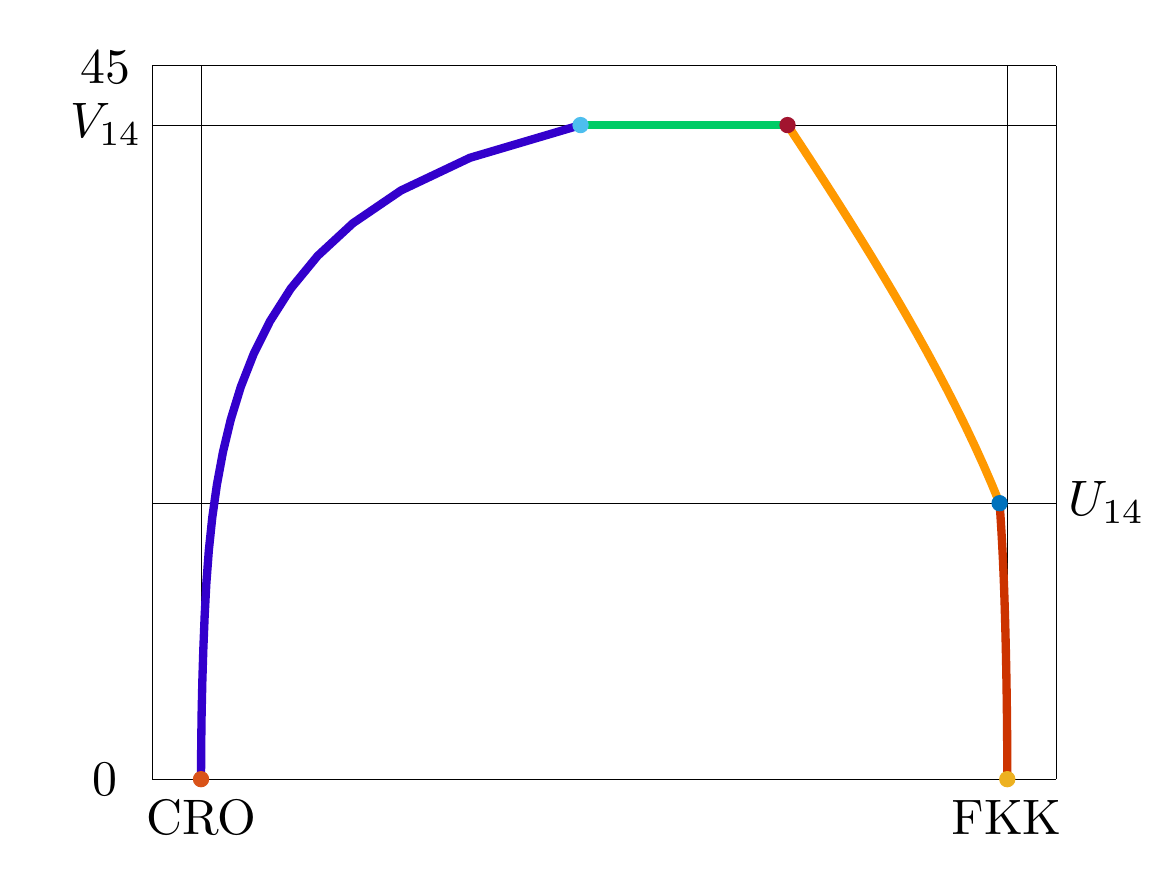}
\includegraphics[width=5.4cm]{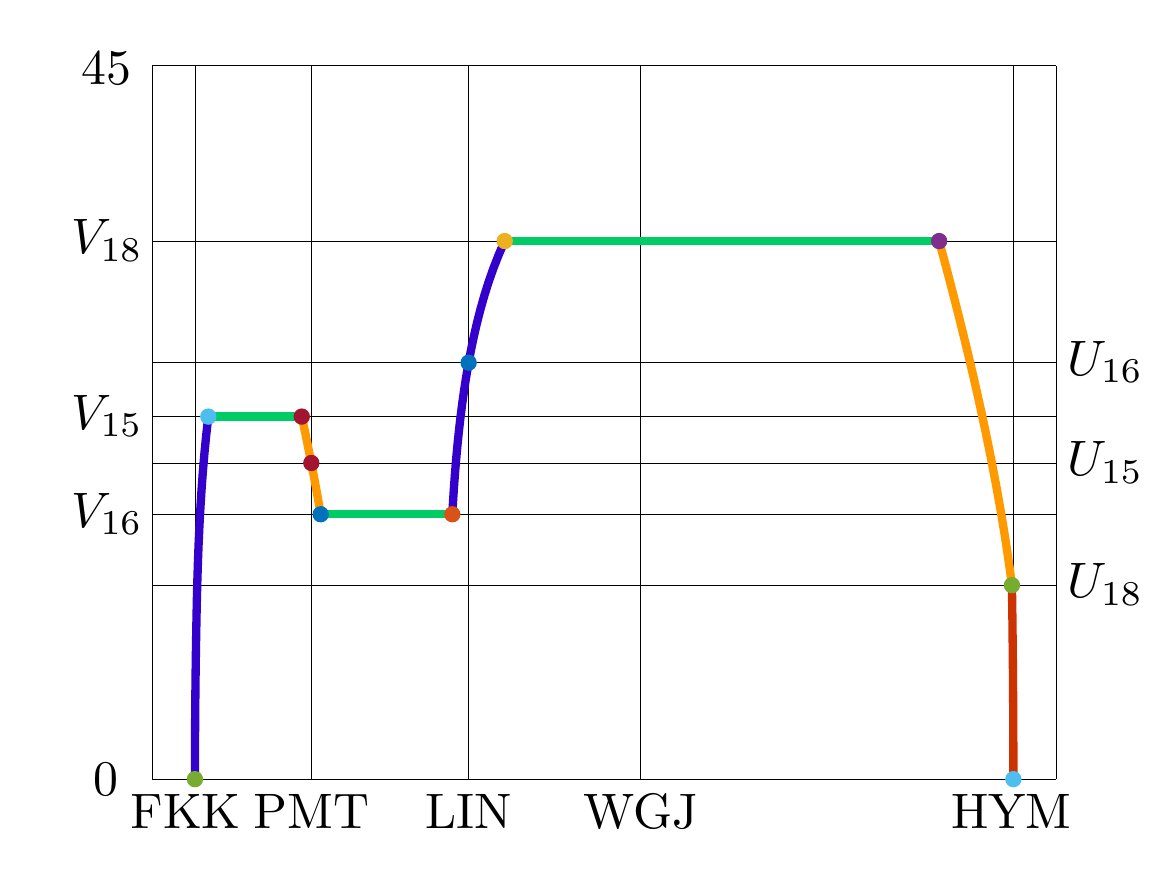}
\end{center}
\caption{\boldmath \small \bf Case Study~\ref{cs3}.  Optimal strategies for ${\mathfrak T}_1$ on the segments from GLQ to CRO, from CRO to FKK and from FKK to HYM passing through PMT, LIN and WGJ.  The horizontal scales are different on each graph.  Distances are $18350$~m from GLQ to CRO, $16470$~m from CRO to FKK and $38190$~m from FKK to HYM.}
\label{fig4}
\end{figure}

For the segment from CRO to FKK on $(x_3,x_4)$ train ${\mathfrak T}_1$ uses a long-haul strategy of maximum acceleration, speedhold at speed $V_{14}$, coast to speed $U_{14} = U_b(V_{14})$ and maximum brake.  The time allowed is $t_{1,4} - \sigma - t_{1,3} = 519$~s.   The time taken on $(x_3,x_4)$ is
\begin{eqnarray*}
\tau_{1,4}(V_{14}) & = & \int_0^{V_{14}} dt_a(v) + \int_{U_{14}}^{V_{14}} |dt_c(v)| + \int_0^{U_{14}} |dt_b(v)| \\
& & \hspace{3cm} \frac{1}{V_{14}} \left[ \int_0^{V_{14}} dx_a(v) + \int_{U_{14}}^{V_{14}} |dx_c(v)| + \int_0^{U_{14}} |dx_b(v)| \right].
\end{eqnarray*}
We use a Newton iteration to solve $\tau_{1,4}(V_{14}) - t_{1,4} + \sigma + t_{1,3} = 0$ with initial value $V_{14} = 40$.   We obtain $V_{14} \approx 41.2507$~ms$^{-1}$.  We calculate $U_{14} \approx 17.4032$~ms$^{-1}$.  The cost is $J_{1,34} \approx 3001.7179$ {J kg}$^{-1}$.   The speed profile is shown in Figure~\ref{fig4}.
 
For the segment from FKK to HYM on $(x_4, x_8)$ train ${\mathfrak T}_1$ uses a strategy of maximum acceleration, speedhold at speed $V_{15}$, coast through PMT with speed $U_{15} = U_s(V_{15},V_{16})$, speedhold at speed $V_{16}$, maximum acceleration through LIN with speed $U_{16} = U_s(V_{16}, V_{17})$, speedhold at speed $V_{17} = V_{18}$ through WGJ, coast to speed $U_{18} = U_b(V_{18})$ and maximum brake.  The time constraints are $t_{1,5} - t_{1,4} = 257$, $t_{1,6} - t_{1,5} = 443$, and $t_{1,8} - \sigma - t_{1,6} = 802$~s.  The respective times are
$$
\tau_{1,5}(V_{15}, V_{16}) = \int_0^{V_{15}} \hspace{-2mm} dt_a(v) + \int_{U_{15}}^{V_{15}} |dt_c(v)| + \frac{1}{V_{15}} \left[ x_5 - x_4 - \int_0^{V_{15}} \hspace{-2mm} dx_a(v) - \int_{U_{15}}^{V_{15}} |dx_c(v)| \right],
$$
$$
\tau_{1,6}(V_{15}, V_{16}, V_{18}) = \int_{V_{16}}^{U_{15}} |dt_c(v)| + \int_{V_{16}}^{U_{16}} \hspace{-2mm}  dt_a(v) + \frac{1}{V_{16}} \left[x_6 - x_5 - \int_{V_{16}}^{U_{15}} |dx_c(v)| - \int_{V_{16}}^{U_{16}} \hspace{-2mm} dx_a(v)  \right],
$$
$$
\tau_{1,8}(V_{16},V_{18}) = \int_{U_{18}}^{V_{18}} |dt_c(v)| + \int_0^{U_{18}}\hspace{-2mm} |dt_b(v)| + \frac{1}{V_{18}} \left[ x_8 - x_6 - \int_{U_{18}}^{V_{18}} |dx_c(v)| - \int_0^{U_{18}} \hspace{-2mm} |dx_b(v)| \right].
$$
We apply a Newton iteration to solve the system of equations
\begin{eqnarray*}
\tau_{1,5}(V_{15},V_{16}) - t_{1,5} + t_{1,4} & = & 0 \\
\tau_{1,6}(V_{15}, V_{16}, V_{18}) - t_{1,6} + t_{1,5} & = & 0 \\
\tau_{1,8}(V_{16}, V_{18}) - t_{1,8} + \sigma + t_{1,6} & = & 0 
\end{eqnarray*}
using the initial value $(V_{15}, V_{16}, V_{18}) = (21, 18, 30)$.  We obtain $V_{15} \approx 22.9024$, $V_{16} \approx 15.9710$ and $V_{17} = V_{18} \approx 34.5262$~ms$^{-1}$.  We calculate $U_{15} \approx 19.6342$, $U_{16} \approx 26.3486$ and $U_{18} \approx 12.6385$~ms$^{-1}$.  The cost is $J_{1,48} \approx 6070.0102$ {J kg}$^{-1}$.  The speed profile is shown in Figure~\ref{fig4}.

For the segment from HYM to EDB on $(x_8,x_9)$ train ${\mathfrak T}_1$ uses a long-haul strategy of maximum acceleration, speedhold at speed $V_{19}$, coast to speed $U_{19} = U_b(V_{19})$ and maximum brake.  The time constraint is $t_{1,9} - t_{1,8} = 209$~s.  The time taken is
\begin{eqnarray*}
\tau_{1,9}(V_{19}) & = & \int_0^{V_{19}} dt_a(v) + \int_{U_{19}}^{V_{19}} |dt_c(v)| + \int_0^{U_{19}} |dt_b(v)| \\
& & \hspace{3cm} + \frac{1}{V_{19}} \left[ x_9 - x_8 - \int_0^{V_{19}} dx_a(v) - \int_{U_{19}}^{V_{19}} |dx_c(v)| - \int_0^{U_{19}} |dx_b(v)| \right].
\end{eqnarray*}
We use a Newton iteration to solve $\tau_{1,9}(V_{19}) - t_{1,9} + t_{1,8} = 0$ with initial value $V_{19} = 15$.  We obtain  $V_{19} \approx 19.0591$~ms$^{-1}$.  We calculate $U_{19} \approx 3.5255$~ms$^{-1}$ and $J_{1,89} \approx 359.5285$ {J kg}$^{-1}$.  The speed profile is shown in Figure~\ref{fig5}.  The total cost for ${\mathfrak T}_1$ is
\begin{equation}
\label{t1cost1}
J_1 = J_{1,03} + J_{1,34} + J_{1,48} + J_{1,89} \approx 12544.4832\ \mbox{{J kg}}^{-1}.
\end{equation}

\begin{figure}[htb]
\begin{center}
\includegraphics[width=5.4cm]{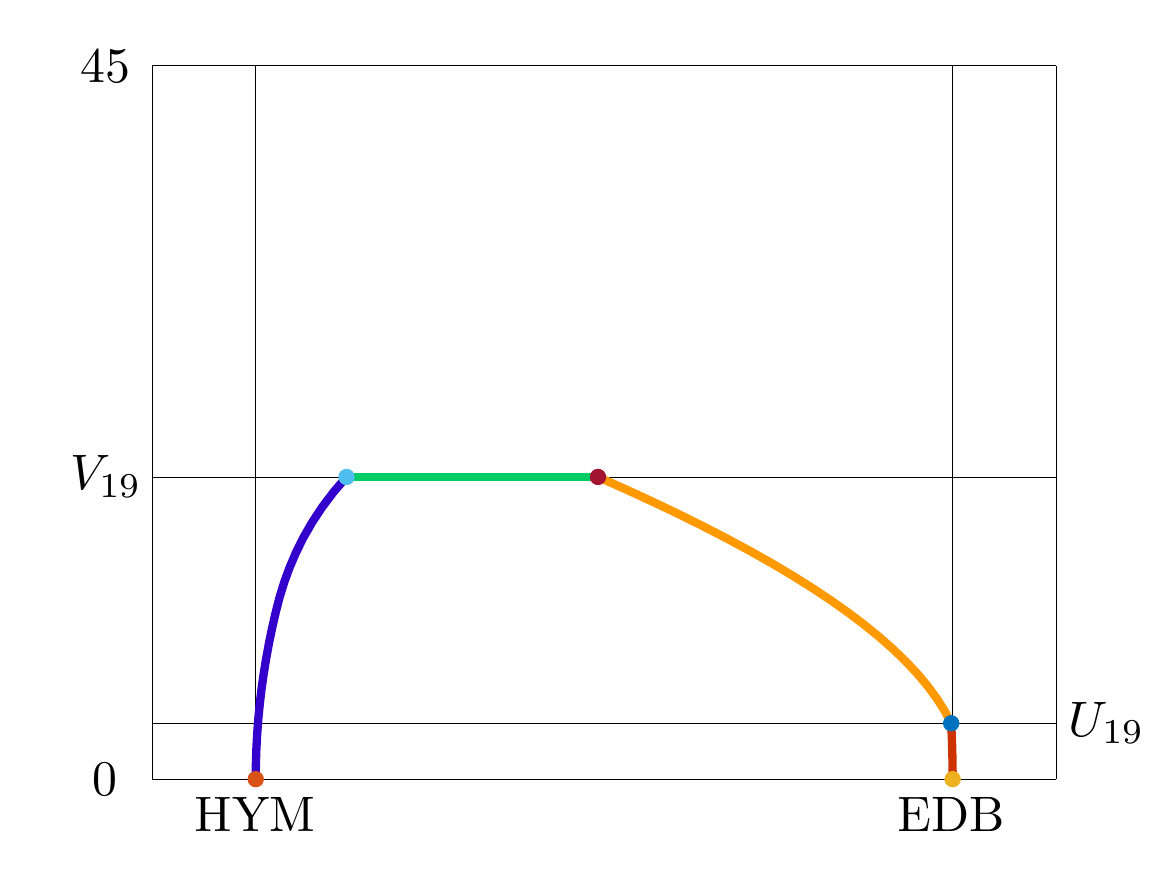}
\includegraphics[width=5.4cm]{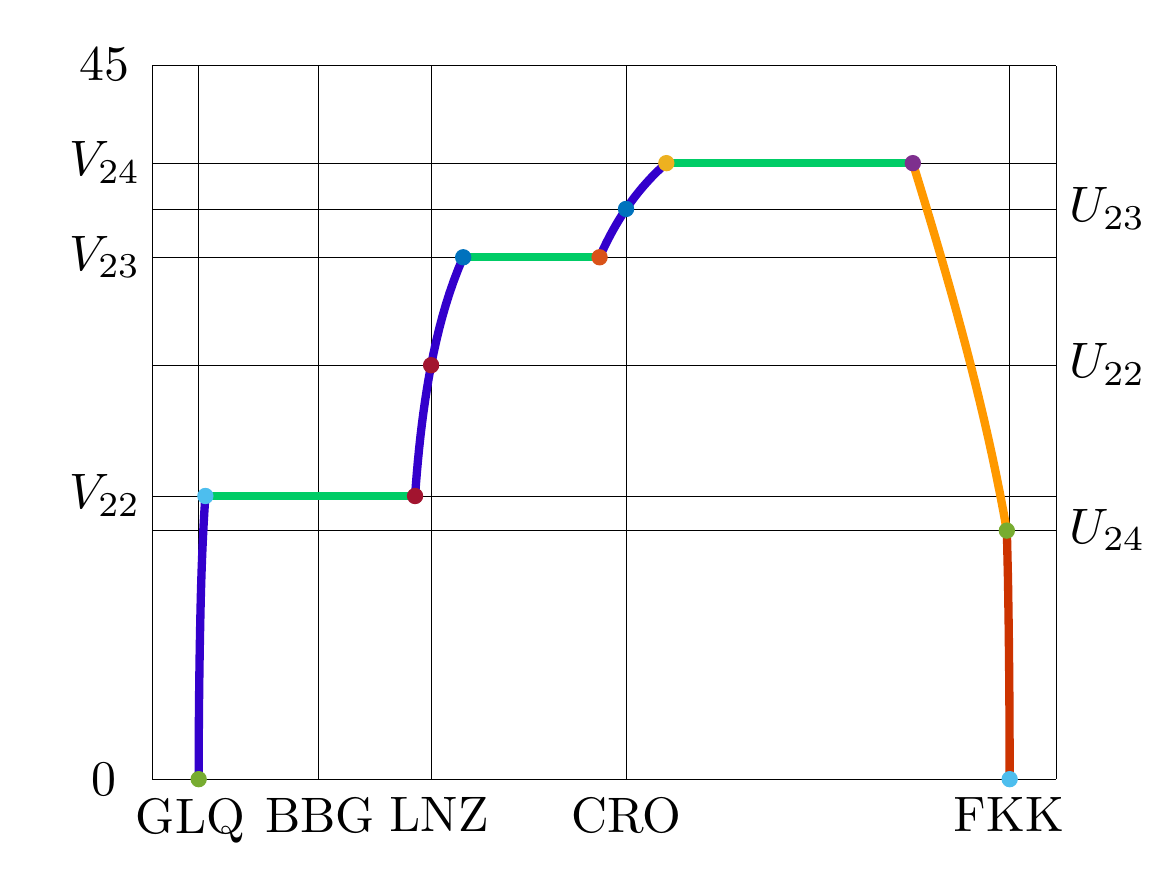}
\includegraphics[width=5.4cm]{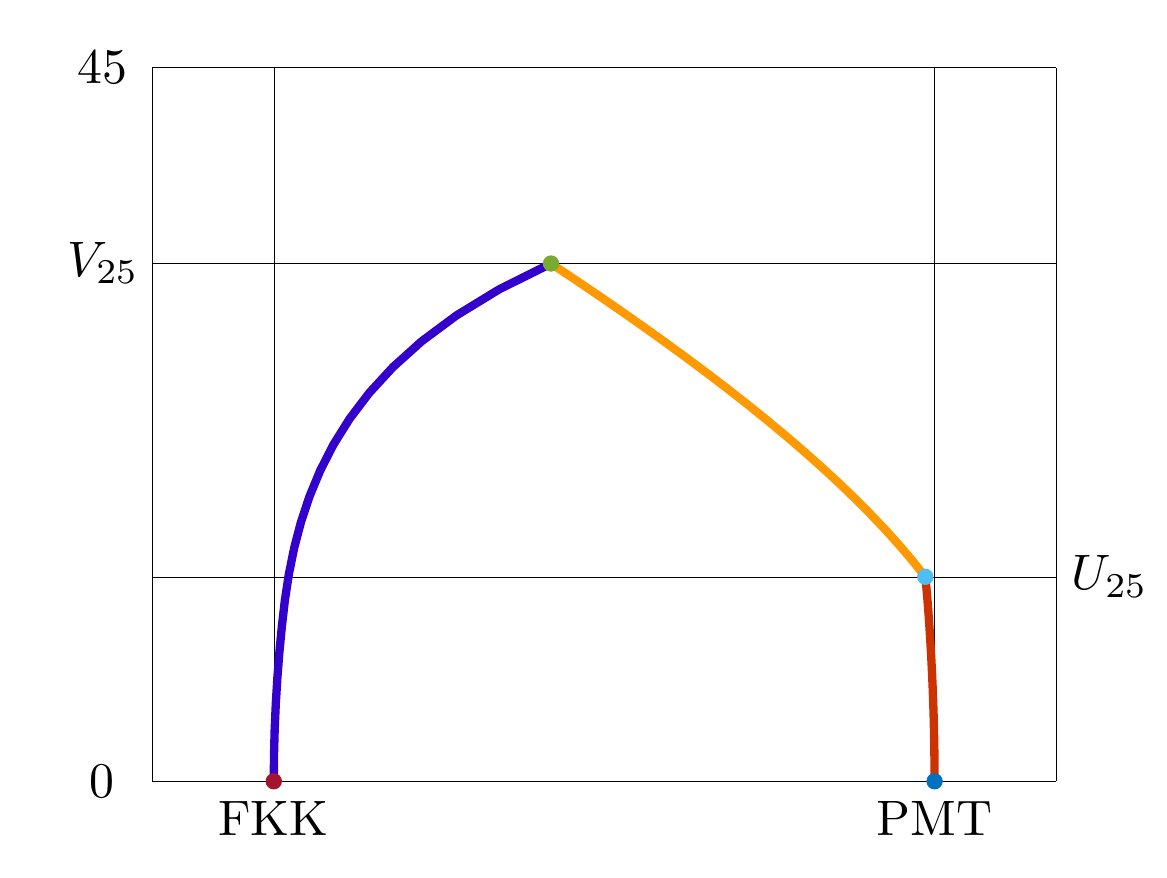}
\end{center}
\caption{\boldmath \small \bf Case Study~\ref{cs3}.  Optimal strategies for ${\mathfrak T}_1$ on the segment from HYM to EDB (left) and ${\mathfrak T}_2$ on the segments from GLQ to FKK (centre) and FKK to PMT (right).  The horizontal scales are different on each graph.  Distances are $2690$~m from HYM to EDB, $34820$~m from GLQ to FKK and $5430$~m from FKK to PMT.}
\label{fig5}
\end{figure}

{\bf \boldmath Train ${\mathfrak T}_2$:}  The strategies in Case Study~\ref{cs2} suggest that on the segment from GLQ  to FKK train ${\mathfrak T}_2$ should use a strategy of maximum acceleration, speedhold at speed $V_{21} = V_{22}$ through BBG at $x_1$, maximum acceleration through LNZ with speed $U_{22} = U_s(V_{22}, V_{23})$ at $x_2$, speedhold at speed $V_{23}$, maximum acceleration through CRO with speed $U_{23} = U_s(V_{23},V_{24})$ at $x_3$, speedhold at speed $V_{24}$, coast to speed $U_{24} =  U_b(V_{24})$ and maximum brake.  The time constraints are $t_{2,2} - t_{2,0} = 569$, $t_{2,3} - t_{2,2} = 257$, and $t_{2,4} - \sigma - t_{2,3} = 485$~s.  The respective times taken are
$$
\tau_{2,2}(V_{22}, V_{23}) = \int_0^{U_{22}} \hspace{-2mm} dt_a(v) +  \frac{1}{V_{22}} \left[ x_2 - x_0 - \int_0^{U_{22}} \hspace{-2mm} dx_a(v) \right],
$$
$$
\tau_{2,3}(V_{22},V_{23},V_{24}) = \int_{U_{22}}^{U_{23}} \hspace{-2mm} dt_a(v) + \frac{1}{V_{23}} \left[ x_3 - x_2 - \int_{U_{22}}^{U_{23}} \hspace{-2mm} dx_a(v) \right],
$$
\begin{eqnarray*}
\tau_{2,4}(V_{23},V_{24}) & = & \int_{U_{23}}^{V_{24}} \hspace{-2mm} dt_a(v) + \int_{U_{24}}^{V_{24}} \hspace{-2mm} |dt_c(v)| + \int_0^{U_{24}} \hspace{-2mm} |dt_b(v)| \\
& &\hspace{2cm} + \frac{1}{V_{24}} \left[ x_4 - x_3 - \int_{U_{23}}^{V_{24}} \hspace{-2mm} dx_a(v) - \int_{U_{24}}^{V_{24}} \hspace{-2mm} |dx_c(v)| - \int_0^{U_{24}} \hspace{-2mm} |dx_b(v)| \right].
\end{eqnarray*}
We apply a Newton iteration to solve the system of equations
\begin{eqnarray*}
\tau_{2,2}(V_{22},V_{23}) - t_{2,2} + t_{2,0} & = & 0 \\
\tau_{2,3}(V_{22},V_{23},V_{24}) - t_{2,3} + t_{2,2} & = & 0 \\
\tau_{2,4}(V_{23},V_{24}) - t_{2,4} + \sigma + t_{2,3} & = & 0
\end{eqnarray*}
using initial values $(V_{22},V_{23},V_{24}) = (20, 30, 35)$.  We obtain $V_{21} = V_{22} \approx 17.6871$, $V_{23} \approx 32.9234$ and $V_{24} \approx 38.8491$~ms$^{-1}$.  We calculate $U_{22} \approx 26.0454$, $U_{23} \approx 35.9660$, $U_{24} \approx 15.6746$~ms$^{-1}$ and $J_{2,04} \approx 5750.0444$ {J kg}$^{-1}$.  The speed profile is shown in Figure~\ref{fig5}.

\begin{figure}[htb]
\begin{center}
\includegraphics[width=5.4cm]{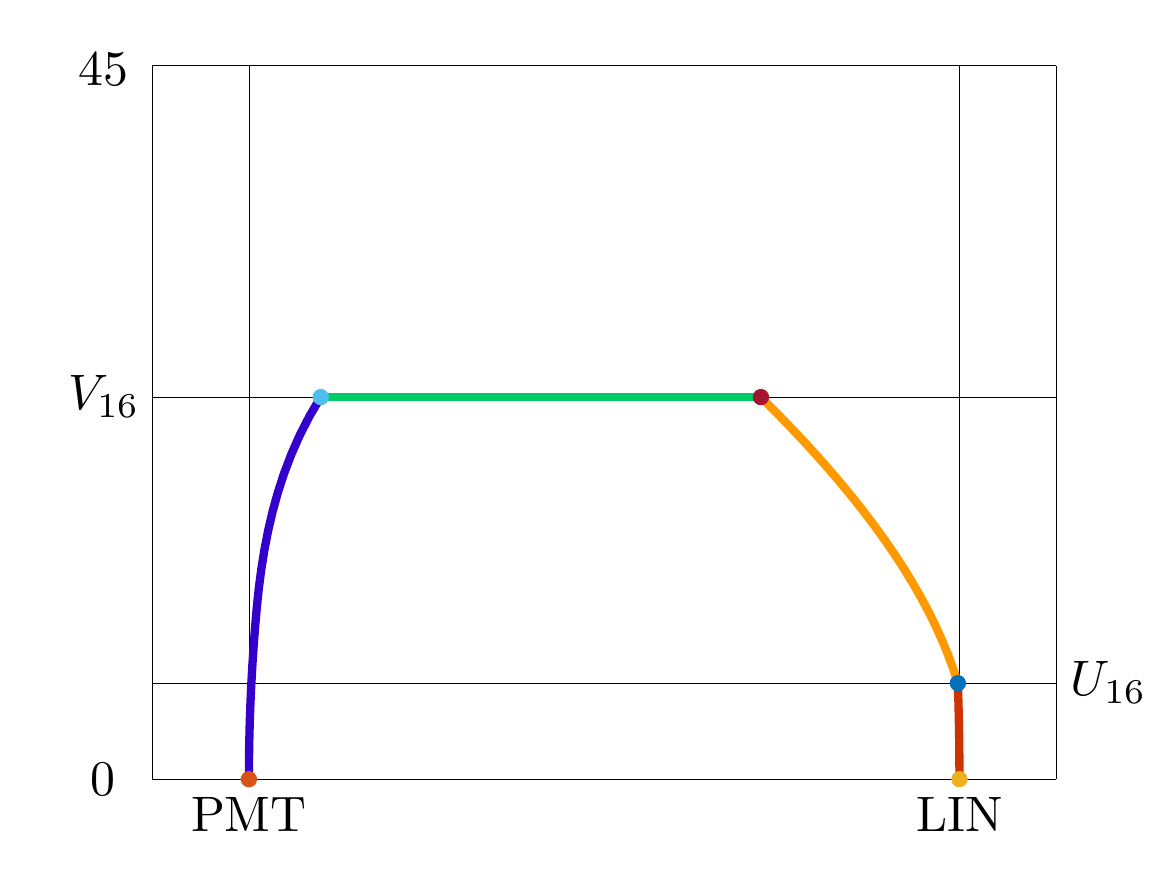}
\includegraphics[width=5.4cm]{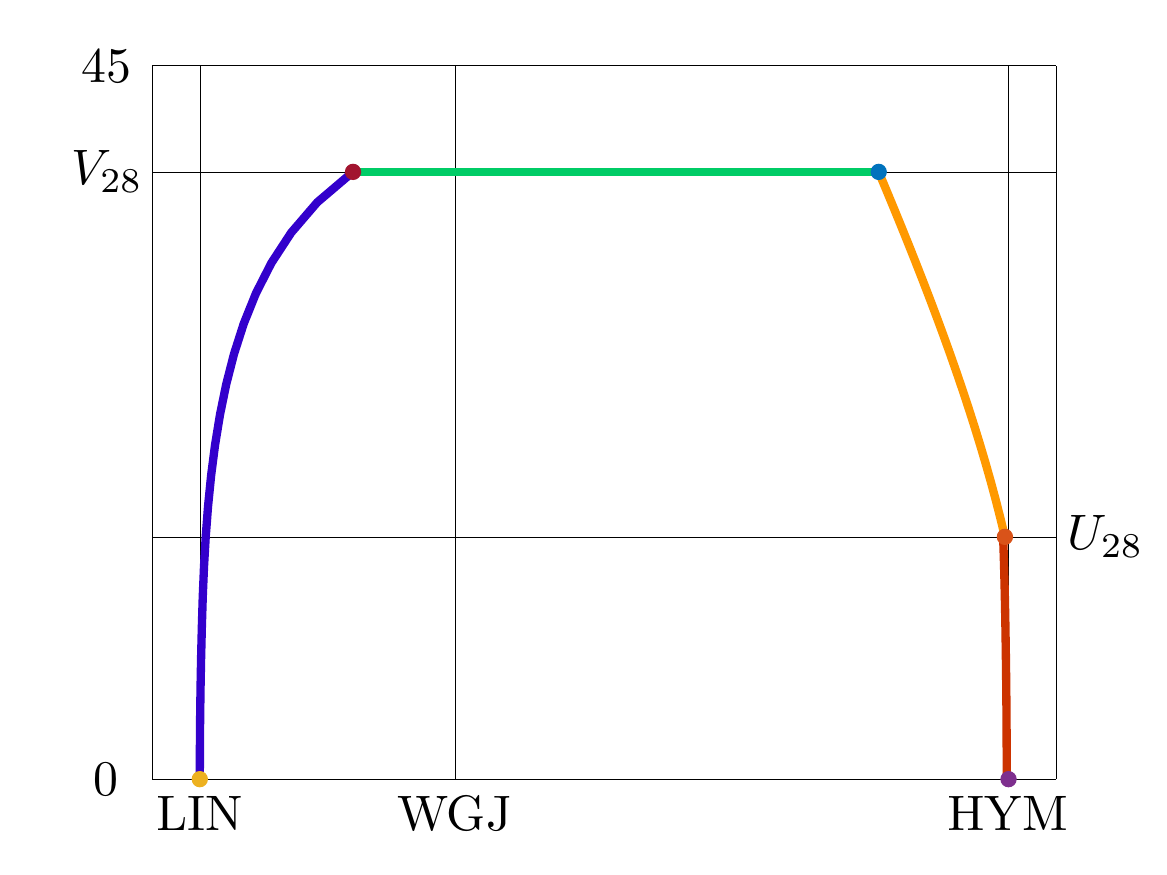}
\includegraphics[width=5.4cm]{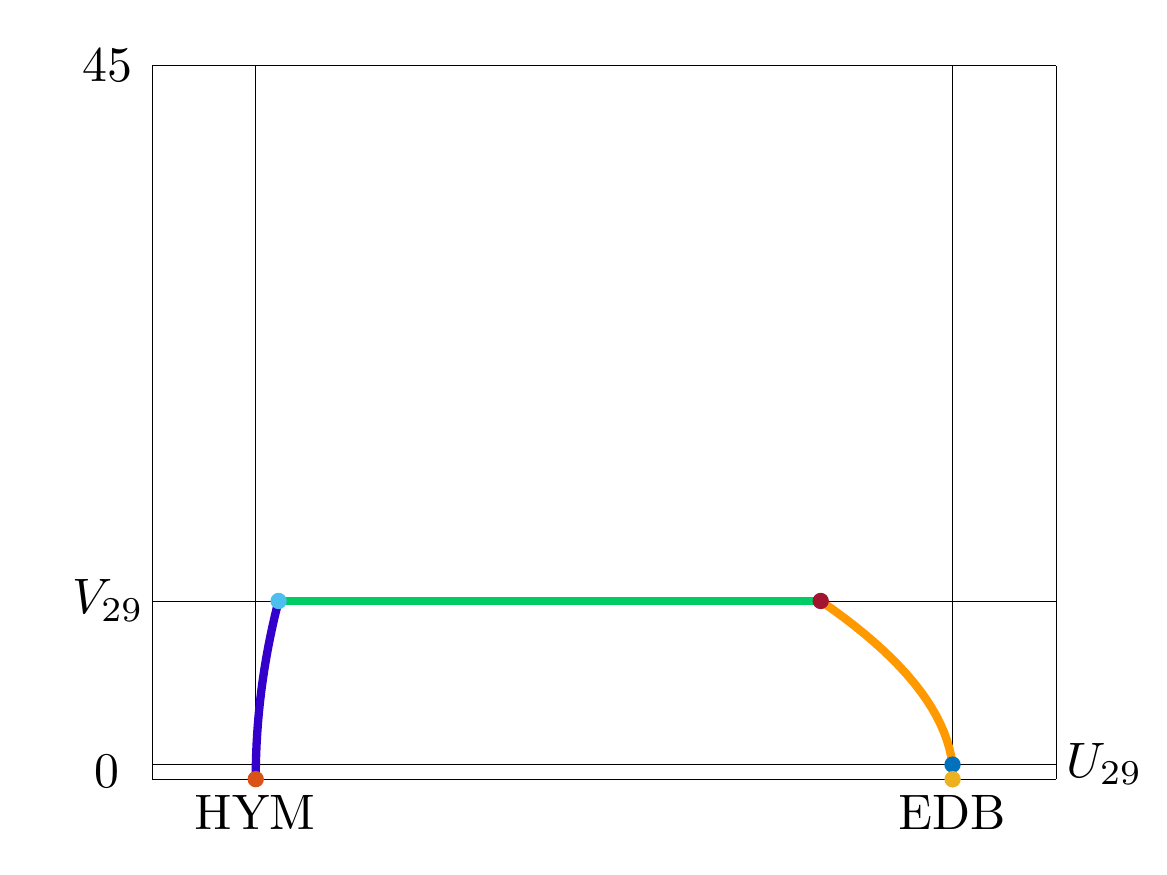}
\end{center}
\caption{\boldmath \small \bf Case Study~\ref{cs3}.  Optimal strategies for ${\mathfrak T}_2$ on the segments from PMT to LIN (left), from LIN to HYM (centre) and from HYM to EDB (right).   The horizontal scales are different on each graph.  Distances are $7340$~m from PMT to LIN, $25420$~m from LIN to HYM and $2690$~m from HYM to EDB.}
\label{fig6}
\end{figure}

The optimal strategy for ${\mathfrak T}_2$ from FKK to PMT is a rapid-transit strategy with $U_{25} \approx  
12.8945$ and $V_{25} \approx 32.6508$~ms$^{-1}$ and $J_{2,45} \approx 874.9529$ {J kg}$^{-1}$.  The speed profile is shown in Figure~\ref{fig5}.  The remaining strategies are all standard long-haul strategies from PMT to LIN with $V_{26} \approx 23.6455$~ms$^{-1}$, $U_{26} = U_b(V_{26}) \approx 5.8051$~ms$^{-1}$ and $J_{2,56} \approx 1047.6125$ {J kg}$^{-1}$, from LIN to HYM with $V_{28} \approx 38.7348$~ms$^{-1}$, $U_{28} = U_b(V_{28}) \approx 15.5929$~ms$^{-1}$ and $J_{2,68} \approx 4563.8225$ {J kg}$^{-1}$, and from HYM to EDB with $V_{29} \approx 11.4642$~ms$^{-1}$, $U_{29} = U_b(V_{29}) \approx 0.9682$~ms$^{-1}$ and $J_{2,89} \approx 338.5729$ {J kg}$^{-1}$.  The speed profiles are shown in Figure~\ref{fig6}.  The total cost for ${\mathfrak T}_2$ is
\begin{equation}
\label{t2cost1}
J_2 = J_{2,04} + J_{2,45} + J_{2,56} + J_{2,68} + J_{2,89} \approx 12575.0052\ \mbox{{J kg}}^{-1}. 
\end{equation}

{\bf \boldmath Train ${\mathfrak T}_3$:}  The strategies for ${\mathfrak T}_3$ are similar to those for ${\mathfrak T}_1$.  The optimal driving speeds from GLQ to CRO are $V_{31} = V_{32} \approx 17.1769$, $V_{33} \approx 41.1474$, $U_{33} \approx 23.3292$, $U_{32} = U_s^{\dag}(U_{33},V_{23},V_{33}) \approx 34.8388$~ms$^{-1}$ and $J_{3,03} \approx 3086.6375$ {J kg}$^{-1}$.  The speed profile is shown in Figure~\ref{fig7}.  
The optimal strategy from CRO to FKK is a rapid-transit strategy with $U_{34} \approx 22.7752$ and $V_{34} \approx 43.2069$~ms$^{-1}$ and with $J_{3,34} \approx 3139.5673$ {J kg}$^{-1}$.  The optimal strategy for ${\mathfrak T}_3$ from FKK to HYM takes the same form as the strategy for ${\mathfrak T}_1$.  The speeds are $V_{35} \approx 23.0633$, $V_{36} \approx 15.1369$ and $V_{37} = V_{38} \approx 37.7486$~ms$^{-1}$.  We calculate $U_{35} = U_s(V_{35},V_{36}) \approx 19.3628$, $U_{36} = U_s(V_{36},V_{38}) \approx 28.0048$ and $U_{38} = U_b(V_{38}) \approx 14.8912$~ms$^{-1}$.  The cost is $J_{3,48} \approx 6296.2674$ {J kg}$^{-1}$.  The speed profile is shown in Figure~\ref{fig7}.  The optimal strategy for the segment from HYM to EDB is a rapid-transit strategy with $U_{39} \approx 21.2930$ and $V_{39} \approx 28.12294$~ms$^{-1}$ and $J_{3,89} \approx 576.9359$ {J kg}$^{-1}$.   The total cost for ${\mathfrak T}_3$ is
\begin{equation}
\label{t3cost1}
J_3 = J_{3,03} + J_{3,34} + J_{3,48} + J_{3,89} \approx 13099.4081\ \mbox{{J kg}}^{-1}. 
\end{equation}

\begin{figure}[htb]
\begin{center}
\includegraphics[width=5.4cm]{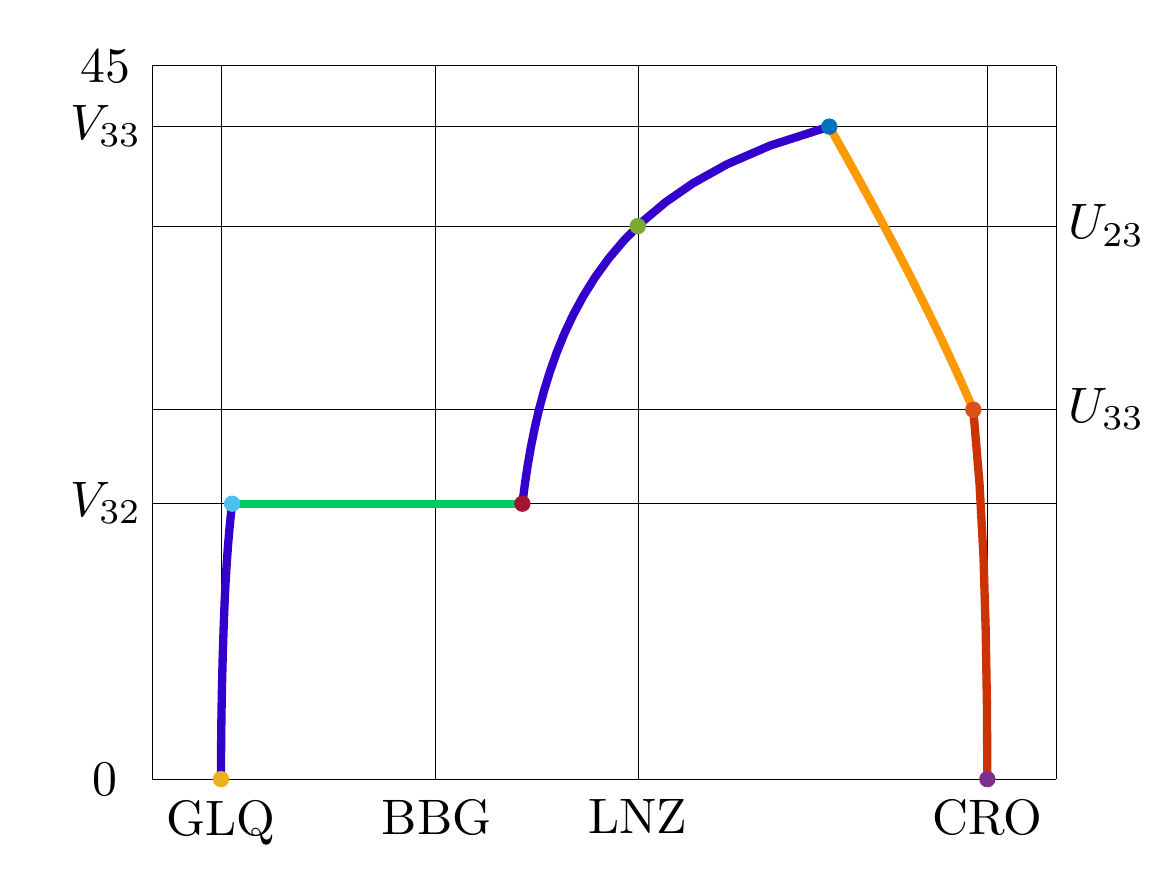}
\includegraphics[width=5.4cm]{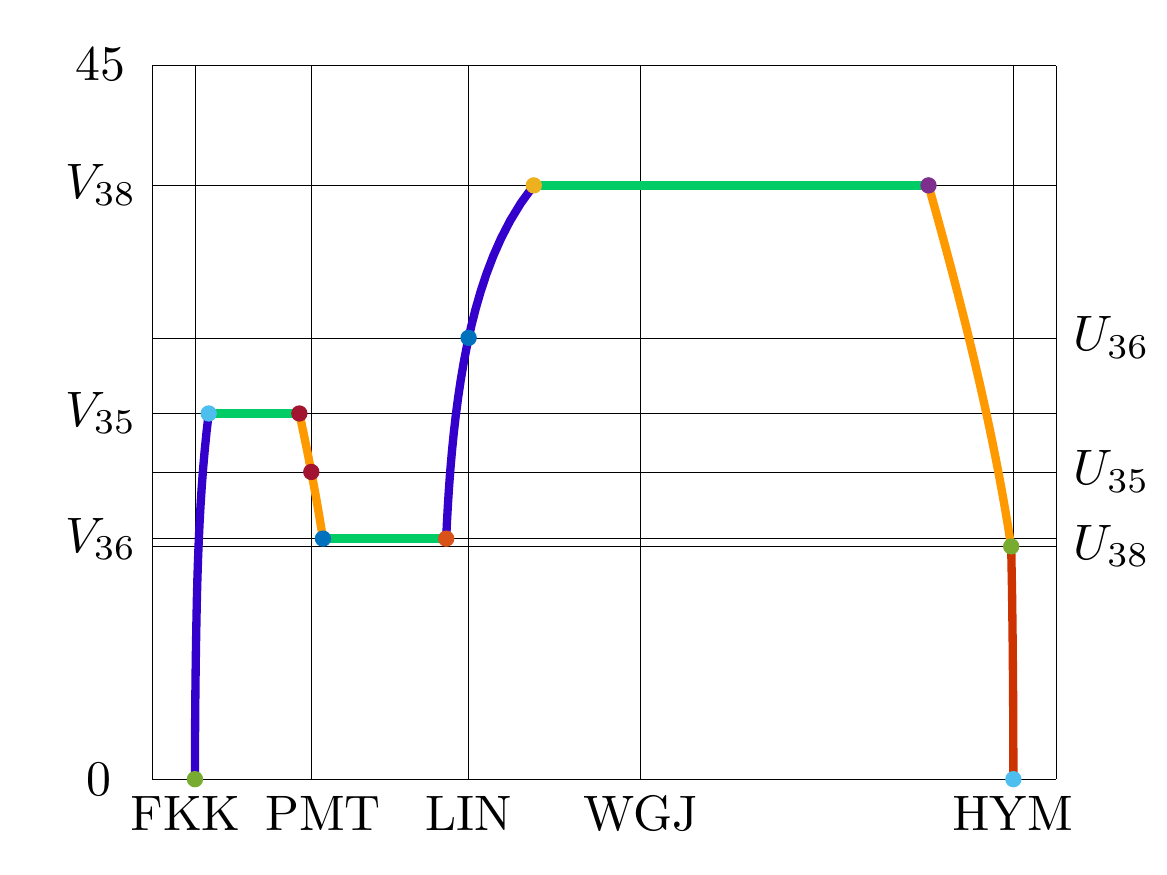}
\includegraphics[width=5.4cm]{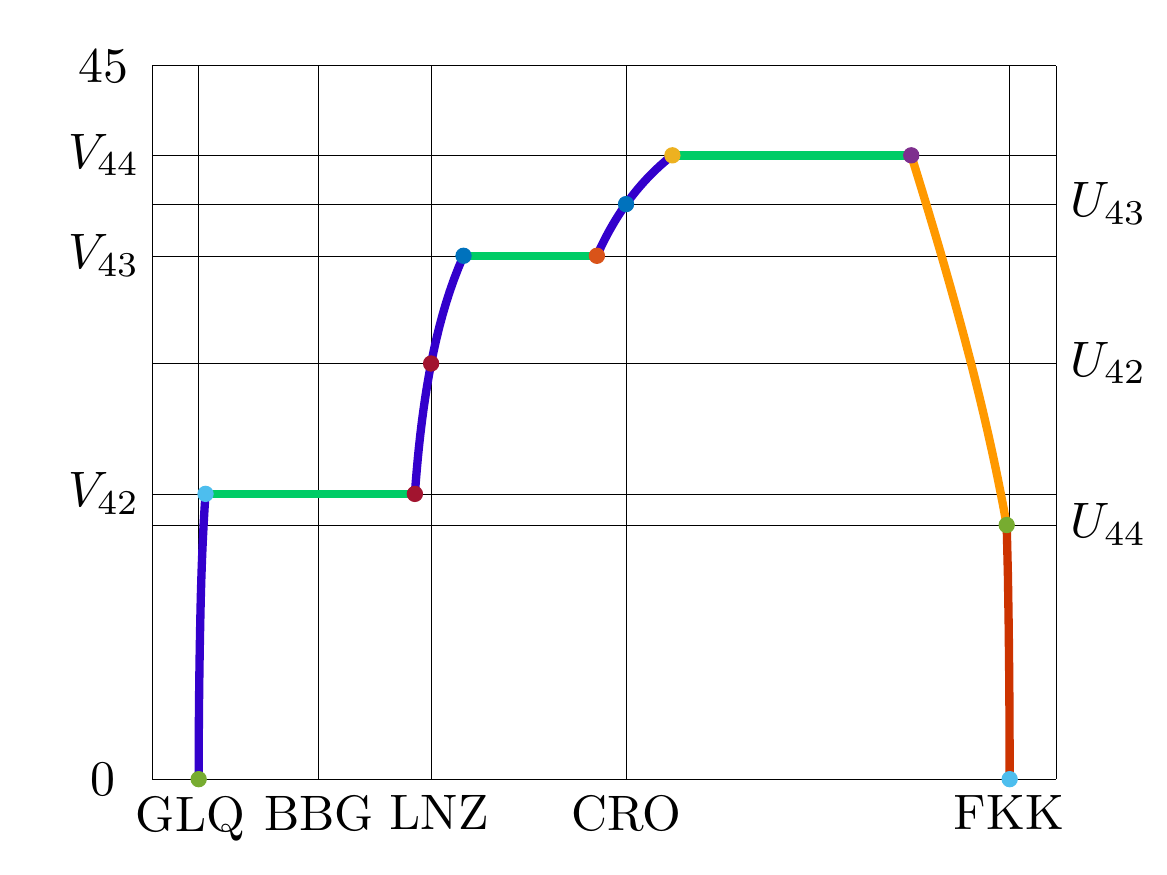}
\end{center}
\caption{\boldmath \small \bf Case Study~\ref{cs3}.  Optimal strategies for ${\mathfrak T}_3$ from GLQ to CRO (left) and from FKK to HYM (centre) and for ${\mathfrak T}_4$ from LIN to HYM (right).  The horizontal scales are different on each graph.  Distances are $18350$ m from GLQ to CRO, $38190$~m from FKK to HYM and $25420$~m from LIN to HYM.}
\label{fig7}
\end{figure}

{\bf \boldmath Train ${\mathfrak T}_4$:}  The strategies for ${\mathfrak T}_4$ are similar to those for ${\mathfrak T}_2$.  For the first segment from GLQ to FKK the optimal driving speeds are $V_{41} = V_{42} \approx 17.8532$, $V_{43} \approx 33.0054$ and $V_{44} \approx 39.4441$~ms$^{-1}$.  We calculate $U_{42} = U_s(V_{42},V_{43}) \approx 26.1578$, $U_{43} = U_s(V_{43},V_{44}) \approx 36.3180$ and $U_{44} = U_b(V_{44}) \approx 16.1006$~ms$^{-1}$.  The cost is $J_{4,04} \approx 5785.1358$ {J kg}$^{-1}$.  The speed profile is shown in Figure~\ref{fig7}.  The strategy from FKK to PMT is a rapid-transit strategy with $U_{45} \approx 13.3900$ and $V_{45} \approx 32.7408$~ms$^{-1}$.  The cost is $J_{4,45} \approx 882.2498$ {J kg}$^{-1}$.  The strategy from PMT to LIN, is a long-haul strategy with $V_{46} \approx 27.5398$ and $U_{46} = U_b(V_{46}) \approx 8.0677$~ms$^{-1}$.  The cost is $J_{4,56} \approx 1099.2413$ {J kg}$^{-1}$.  The strategy from LIN to HYM is also a long-haul strategy with $V_{47} = V_{48} \approx 35.9119$ and $U_{48} = U_b(V_{48}) \approx 13.5989$~ms$^{-1}$.  The cost is $J_{4,68} \approx 4379.5172$ {J kg}$^{-1}$.  The optimal strategy from HYM to EDB is a rapid-transit strategy with $U_{49} \approx 9.7601$ and $V_{49} \approx 24.4923$~ms$^{-1}$.  The cost is $J_{4,89} \approx 406.9617$ {J kg}$^{-1}$.  The total cost for ${\mathfrak T}_4$ is
\begin{equation}
\label{t4cost1}
J_4 = J_{4,04} + J_{4,45} + J_{4,56} + J_{4,68} + J_{4.89} = 12553.1058\ \mbox{{J kg}}^{-1}.
\end{equation}

This means that the overall cost for a realistic implementation of the weighted optimal constant-speed schedule from Case Study~\ref{cs2} is
\begin{equation}
\label{overallcost1}
J = J_1 + J_2 + J_3 + J_4 \approx 50772.0023\ \mbox{{J kg}}^{-1}.
\end{equation} 
In adding the costs we have once again assumed that the masses of the trains are all equal. $\hfill \Box$
\end{cs}

\begin{remark}
\label{cs3r}
Case Study~\ref{cs3} shows that the optimal weighted constant-speed timetable from Case Study~\ref{cs2} provides a feasible schedule for the optimal realistic strategies.  $\hfill \Box$
\end{remark}

\begin{cs}
\label{cs4}

We can use the optimal driving speeds calculated in Case Study~\ref{cs3} to find an improved timetable for the realistic strategies.  The journey times used in Case Study~\ref{cs3} are obtained from the vector $\bfh = [h_j]_{j=1}^{12}$ defined in Table~\ref{tab7} with values given in Table~\ref{tab8}.  We have
$$
\bfh^T = [519, 830, 1409, 1666, 2109, 2971, 2271, 2584, 3857, 3204, 3460, 4720].
$$
The cost gradient vector with respect to $\bfh$ at this point is given by
$$
\bnab J = \left[ \begin{array}{l}
\psi_{13}^{\dag} - \psi_{12} + \psi_{45}^{\dag} - \psi_{44} \\
\psi_{14} - \psi_{13}^{\dag} + \psi_{46} - \psi_{45}^{\dag} \\
\psi_{15} - \psi_{14} + \psi_{23} - \psi_{22} \\
\psi_{16} - \psi_{15} + \psi_{24} - \psi_{23} \\
\psi_{18} - \psi_{16} + \psi_{49}^{\dag} - \psi_{48} \\
\psi^{\dag}_{19} - \psi_{18} + \psi_{28} - \psi_{26} \\
\psi_{25}^{\dag} - \psi_{24} + \psi_{33}^{\dag} - \psi_{32} \\
\psi_{26} - \psi_{25}^{\dag} + \psi_{34}^{\dag} - \psi_{33}^{\dag} \\
\psi_{29} - \psi_{28} + \psi_{38} - \psi_{36} \\
\psi_{35} - \psi_{34}^{\dag} + \psi_{43} - \psi_{42} \\
\psi_{36} - \psi_{35} + \psi_{44} - \psi_{43} \\
\psi_{39}^{\dag} - \psi_{38} + \psi_{48} - \psi^{\dag}_{46} \end{array} \right]
$$
where $\psi_{i,j} = \psi(V_{i,j})$ for the standard entries and $\psi^{\dag}_{i,j} = \varphi(V_{i,j})U_{i,j}/(V_{i,j} - U_{i,j})$ for the non-standard entries corresponding to sections with no speedhold segment where $V_{i,j}$ is a maximum speed rather than an optimal driving speed.  The calculations in Case Study~\ref{cs3} give
$$
\bnab J^T \approx [9.16, -7.12, -2.29, 1.22, 1.61, 0.86, 8.60, 2.38, -0.52, -6.06, 1.34, 11.17 ].
$$
We calculate $\| \bnab J \| \approx 19.79$.  We can now compute an improved timetable by defining a new value of $\bfh$ using the formula $\bfh_{\mbox{\rm \scriptsize new}} = \bfh_{\mbox{\rm \scriptsize old}} - r \bnab J(\bfh_{\mbox{\rm \scriptsize old}})$ where $r = r_0 > 0$ is chosen so that 
$$
J(\bfh_{\mbox{\rm \scriptsize new}}) = \min_{r \geq 0} J(\bfh_{\mbox{\rm \scriptsize old}} - r \bnab J(\bfh_{\mbox{\rm \scriptsize old}})) = J(\bfh_{\mbox{\rm \scriptsize old}} - r_0 \bnab J(\bfh_{\mbox{\rm \scriptsize old}})).
$$
Some preliminary calculations reveal that $r_0 \approx 1.5$.  This gives 
$$
\bfh_{\mbox{\rm \scriptsize new}}^T \approx [505, 841, 1412, 1664, 2107, 2970, 2258, 2580, 3858, 3213, 3458, 4703]
$$
where we have taken the components of $\bfh_{\mbox{\rm \scriptsize new}}$ as the nearest integer to the computed value. We substitute $\bfh_{\mbox{\rm \scriptsize new}}$ into Table~\ref{tab7} to obtain the improved timetable shown in Table~\ref{tab9}.

\begin{table}[htb]
\begin{center}
\begin{tabular}{|c|c|c|c|c|} \hline
Station & ${\mathfrak T}_1$ & ${\mathfrak T}_2$ & ${\mathfrak T}_3$ & ${\mathfrak T}_4$ \\ \hline
GLQ & $t_{1,0} = 0000$ & $t_{2,0} = 0900$ & $t_{3,0} = 1800$ & $t_{4,0} = 2700$ \\ \hline
BBG & \textemdash & \textemdash & \textemdash & \textemdash \\ \hline
LNZ & ($t_{1,2} = 0505)$ & $(t_{2,2} = 1472)$ & $(t_{3,2} = 2318)$ & $(t_{4,2} = 3273)$ \\ \hline
CRO & $t_{1,3} = 0841$ & $(t_{2,3} = 1724)$ & $t_{3,3} = 2640$ & $(t_{4,3} = 3518)$ \\ \hline
FKK & $t_{1,4} = 1412$ & $t_{2,4} = 2258$ & $t_{3,4} = 3213$ & $t_{4,4} = 4045$ \\ \hline
PMT & $(t_{1,5} = 1664)$ & $t_{2,5} = 2580$ & $(t_{3,5} = 3458)$ & $t_{4,5} = 4381$ \\ \hline
LIN & $(t_{1,6} = 2107)$ & $t_{2,6} = 3030$ & $(t_{3,6} = 3918)$ & $t_{4,6} = 4763$ \\ \hline
WGJ & \textemdash & \textemdash & \textemdash & \textemdash \\ \hline
HYM & $t_{1,8} = 2970$ & $t_{2,8} = 3858$ & $t_{3,8} = 4703$ & $t_{4,8} = 5647$ \\ \hline
EDB & $t_{1,9} = 3180$ & $t_{2,9} = 4140$ & $t_{3,9} = 4860$ & $t_{4,9} = 5820$ \\ \hline
\end{tabular}
\end{center}
\vspace{0.4cm}
\caption{\boldmath \small \bf Case Study~\ref{cs4}.  Scheduled departure times for the improved timetable.}
\label{tab9}
\end{table}
    
We now calculate revised optimal strategies for the new schedule.  The calculations are similar to the previous calculations but not identical.  The calculations are summarised below.

{\bf Train \boldmath ${\mathfrak T}_1$:}  From GLQ to CRO the strategy takes the same form as the original strategy.  We have $V_{11} = V_{12} \approx 19.4930$, $V_{13} \approx 39.8883$, $U_{13} \approx 17.2814$~ms$^{-1}$.  We calculate $U_{12} \approx 31.5081$~ms$^{-1}$ and $J_{1,03} \approx 2923.8157$ {J kg}$^{-1}$.  From CRO to FKK the optimal strategy is a rapid-transit strategy defined by $U_{14} \approx 19.0041$ and $V_{14} \approx 43.0894$~ms$^{-1}$ with  $J_{1,34} \approx 3051.2727$ {J kg}$^{-1}$.  From FKK to HYM the strategy is similar to the original strategy with $V_{15} \approx 23.4916$, $V_{16} \approx 15.9525$ and $V_{17} = V_{18} \approx 34.4778$~ms$^{-1}$.  We calculate $U_{15} = U_s(V_{15},V_{16}) \approx 19.9525$ and $U_{16} = U_s(V_{16},V_{17}) \approx 26.3131$~ms$^{-1}$ and $J_{1,48} \approx 6071.8661$ {J kg}$^{-1}$.  From HYM to EDB the optimal strategy is a long-haul strategy with $V_{19} \approx 18.8740$, $U_{19} = U_b(V_{19}) \approx 3.4440$~ms$^{-1}$ and $J_{1,89} \approx 358.9467$ {J kg}$^{-1}$.  Total cost for ${\mathfrak T}_1$ is $J_1 \approx 12405.9012$ {J kg}$^{-1}$.

{\bf Train \boldmath ${\mathfrak T}_2$:}  From GLQ to FKK the strategy takes the same form as the original strategy with $V_{21} = V_{22} \approx 17.5539$, $V_{23} \approx 33.6036$ and $V_{24} \approx 39.9312$~ms$^{-1}$.  We calculate $U_{22} = U_s(V_{22},V_{23}) \approx 26.3195$, $U_{23} = U_s(V_{23},V_{24}) \approx 36.8562$, $U_{24} = U_b(V_{24}) \approx 16.4505$~ms$^{-1}$ and $J_{2,04} \approx 5819.3861$ {J kg}$^{-1}$.  The remaining strategies are all long-haul strategies. From FKK to PMT we have $V_{25} \approx 31.8910$~ms$^{-1}$ with $U_{25} = U_b(V_{25}) \approx  10.8543$~ms$^{-1}$ and $J_{2,45} \approx 847.1220$ {J kg}$^{-1}$.  From PMT to LIN we have $V_{26} \approx 23.3843$, $U_{26} = U_b(V_{26}) \approx 5.6630$~ms$^{-1}$ and $J_{2,56} \approx 1044.3267$ {J kg}$^{-1}$.  From LIN to HYM we have $V_{27} = V_{28} \approx 38.5875$~ms$^{-1}$ and with $U_{28} = U_b(V_{28}) \approx 15.4879$~ms$^{-1}$ and $J_{2,68} \approx 4554.2779$ {J kg}$^{-1}$.  From HYM to EDB we have $V_{19} \approx 11.5225$ and $U_{29} = U_b(V_{29}) \approx 0.9814$~ms$^{-1}$.  The cost is $J_{2,89} \approx 338.7075$ {J kg}$^{-1}$.  The total cost for ${\mathfrak T}_2$ is $J_2 \approx 12603.8202$ {J kg}$^{-1}$.

{\bf Train \boldmath ${\mathfrak T}_3$:} From GLQ to CRO the strategy takes the same form as the original strategy.  We have $V_{31} = V_{32} \approx 18.3070$, $V_{33} \approx 40.6212$ and $U_{33} \approx 20.6766$~ms$^{-1}$.  We calculate $U_{32} = U_s^{\dag}(U_{33},V_{31},V_{32}) \approx 33.5071$~ms$^{-1}$ and $J_{3,03} \approx 3011.0629$ {J kg}$^{-1}$.  From CRO to FKK the optimal strategy is a rapid-transit strategy with $U_{34} \approx 18.4111$ and $V_{34} = 43.0711$~ms$^{-1}$ and $J_{3,34} \approx 3038.1430$ {J kg}$^{-1}$.  From FKK to HYM the strategy is similar to the original strategy with $V_{35} \approx 24.4651$, $V_{36} = 14.8277$ and $V_{37} = V_{38} \approx 38.8945$~ms$^{-1}$.  We calculate $U_{35} = U_s(V_{35},V_{36}) \approx 20.0243$ and $U_{36} = U_s(V_{36},V_{37}) \approx 28.6040$~ms$^{-1}$ and $J_{3,48} \approx 6391.4530$ {J kg}$^{-1}$.  From HYM to EDB the optimal strategy is a rapid-transit strategy with $U_{39} \approx 13.5260$ and $V_{39} \approx 25.5566$~ms$^{-1}$ and $J_{3,89} \approx 451.7930$ {J kg}$^{-1}$.  Total cost for ${\mathfrak T}_3$ is $J_3 \approx 12882.4519$ {J kg}$^{-1}$. 

{\bf Train \boldmath ${\mathfrak T}_4$:}  From GLQ to FKK the strategy takes the same form as the original strategy with $V_{41} = V_{42} \approx 17.4604$, $V_{43} \approx 34.7200$ and $V_{44} \approx 40.6221$~ms$^{-1}$.  We calculate $U_{42} = U_s(V_{42},V_{43}) \approx 27.0122$, $U_{43} = U_s(V_{43},V_{44}) \approx 37.7464$ and $U_{44} = U_b(V_{44}) \approx 16.9485$~ms$^{-1}$ and $J_{4,04} \approx 5879.7786$ {J kg}$^{-1}$.  The next three strategies are all long-haul strategies. From FKK to PMT we have $V_{45} \approx 28.4459$~ms$^{-1}$ with $U_{45} = U_b(V_{45}) \approx  8.6284$~ms$^{-1}$ and $J_{4,45} \approx 815.3487$ {J kg}$^{-1}$.  From PMT to LIN we have $V_{46} \approx 31.8700$, $U_{46} = U_b(V_{46}) \approx 10.8403$~ms$^{-1}$ and $J_{4,56} \approx 1159.8224$ {J kg}$^{-1}$.  From LIN to HYM we have $V_{47} = V_{48} \approx   35.0587$~ms$^{-1}$ and we calculate $U_{48} = U_b(V_{48}) \approx 13.0060$~ms$^{-1}$ and $J_{4,68} \approx 4324.0217$ {J kg}$^{-1}$.  From HYM to EDB we have a rapid-transit strategy with $U_{49} \approx 9.3007$ and $V_{49} \approx 24.3787$~ms$^{-1}$ and $J_{4,89} \approx 402.3895$ {J kg}$^{-1}$.  The total cost for ${\mathfrak T}_4$ is $J_2 \approx 12581.3709$ {J kg}$^{-1}$.

The overall cost for all trains using the improved timetable is
\begin{equation}
\label{overallcost2}
J = J_1 + J_2 + J_3 + J_4 \approx 50473.8853\ \mbox{{J kg}}^{-1}
\end{equation}
which is a modest improvement on the original cost.  Once again we have assumed that the masses of the trains are equal.  We could continue to reduce the overall cost by repeated application of the method of steepest descent. $\hfill \Box$
\end{cs}

\begin{remark}
\label{cs4r}
Case Study~\ref{cs4} shows that we can use the elegant cost gradient formula to find an improved clearance time vector $\bfh_{\mbox{\scriptsize new}}$ and a corresponding reduced cost $J = J(\bfh_{\mbox{\scriptsize new}})$.  Although it seems that the major task in this procedure is calculation of optimal speed profiles for each train it is pertinent to realise that the {\em Energymiser}\textsuperscript{\textregistered} system (known as {\em Opti-Conduite} in France: https://www.sncf.com/fr/groupe/newsroom/opticonduite-energie-economisee) is currently used on-board the famed {\em TGV} service operated by {\em SNCF} in France to continually update optimal driving strategies for journeys in excess of $100$ km in a matter of a few seconds.  $\hfill \Box$
\end{remark}

\begin{cs}
\label{cs5}
We will show that trains with similar but not identical performance functions can use the same timetable.  We assume that the alternative train satisfies (\ref{em:vx}) and (\ref{em:tx}) with $H(v) = \min\{ P_0, P_1/v \}$ ms$^{-2}$ where $P_0 = 1.00$ ms$^{-2}$ and $P_1 = 7.00$ m$^2$s$^{-3}$ and $K(v) = \max\{ -Q_0, -Q_1/v \}$ ms$^{-2}$ where $Q_0 = 1.5$ ms$^{-2}$ and $Q_1 = 9.00$ m$^2$s$^{-3}$ and with resistance $r(v) = r_0 + r_1v + r_2v^2$ ms$^{-2}$ where $r_0 = 0.070$ ms$^{-2}$, $r_1 = 0$ s$^{-1}$ and $r_2 = 0.00005$ m$^{-1}$.  We use the schedule in Case Study~\ref{cs4}.

On the segment from FKK to HYM train ${\mathfrak T}_{1,a}$ uses a similar strategy to ${\mathfrak T}_1$ with $V_{15} \approx 23.6161$, $V_{16} \approx 15.7257$, $V_{17} = V_{18} \approx 34.5982$,  $U_{15} = U_s(V_{22},V_{23}) \approx 19.9346$, $U_{16} = U_s(V_{16},V_{17}) \approx 26.3416$ and $U_{18} = U_b(V_{18}) \approx 16.5956$~ms$^{-1}$ and with $J_{1,a,48} \approx 4423.5502$ {J kg}$^{-1}$.  

Train ${\mathfrak T}_{2,a}$ uses a similar strategy to train ${\mathfrak T}_2$ from GLQ to FKK with $V_{21} = V_{22} \approx 17.4394$, $V_{23} \approx 33.6715$, $V_{24} \approx 40.0411$, $U_{22} = U_s(V_{22},V_{23}) \approx 26.4147$, $U_{23} = U_s(V_{23},V_{24}) \approx 36.9481$ and $U_{24} = U_b(V_{24}) \approx 20.6759$~ms$^{-1}$ and with $J_{2,a,14} \approx 4379.3064$ {J kg}$^{-1}$.

The optimal strategy for ${\mathfrak T}_{3,a}$ from FKK to HYM is similar to the strategy for ${\mathfrak T}_3$ with $V_{35} \approx 24.7406$, $V_{36} \approx 14.3727$ and $V_{37} = V_{38} \approx 39.1359$ ms$^{-1}$ with $U_{32} = U_s(V_{32},V_{33}) \approx 20.0147$, $U_{33} = U_s(V_{33},V_{34}) \approx 28.6643$ and $U_{24} = U_b(V_{34}) \approx 19.9975$~ms$^{-1}$ and with $J_{3,a,48} \approx 4815.7082$ {J kg}$^{-1}$.  

\begin{figure}[htb]
\begin{center}
\includegraphics[width=5.4cm]{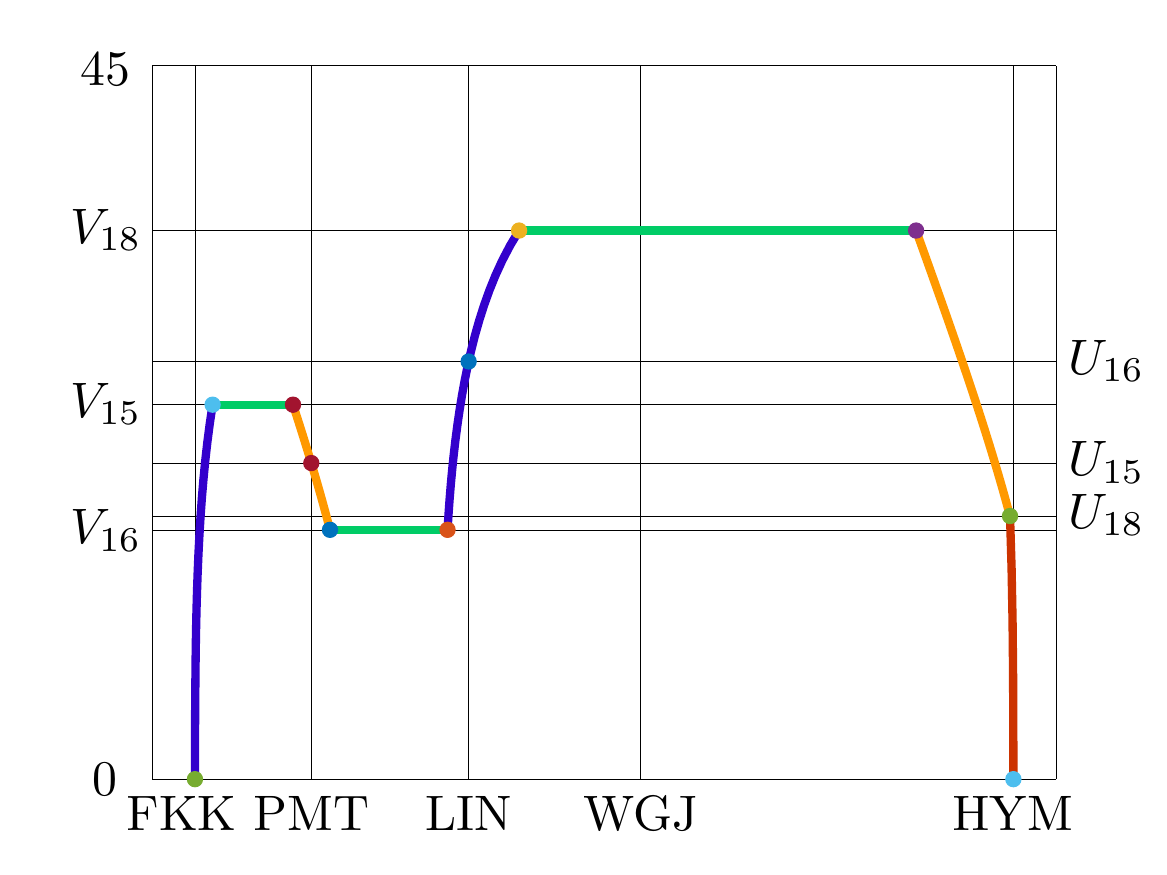}
\includegraphics[width=5.4cm]{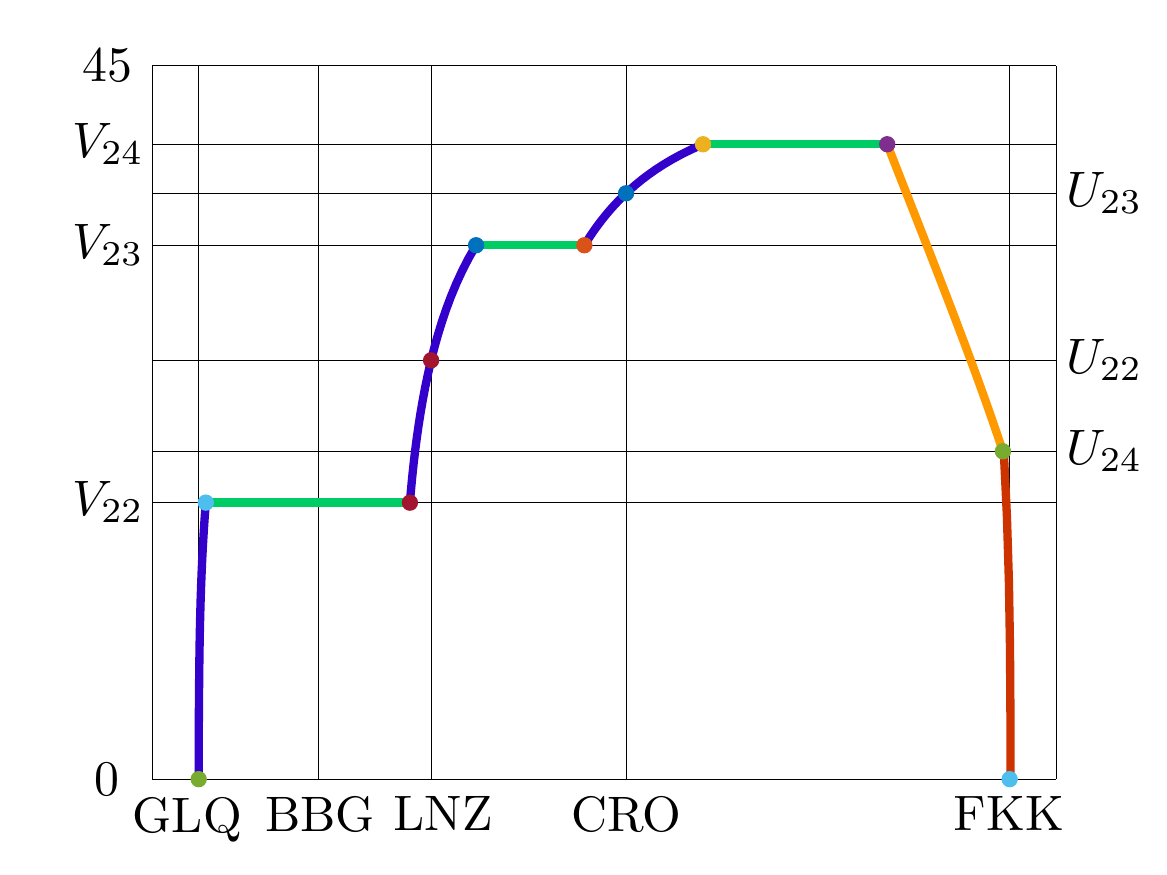}
\includegraphics[width=5.4cm]{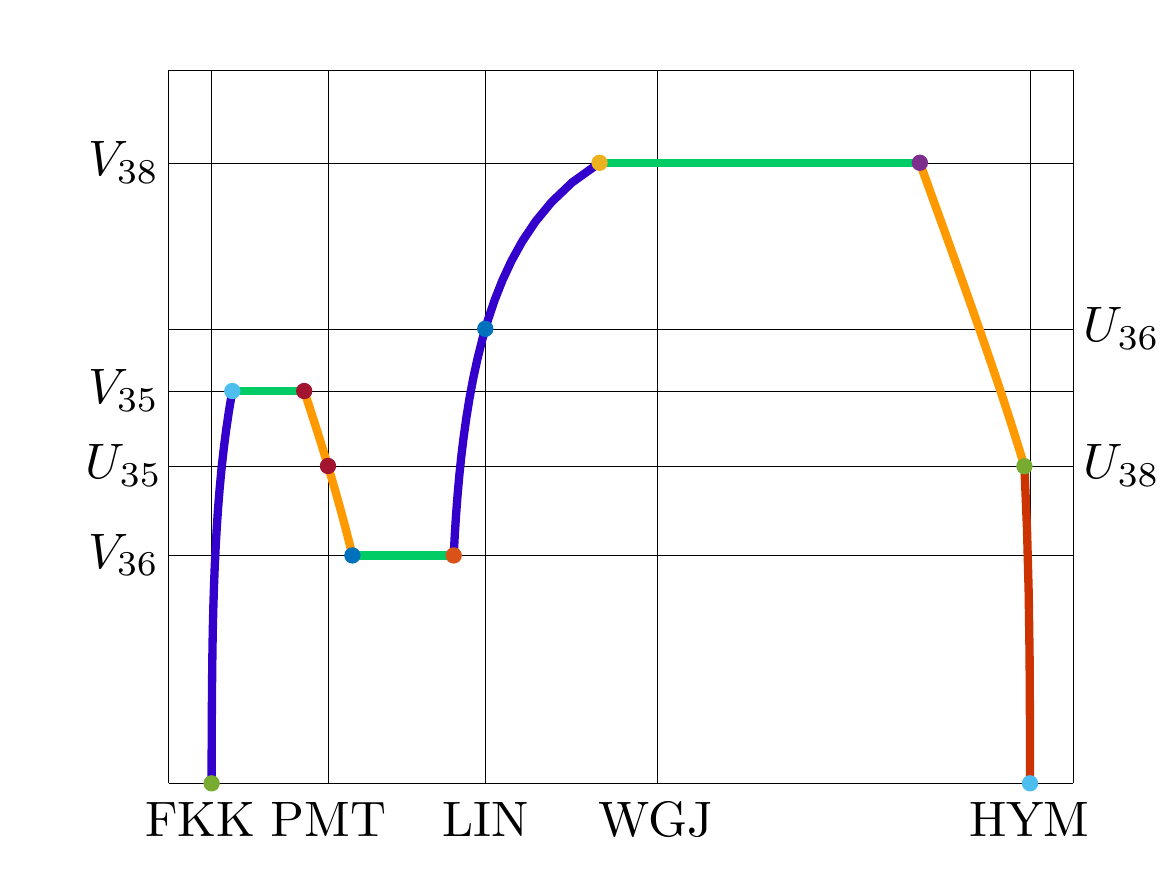}
\end{center}
\caption{\boldmath \small \bf Case Study~\ref{cs5}.  Optimal strategies for ${\mathfrak T}_1$ from FKK to HYM (left), for ${\mathfrak T}_2$ from GLQ to FKK (centre) and for ${\mathfrak T}_3$ from FKK to HYM using the improved timetable in Case Study~\ref{cs4} and an alternative train model.  The horizontal scale on the centre graph is different but the horizontal scale of the other graphs is the same.  Distances are $38190$~m from FKK to HYM and $34820$~m from GLQ to FKK.}
\label{fig8}
\end{figure}
  
The speed profiles are shown in Figure~\ref{fig8}.  Although the profiles are similar to the corresponding profiles for ${\mathfrak T}_1$, ${\mathfrak T}_2$ and ${\mathfrak T}_3$ the cost for ${\mathfrak T}_{1,a}$, ${\mathfrak T}_{2,a}$ and ${\mathfrak T}_{3,a}$ is significantly less because the static resistance is lower.  The key observed impact of the static resistance defined by $r_0$ is the effect on cost whereas the key observed impact of $r_2$ is to limit the maximum achievable speed. $\hfill \Box$
\end{cs}

\begin{remark}
\label{cs5r}
Case Study~\ref{cs5} shows that trains with similar but not identical performance functions can drive to the same timetable.  The entire process of timetable development can be implemented in precisely the same way if the trains ${\mathfrak T}_i$ are not identical.  The speed profile and the costs for ${\mathfrak T}_i$ depend only on the performance functions and the scheduled times.  $\hfill \Box$
\end{remark}
  
\section{Train separation with stochastic journey evolution}
\label{s:tssje}

It is important that theoretical strategies can be implemented effectively in practice.  When pre-planned strategies are implemented on real trains there will inevitably be small discrepancies between planned and actual schedules.  This can cause unwanted difficulties.  When safe separation requires coordinated clearance times and ${\mathfrak T}_i$ is slightly behind schedule with $t_{i, j+1} = h_{i, j+1} + \epsilon_{i, j}$ or ${\mathfrak T}_{i+1}$ is slightly ahead of schedule with $t_{i+1, j-1} = h_{i+1, j-1} - \epsilon_{i+1, j}$ where $h_{i,j+1} = h_{i+1,j-1}$ the signal for ${\mathfrak T}_{i+1}$ at $x_{j-1}$ will be yellow and the signal for ${\mathfrak T}_{i+1}$ at $x_j$ will be red.  Consequently ${\mathfrak T}_{i+1}$ must slow down so that it can stop at signal $x_j$ if necessary.  Yellow signals from a minor time violation can cause significant delays that propagate back through the system.  These delays can be largely avoided by specifying adequate buffer times at each signal.  To model unplanned discrepancies we must find a way to allow for normal stochastic variation.  

In Section~\ref{s:ispit} we begin by analysing the relationship between train position and elapsed journey time.  In Section~\ref{s:esm} we propose a model that generates stochastic differences between the actual elapsed time and the scheduled elapsed time.  In Section~\ref{s:asmes} we apply the stochastic model to a case study looking at implementation of an optimal schedule.  

\subsection{Implementation of strategies with prescribed intermediate times}
\label{s:ispit}

Consider an optimal strategy on a section of level track $[x_{j-1},x_j]$ using a phase of maximum acceleration with $u = H(v)$ from speed $v(x_{j-1}) = U_{j-1} = U_s(V_{j-1},V_j)$ to speed $V_j$, followed by a phase of speedhold at speed $V_j$ with $u = r(V_j)$ and a phase of maximum acceleration to speed $v(x_j) = U_j = U_s(V_j,V_{j+1})$ with $u = H(v)$.  The speed $v = v(x)$ can be found as a function of position by solving the differential equation $vdv/dx = u(v) - r(v)$ subject to the given boundary conditions and the appropriate control function. For the first phase of maximum acceleration define
$$
x_a(v) = x_{j-1} + \int_{U_{j-1}}^v wdw/[H(w) - r(w)]
$$
for each $v \in [U_{j-1},V_j)$ and let $a_j = \lim_{v \uparrow V_j} x_a(v)$ be the point where the speed reaches $v = V_j$ and the first phase ends.  For the second phase of maximum acceleration define
$$
x_a(v) = x_j - \int_v^{U_j} wdw/[H(w) - r(w)]
$$
for each $v \in (V_j,U_j]$ and let $b_j = \lim_{v \downarrow V_j} x_a(v)$ be the point where the second phase begins.  Now the speed function $v = v(x)$ is given by
$$
v(x) = \left\{ \begin{array}{ll}
x_a^{-1}(x) & \mbox{for}\ x \in [x_{j-1}, \xi_{j,1}) \cup (\xi_{j,2}, x_j] \\
V_j & \mbox{for}\ x \in [a_j, b_j]. \end{array} \right.
$$
The above calculations show that we can find the speed $v = v(x)$ from the position $x$ with no knowledge of the elapsed journey time $t$.  Nevertheless $x$ and $t$ are related by
\begin{equation}
\label{tet}
t = t_j + \int_{x_{j-1}}^x d\xi/v(\xi)
\end{equation}
for all $x \in [x_{j-1},x_j]$.  If $x$ is known we can determine $t = t(x)$ directly from (\ref{tet}).  If $t$ is known we can find $x = x(t)$ by solving (\ref{tet}).  The ideas underlying the above discussion remain true for all sections of the journey and all combinations of optimal controls.

\subsection{An elementary stochastic model}
\label{s:esm}

The {\em Energymiser}\textsuperscript{\textregistered} system \cite[Section 1.2]{alb9} and other modern Driver Advisory Systems (DAS) assist drivers to follow an optimal speed profile by displaying the current position $x \in [x_{j-1},x_j]$, the recommended speed $v = v(x)$ and the recommended applied acceleration $u = u(x)$.  In practice, for various reasons, implementation errors will occur and the actual speed profile will differ slightly from the planned speed profile with $v_{\theta}(x,\omega) = v(x) + \delta_{\theta}(x, \omega)$ where the error depends on some scale parameter $\theta \in {\mathbb R}$ and the outcome $\omega \in \Omega$ of a random process.  The driver will always try to correct these errors and there will typically be a succession of intervals $[c,d]$ with $v = v(c)$ at $x=c$ and $v=v(d)$ at $x=d$.  However we may also have
\begin{equation}
\label{cejt}
t_{\theta}(d,\omega) - t_{\theta}(c, \omega) = \int_c^d dx/v_{\theta}(x,\omega) \neq \int_c^d dx/v(x) = t(d) - t(c)
\end{equation}  
Thus the train speed will be correct at $x=c$ and $x=d$ but the elapsed travel time may be wrong.  This means we could model the errors by assuming that the driver implements the correct speed at the wrong time.  Hence we propose a model where the elapsed time function $t_{\theta}(t,\omega) = t + \epsilon_{\theta}(t,\omega)$ is a random walk with drift defined by the stochastic differential equation
\begin{equation}
\label{rwdrift}
dt_{\theta} = dt + \theta dW(t)
\end{equation}
where $t \in [t_{j-1}, t_j]$ is the true time, $[t_{j-1}, t_j]$ is the prescribed time interval and $W(t) = W(t,\omega)$ is a standard Wiener process.  We have followed convention in (\ref{rwdrift}) and suppressed the dependence on $\omega \in \Omega$.  The constant $\theta > 0$ is determined by the observed standard deviation in journey times.  This constant may depend on the train, the train driver, the DAS or any other relevant factor.  The equation (\ref{rwdrift}) has an analytic solution given by
\begin{equation}
\label{rwdriftsol}
t_{\theta}(t) = t_{\theta, j-1} + t - t_{j-1} + \theta \int_{s=t_0}^t dW(s).
\end{equation}
for all $t \in [t_{j-1},t_j]$.  In this model the equations of motion (\ref{em:vx}) and (\ref{em:tx}) do not change and the solutions to these equations do not change.   If we want to match a particular outcome $t_{\theta}(t, \omega)$ to a position then we simply solve
\begin{equation}
t_{\theta}(t,\omega) = t_{\theta}(t_{j-1},\omega) + \int_{x_{j-1}}^{x_{\theta}} d\xi/v(\xi)
\end{equation}
to find $x_{\theta} = x_{\theta}(t,\omega) = x[t_{\theta}(t,\omega)]$.  Now we can calculate the corresponding speed $v_{\theta} = v_{\theta}(t, \omega) = v[x_{\theta}(t,\omega)]$.  Each realization of the random process generates a different elapsed journey time $t_{\theta}(t,\omega) - t_{\theta}(t_{j-1}, \omega)$.  In general $t_{\theta} \neq t$ with $x_{\theta} \neq x$ and $v_{\theta} \neq v$.  The essence of the disparity is that the driver implements the planned speed $v(x)$ at the right point $x$ but at the wrong time $t_{\theta}$.  From (\ref{rwdrift}) it is known that the random variable $\rho = \theta \cdot z = (t_{\theta}(t) - t_{\theta}(t_{j-1})) - (t - t_{j-1})$ is generated by a normal probability density
\begin{equation}
\label{npd}
f(\rho) = (1/\theta \sqrt{2 \pi t}) \exp[ - \rho^2/(\theta^2 t)].
\end{equation}
with mean $0$ and standard deviation $\theta \cdot t^{1/2}$.  The times generated by (\ref{rwdriftsol}) follow a Wiener distribution.  The Wiener process $W(t) = W(t, \omega)$ can be simulated by a random walk with small discrete time steps defined by the random function
\begin{equation}
\label{simwd}
S_n(t, \omega) = \left( 1/\sqrt{n} \right) \sum_{1 \leq k \leq \lfloor nt \rfloor} \xi_k(\omega)
\end{equation}
where $\{ \xi_k\}_{k \in {\mathbb N}}$ are independent identically distributed normal random variables with mean $0$ and variance $1$ and $n \in {\mathbb N}$ is large.  Some basic facts about the Wiener process are outlined in Appendix~\ref{s:dip}.

\subsection{Application of the stochastic model to an existing schedule}
\label{s:asmes}

We will use the stochastic model to examine implementation of the Schedule in Case Study~\ref{cs4}. 

\begin{cs}
\label{cs6}
We assume that an on-board DAS provides continually updated driving advice to train drivers to encourage energy-efficient driving strategies and assist in on-time arrival.  For this reason, in normal operation, we expect variations to section traversal times to be relatively small.  For the sake of argument we assume that on the GLQ to EDB service the observed standard deviation for a scheduled journey time $T = 3120$~s is $s_T \approx 30$~s.  If $t + \epsilon(t)$ denotes the actual time that train ${\mathfrak T}$ reaches position $x$ when the scheduled arrival time is $t = t(x)$ and if the error $\epsilon(t)$ is a Wiener process with scale parameter $\theta$ then the increments $\epsilon(t+\Delta t) - \epsilon(t)$ are independently distributed normal random variables with mean $\mu_{\epsilon} = 0$ and standard deviation $\sigma_{\epsilon} = \theta (\Delta t)^{1/2}$.  We write $\epsilon(t+\Delta t) - \epsilon(t) = \theta (\Delta t)^{1/2} \cdot z \sim {\mathcal N}(0, \theta^2 \Delta t)$.  In our simulation $\theta = s_T/T^{1/2} \approx 0.5371$.

We will use the timetable obtained in Case Study~\ref{cs4}.  The general form of the timetable is displayed in Table~\ref{tab7} with departure times defined by the parameter $\bfh = [h_1,\ldots,h_{12}]$ with value
$$
\bfh = [505, 841, 1412, 1664, 2107, 2970, 2258, 2580, 3858, 3213, 3458, 4703].
$$
The scheduled departure times are displayed in Table~\ref{tab9}.  When the timetables are implemented in practice we assume that updated optimal strategies will be calculated using an on-board computer \cite[Section 9, pp 534\textendash 535]{alb5} in order to reach the next target on time.  In general the updated strategies will take the same form as the strategies calculated in Case Studies~\ref{cs3} and \ref{cs4} but the optimal driving speeds and control switching locations may change.  Let $\bfz_i \sim {\mathcal N}(0,1)^7 \in {\mathbb R}^7$  for each $i=1,\ldots,4$ and $\bfz_{1,\star} \sim {\mathcal N}(0,1)^7$ be random vectors with independently generated standard normal components.  The scheduled signal locations are defined by the vector $\bfx_s = [x_0,x_2,x_3,x_4,x_5,x_6,x_8,x_9]$.  The scheduled times for the trains are defined by 
\begin{eqnarray}
\bfh_1 &= & [0, h_1, h_2, h_3, h_4, h_5, h_6, T_1] \label{h1s} \\
\bfh_2 & = & [900, h_3+60, h_4+60, h_7, h_8, h_6+60, h_9, 900+T_2] \label{h2s} \\
\bfh_3 & = & [1800, h_7+60,h_8+60, h_{10}, h_{11}, h_9+60, h_{12}, 1800+T_3] \label{h3s} \\
\bfh_4 & = & [2700, h_{10}+60, h_{11}+60, h_1+3540, h_2+3540,h_{12}+60,h_5+3540,2700+T_4] \label{h4s} \\
\bfh_{1,\star} & = & [h_{1,\star,j}] = \bft_1 + 3600 \cdot \bfone \label{t1rs}
\end{eqnarray}
where $\bfone = [1] \in {\mathbb R}^{1 \times 7}$.  We write $\bfh_i = [h_{i,j}]_{j=1}^8$ and define
\begin{equation}
\label{dhis}
\Delta^{1/2} \bfh_i = \left[ 0, (h_{i,2}-h_{i,1})^{1/2}, (h_{i,3} - h_{i,2})^{1/2}, \ldots, (h_{i,8} - h_{i,7})^{1/2} \right]
\end{equation}
for each $i=1,\ldots,4$.  The actual departure times will be generated according to the formul{\ae}
\begin{equation}
\label{h1a}
\bfh_{i,a} = [h_{i,a,j}]_{j=1}^8 = \bfh_i + \theta \bfz_i \circ \Delta^{1/2} \bfh_i = \bfh_i + \beps_i
\end{equation}
for each $i=1,\ldots,4$ and 
\begin{equation}
\label{h1sa}
\bfh_{1,\star,a} = \bfh_{1,\star} + \theta \bfz_{1,\star} \circ \Delta^{1/2} \bfh_1 = \bfh_{1,\star} + \beps_{1,\star}.
\end{equation}
where $\bfh_i$ for each $i=1,\ldots,4$ and $\bfh_{1,*}$ are the vectors of scheduled departure times defined above and where we remind readers that $\bfx \circ \bfy$ denotes the Hadamard product of the vectors $\bfx, \bfy$. 

We ran $10000$ simulated trials to calculate actual journey times subject to the normal stochastic variations described above.  There were $48$ trials where the separation conditions were violated.  The mean value of the minimum separation over all $10000$ trials was $m_{\mbox{\scriptsize \rm msep}} \approx 32.2513$~s and the standard deviation was $s_{\mbox{\scriptsize \rm msep}} \approx 10.3963$~s.  The minimum observed value was $min_{\mbox{\scriptsize \rm msep}} \approx -26.9082$~s.  The histogram of observed minimum separations (measured in seconds) is shown on the left in Figure~\ref{fig9}.

\begin{figure}[htb]
\begin{center}
\includegraphics[width=8cm]{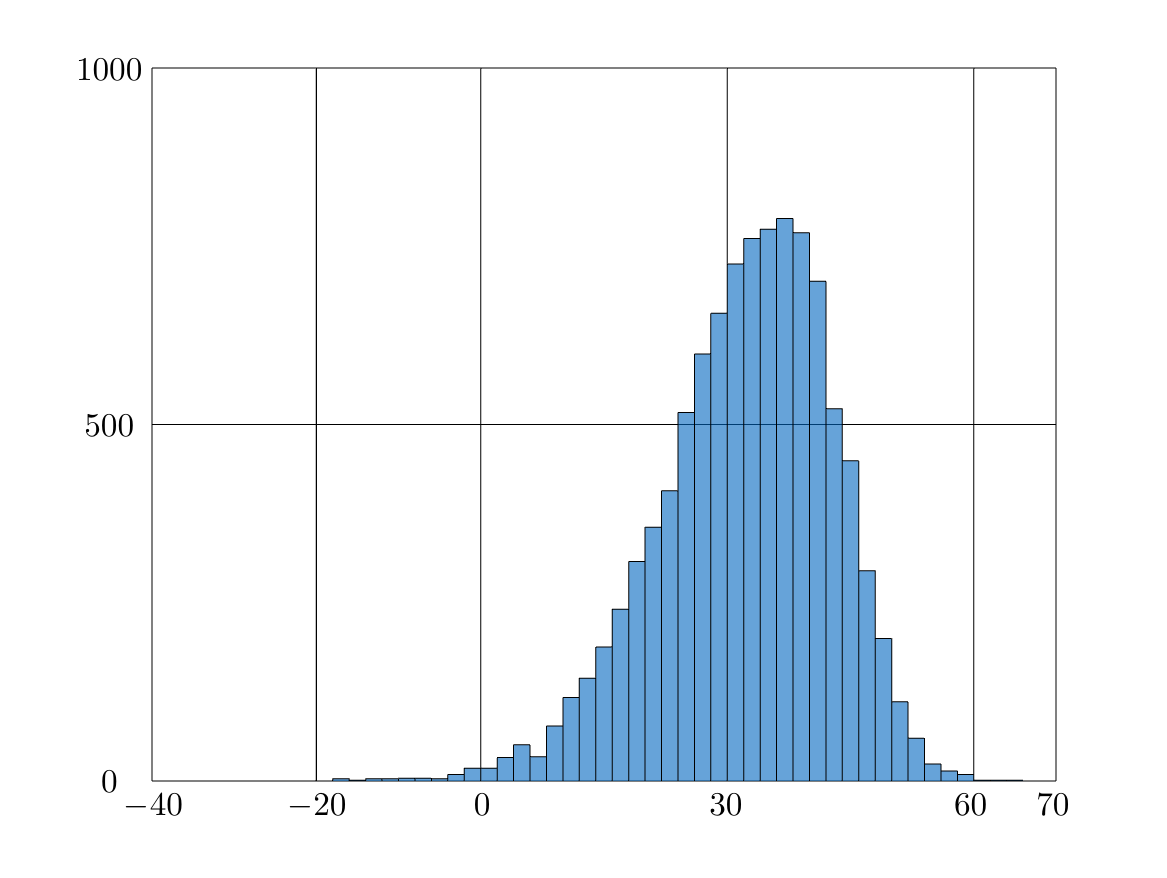}
\includegraphics[width=8cm]{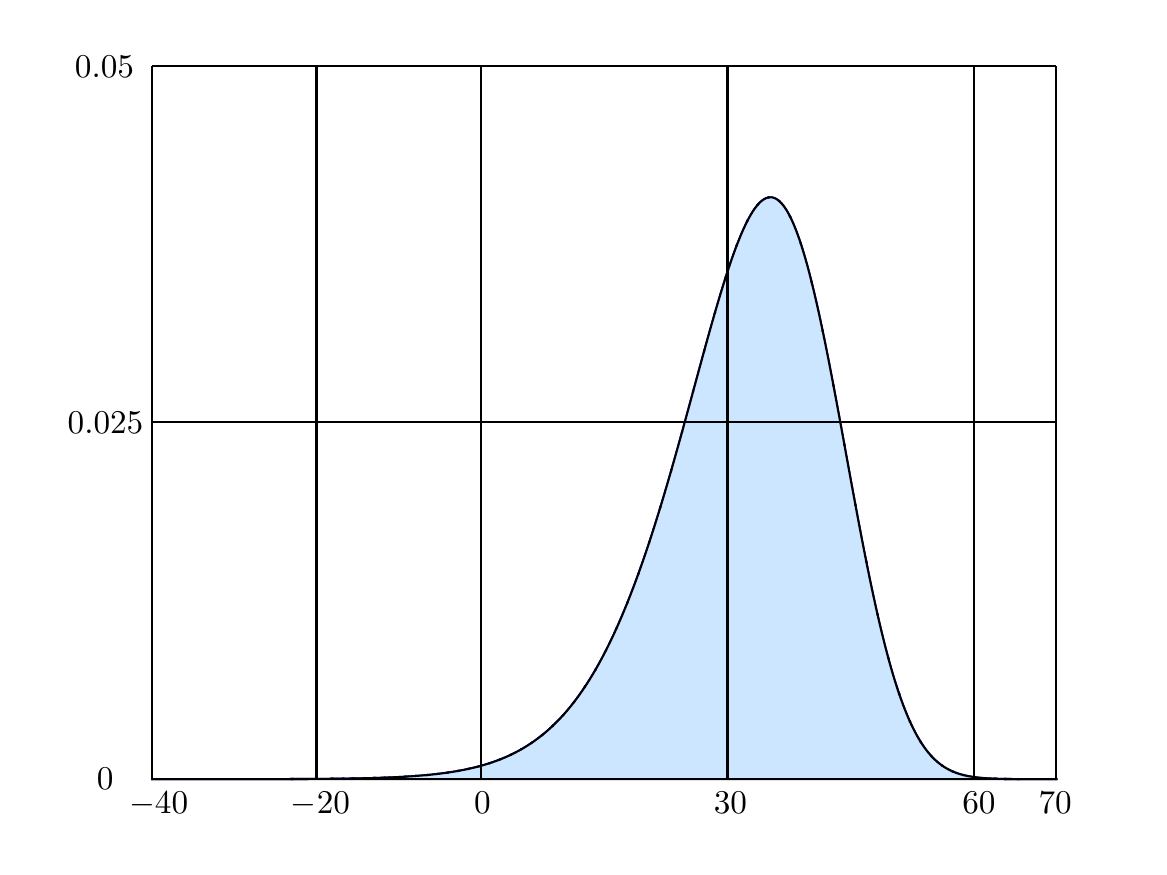}
\end{center}
\caption{\boldmath \small \bf Case Study~\ref{cs6}.  Histogram of simulated minimum separation times from 10000 trials (left) and probability density function for the minimum separation time (right).   In $10000$ trials there were $48$ violations of the required separation conditions due to normal stochastic variation in journey times.  Time is measured in seconds.  The maximum observed violation was approximately $27$~s.  The probability calculated for a violation using the probability density function is ${\mathbb P}[r < 0] \approx 0.0051$.}
\label{fig9}
\end{figure}

\begin{table}[htb]
\begin{center}
\begin{tabular}{|c|c|c|c|c|c|} \hline
Station & ${\mathfrak T}_1$ & ${\mathfrak T}_2$ & ${\mathfrak T}_3$ & ${\mathfrak T}_4$ & ${\mathfrak T}_{1,\star}$ \\ \hline
GLQ & $t_{1,0} = 0000$ & $t_{2,0} = 0900$ & $t_{3,0} = 1800$ & $t_{4,0} = 2700$ & $t_{1,\star,0} = 3600$ \\ \hline
BBG & \textemdash & \textemdash & \textemdash & \textemdash & \textemdash \\ \hline
LNZ & ($t_{1,2} = 0487)$ & $(t_{2,2} = 1472)$ & $(t_{3,2} = 2347)$ & $(t_{4,2} = 3266)$ & $(t_{1,\star,2} = 4106)$ \\ \hline
CRO & $t_{1,3} = 0844$ & $(t_{2,3} = 1716)$ & $t_{3,3} = 2636$ & $(t_{4,3} = 3511)$ & $t_{1,\star,3} = 4454$ \\ \hline
FKK & $t_{1,4} = 1408$ & $t_{2,4} = 2260$ & $t_{3,4} = 3205$ & $t_{4,4} = 4034$ & $t_{1,\star,4} = 5012$ \\ \hline
PMT & $(t_{1,5} = 1666)$ & $t_{2,5} = 2587$ & $(t_{3,5} = 3443)$ & $t_{4,5} = 4390$ & $(t_{1,\star,5} = 5265)$ \\ \hline
LIN & $(t_{1,6} = 2103)$ &{\color{red}$\dag$} $t_{2,6} = 3025$ & $(t_{3,6} = 3923)$ & $t_{4,6} = 4753$ &  $(t_{1,\star,6} = 5688)$ \\ \hline
WGJ & \textemdash & \textemdash & \textemdash & \textemdash & \textemdash \\ \hline
HYM & {\color{red}$\dag$} $t_{1,8} = 3034$ & $t_{2,8} = 3872$ & $t_{3,8} = 4697$ &  $t_{4,8} = 5652$ & $t_{1,\star,8} = 6548$ \\ \hline
EDB & $t_{1,9} = 3182$ & $t_{2,9} = 4149$ & $t_{3,9} = 4861$ & $t_{4,9} = 5818$ & $t_{1,\star,9} = 6772$ \\ \hline
\end{tabular}
\end{center}
\vspace{0.4cm}
\caption{\boldmath \small \bf Case Study~\ref{cs6}.  Trial $\# 107$ departure times.  Trial $\# 107$ was the first trial with a violation of the separation conditions.  The daggers {\color{red}\dag} mark the violated separation condition $t_{1,8} \leq t_{2,6}$.  Train ${\mathfrak T}_2$ would be required to delay departure from LIN until train ${\mathfrak T}_1$ has departed from HYM.  This delay would not impact any other scheduled times.}
\label{tab10}
\end{table}

Trial $\#107$ was the first trial where the separation conditions were violated.  The departure times for this trial are shown in Table~\ref{tab11}.  We also include the results of Trial $\# 1$ which is typical of a trial with no violations of the separation conditions.  The minimum separation in Trial $\# 1$ was $r \approx 34$~s. 

\begin{table}[htb]
\begin{center}
\begin{tabular}{|c|c|c|c|c|c|} \hline
Station & ${\mathfrak T}_1$ & ${\mathfrak T}_2$ & ${\mathfrak T}_3$ & ${\mathfrak T}_4$ & ${\mathfrak T}_{1,\star}$ \\ \hline
GLQ & $t_{1,0} = 0000$ & $t_{2,0} = 0900$ & $t_{3,0} = 1800$ & $t_{4,0} = 2700$ & $t_{1,\star,0} = 3600$ \\ \hline
BBG & \textemdash & \textemdash & \textemdash & \textemdash & \textemdash \\ \hline
LNZ & ($t_{1,2} = 0500)$ & $(t_{2,2} = 1501)$ & $(t_{3,2} = 2318)$ & {\color{blue} $\bullet$} $(t_{4,2} = 3265)$ & $(t_{1,\star,2} = 4089)$ \\ \hline
CRO & $t_{1,3} = 0831$ & $(t_{2,3} = 1721)$ & $t_{3,3} = 2648$ & $(t_{4,3} = 3522)$ & $t_{1,\star,3} = 4449$ \\ \hline
FKK & $t_{1,4} = 1416$ &  $t_{2,4} = 2255$ & {\color{blue} $\bullet$} $t_{3,4} = 3231$ & $t_{4,4} = 4051$ & $t_{1,\star,4} = 5023$ \\ \hline
PMT & $(t_{1,5} = 1665)$ & $t_{2,5} = 2589$ & $(t_{3,5} = 3453)$ & $t_{4,5} = 4382$ & $(t_{1,\star,5} = 5277)$ \\ \hline
LIN & $(t_{1,6} = 2110)$ & $t_{2,6} = 3023$ & $(t_{3,6} = 3907)$ & $t_{4,6} = 4777$ & $(t_{1,\star,6} = 5692)$ \\ \hline
WGJ & \textemdash & \textemdash & \textemdash & \textemdash & \textemdash \\ \hline
HYM & $t_{1,8} = 2955$ & $t_{2,8} = 3860$ & $t_{3,8} = 4717$ & $t_{4,8} = 5645$ & $t_{1,\star,8} = 6583$ \\ \hline
EDB & $t_{1,9} = 3175$ & $t_{2,9} = 4148$ & $t_{3,9} = 4841$ & $t_{4,9} = 5826$ & $t_{1,\star,9} = 6780$ \\ \hline
\end{tabular}
\end{center}
\vspace{0.4cm}
\caption{\boldmath \small \bf Case Study~\ref{cs6}.  Trial $1$.  Typical timetable with no violations of the separation conditions.  The minimum separation is $34$ s.  The bullets {\color{blue}$\bullet$} show where the minimum separation occurs with $t_{4,2} = t_{3,4} + 34$.}
\label{tab11}
\end{table}

We will now show that the observed distribution of minimum separation times can be supported by a theoretical argument.  The minimum separation on each trial is a random variable defined by $r = \min \{r_1, r_2, r_3 \}$ where
\begin{eqnarray*}
r_1 & = & \min \{ t_{2,2}-t_{1,4}, t_{3,2}-t_{2,4}, t_{4,2}-t_{3,4}, t_{1,\star,2}-t_{4,4} \} \\
r_2 & = & \min \{ t_{2,3}-t_{1,5}, t_{3,3}-t_{2,5}, t_{4,3}-t_{3,5}, t_{1,\star,3}-t_{4,5} \} \\
r_3 & = & \min \{ t_{2,6}-t_{1,8}, t_{3,6}-t_{2,8}, t_{4,6}-t_{3,8}, t_{1,\star,6}-t_{4,8} \}
\end{eqnarray*}
and where the individual time differences are normally distributed random variables defined by 
\begin{eqnarray*}
t_{2,2} - t_{1,4} & = & 60 + \theta  \left[ (h_3 + 60 - 900) + (h_3 - h_2) \right]^{1/2} z_{(2,2),(1,4)} \approx 60 + s_{11} z_{(2,2),(1,4)} \\
t_{3,2} - t_{2,4} & = & 60 + \theta  \left[ (h_7 + 60 - 1800) + (h_7 - h_4 - 60) \right]^{1/2} z_{(3,2),(2,4)} \approx 60 + s_{12} z_{(3,2),(2,4)} \\
\vdots & & \vdots \\
t_{4,3} - t_{3,5} & = & 60 + \theta \left[ (h_{11} - h_{10}) + (h_{11} - h_{10}) \right]^{1/2} z_{(4,3),(3,5)} \approx 60 + s_{23} z_{(4,3),(3,5)} \\
\vdots & & \vdots \\
t_{1,\star,6} - t_{4,8} & = & 60 + \theta \left[ (h_5 - h_4) + (h_5 + 3540 - h_{12} - 60) \right]^{1/2} z_{(1,\star,6),(4,8)} \approx 60 + s_{34} z_{(1,\star,6),(4,8)}
\end{eqnarray*}
in which $z_{(2,2),(1,4)} , z_{(3,2),(2,4)},\ldots, z_{(4,3),(3,5)} \ldots z_{(1,\star,6),(4,8)} \sim {\mathcal N}(0,1)$ are standard normal random variables.  The constants $S = [s_{i,j}]$ are given by
$$
S \approx \left[ \begin{array}{cccc}
18.1580 & 17.4201 & 18.1818 & 17.2538 \\
12.0576 & 13.6297 & 11.8889 & 13.9229 \\
19.4615 & 19.2753 & 18.3476 & 19.5650 \end{array} \right].
$$
We can now find a theoretical distribution which is close to the observed frequencies for the minimum separation.  The cumulative distribution function for the standard normal distribution is defined in terms of the standard tabulated error function $\mbox{\rm erf}(x)$ by the formula $\Phi(x) = (1/2) \left[ 1 + \mbox{\rm erf}(x/2^{1/2}) \right]$ for all $x \in {\mathbb R}$.  It follows that the cumulative distribution function for the random variable $r_i$ is
\begin{eqnarray*}
F_i(x) & = & \Phi \left((x - 60)/s_{i1} \right) + \Phi \left((x-60)/s_{i2} \right) \cdot \left[ 1 - \Phi \left((x - 60)/s_{i1} \right) \right] \\
& & \hspace{0.5cm} + \Phi \left((x-60)/s_{i3} \right) \cdot \left[ 1 - \Phi \left((x - 60)/s_{i2} \right) \right] \cdot \left[ 1 - \Phi \left((x - 60)/s_{i1} \right) \right] \\
& & \hspace{0.5cm} + \Phi \left((x-60)/s_{i4} \right)\left[ 1 - \Phi \left((x-60)/s_{i3} \right) \right] \cdot \left[ 1 - \Phi \left((x - 60)/s_{i2} \right) \right] \cdot \left[ 1 - \Phi \left((x - 60)/s_{i1} \right) \right]
\end{eqnarray*}
for each $i = 1,2,3$.  Now the cumulative distribution function for the random variable $r$ is
$$
F(x) = F_1(x) + F_2(x) \cdot \left[ 1 - F_1(x) \right] + F_3(x) \cdot \left[ 1 - F_2(x) \right] \cdot \left[ 1 - F_1(x) \right]
$$
for all $x \in {\mathbb R}$.  We used {\rm M{\sc atlab}} to sketch the probability density function $f(x) = F^{\, \prime}(x)$.  The graph is shown on the right in Figure~\ref{fig9}.  The theoretical mean and standard deviation are given approximately by $\mu_{\mbox{\scriptsize \rm msep}} \approx 32.2897$ and $\sigma_{\mbox{\scriptsize \rm msep}} \approx 10.3397$.  The probability of a violation is ${\mathbb P}[ r < 0] \approx 0.0051$. $\hfill \Box$
\end{cs} 

\begin{remark}
\label{cs6r}
Case Study~\ref{cs6} shows that by choosing a buffer time approximately double the observed standard deviation for the entire journey it is possible to eliminate almost all potential safe-separation violations caused by normal stochastic variation in segment traversal times. $\hfill \Box$
\end{remark}  

\section{Conclusions}
\label{conc}

We found an analytic solution to the problem of minimizing tractive energy consumption for a fleet of similar trains subject to active clearance-time constraints that ensure safe separation and also compress the line-occupancy timespan for the fleet. Compared to the schedule required for safe separation when each train uses the classic single-train optimal strategy our method provides a substantial reduction in line-occupancy timespan with only a minimal increase in journey costs.  We demonstrated our methods using a sequence of case studies for an existing timetable on the Glasgow to Edinburgh line using hypothetical separation constraints.  Our calculations showed that violations of the separation conditions which would occur if all trains use individual optimal strategies can be avoided using our recommended strategies with no change to the scheduled initial departure and final arrival times and only a small increase in costs.  We showed that an optimal timetable for constant-speed strategies with safe separation constraints could be calculated efficiently using a rapidly convergent Newton iteration.  Our algorithm for calculation of the optimal clearance times using the constant-speed strategies is stable and efficient.  Finally we showed that realistic strategies could be implemented using an initial schedule obtained from a weighted optimal constant-speed timetable and that this initial schedule could be improved using a multi-dimensional method of steepest descent to find a near optimal schedule for the realistic strategies.

The calculations for realistic strategies are much more demanding than the calculations for the constant-speed strategies.  However we reiterate that the {\em Energymiser}\textsuperscript{\textregistered} system is already used in practice on very fast trains to continually calculate updated optimal driving strategies every few seconds for journeys of more than one hundred km.  It is envisaged that commercial algorithms could be developed to implement rapid calculation of the schedules devised here.  Significant reductions in line-occupancy timespan for scheduled services and cost reductions of between $5\%$ and $10\%$ obtained by implementing optimal driving strategies can save many millions of dollars each year for large rail organisations. 
\bigskip

\appendix

\section{Appendix: Mathematical background}
\label{mb}

\subsection{Optimal strategies with intermediate time constraints}
\label{s:itc}

A train travels from $x = x_0$ to $x = x_2$ on level track.  The time allowed for the journey is $h_2=T$.  We use the model described in Section~\ref{s:em} to show that the optimal strategy is a long-haul strategy with optimal driving speed $V$.  Now suppose there is an intermediate time constraint and the train must pass through $x_1$ at time $h_1$ without stopping.  If we assume that 
$$
\int_0^V dt_a(v) + (1/V) \left[ x_1 - \int_0^V dx_a(v) \right] < h_1.
$$
the long-haul strategy is no longer feasible.  We will consider two plausible alternatives.  In each case we assume a strategy on $(x_0,x_1)$ of maximum acceleration to speed $v = V_1 < V$, speedhold with $v = V_1$, and maximum acceleration to  $v = U_1$ at $x_1$.  Our first alternative on $(x_1,x_2)$ is a strategy of maximum acceleration starting from $v(x_1) = U_1$ to $v = V_2 > V$, speedhold with $v = V_2$, coast to $v = U_2$ and maximum brake.  If this alternative is feasible then it is optimal.  Our second alternative on $(x_1,x_2)$ is a strategy of maximum acceleration starting from $v(x_1) = U_1$ to a maximum speed $v = V_2 > V$, coast to $v = U_2$ and maximum brake.  The second alternative is more robust but is only optimal if the first strategy is not feasible.

{\bf \boldmath Case $1$:}   We assume the strategy is maximum acceleration, speedhold with $v = V_1$, maximum acceleration to pass through $x_1$ with speed $v(x_1) = U_1$ and reach speed $V_2$, speedhold at speed $V_2$, coast to speed $U_2$ and maximum brake.   The cost of the strategy is
$$
J(\bfU,\bfV) =  \int_0^{V_2} H(v)dx_a(v) + r(V_1) \xi_1(U_1) + r(V_2) \xi_2(U_1,U_2,V_2)
$$ 
where we have defined
\begin{equation}
\label{Del1}
\xi_1(U_1) =  x_1 - x_0 - \int_0^{U_1} dx_a(v)
\end{equation}
and
\begin{equation}
\label{Del2}
\xi_2(U_1,U_2,V_2) = x_2 - x_1 - \int_{U_1}^{V_2} dx_a(v) - \int_{U_2}^{V_2} (-1)dx_c(v) - \int_0^{U_2} (-1)dx_b(v).
\end{equation}
This strategy is valid if $\xi_1(U_1) \geq 0$ and $\xi_2(U_1,U_2,V_2) \geq 0$.  The time taken to traverse $[x_0,x_1]$ is
$$
\tau_1(U_1,V_1) = \int_0^{U_1} dt_a(v) + \xi_1(U_1)/V_1 
$$
and the time taken to traverse $[x_1,x_2]$ is
$$
\tau_2(U_1,U_2,V_2) = \int_{U_1}^{V_2} dt_a(v) + \int_{U_2}^{V_2} (-1)dt_c(v) + \int_0^{U_2} (-1)dt_b(v) + \xi_2(U_1,U_2,V_2)/V_2. 
$$
We wish to minimize $J(\bfU,\bfV)$ subject to $\tau_1(U_1,V_1) - h_1 = 0$ and $\tau_2(U_1,U_2,V_2)-h_2+h_1=0$.  Define
$$
{\mathcal J}(\bfU,\bfV) = J(\bfU,\bfV) + \mu_1[\tau_1(U_1,V_1)-h_1] + \mu_2[\tau_2(U_1,U_2,V_2)-h_2+h_1]
$$
where $\mu_1, \mu_2$ are Lagrange multipliers.  Differentiation and some elementary algebra gives
\begin{eqnarray*}
\frac{\partial {\mathcal J}}{\partial U_1} & = & \frac{U_1[r(V_2)-r(V_1)] + \mu_1(1-U_1/V_1) - \mu_2(1-U_1/V_2)}{H(U_1) - r(U_1)} \\
\frac{\partial {\mathcal J}}{\partial U_2} & = & \left[ \rule{0cm}{0.4cm} U_2r(V_2) - \mu_2(1 - U_2/V_2) \right] \left[ \frac{1}{r(U_2)} - \frac{1}{K(U_2)+r(U_2)} \right] \\
\frac{\partial {\mathcal J}}{\partial V_1} & = & \left[ \rule{0cm}{0.4cm} r^{\, \prime}(V_1) - \mu_1/V_1^2 \right] \xi_1(U_1) \\
\frac{\partial {\mathcal J}}{\partial V_2} & = & \left[ \rule{0cm}{0.4cm} r^{\, \prime}(V_2) - \mu_2/V_2^2 \right] \xi_2(U_1,U_2,V_2). 
\end{eqnarray*}
Setting the third and fourth partial derivatives equal to zero gives $\mu_j = V_j^2 r^{\, \prime}(V_j) = \psi(V_j)$ for each $j=1,2$.  If we set the second partial derivative to zero we get
\begin{equation}
\label{obs}
U_2 = \psi(V_2)/\varphi^{\, \prime}(V_2) = U_b(V_2).
\end{equation}
Finally, if we set the first partial derivative equal to zero we get
\begin{equation}
\label{osps}
U_1 = [\psi(V_2) - \psi(V_1)]/[\varphi^{\, \prime}(V_2) - \varphi^{\, \prime}(V_1)] = U_s(V_1,V_2).
\end{equation}
This strategy of optimal type is defined by two independent variables $V_1$ and $V_2$, and two dependent variables $U_1= U_s(V_1,V_2)$ and $U_2 = U_b(V_2)$.  The values of $V_1$ and $V_2$ are fixed by the signal times $h_1$ and $h_2$. 

{\bf \boldmath Case $2$:}  We assume the strategy is maximum acceleration, speedhold at speed $V_1$, maximum acceleration to pass through $x_1$ with speed $v(x_1)=U_1$ and reach a maximum speed $V_2$, coast to speed $U_2$ and maximum brake to stop at $x_2$.  The cost of the strategy is
$$
J(\bfU,\bfV) =  \int_0^{V_2} H(v)dx_a(v) + r(V_1) \xi_1(U_1).
$$ 
The time taken to traverse $[x_0,x_1]$ is
$$
\tau_1(\bfU,\bfV) = \int_0^{U_1} dt_a(v) + \xi_1(U_1)/V_1 
$$
and the time taken to traverse $[x_1,x_2]$ is
$$
\tau_2(\bfU,\bfV) = \int_{U_1}^{V_2} dt_a(v) + \int_{U_2}^{V_2} (-1)dt_c(v) + \int_0^{U_2} (-1)dt_b(v). 
$$
We wish the minimize $J(\bfU,\bfV)$ subject to the distance constraint $\xi_2(U_1,U_2,V_2) = 0$ and the time constraints $\tau_1(U_1,V_1) - h_1 = 0$ and $\tau_2(U_1,U_2,V_2)-h_2+h_1=0$.  Define
$$
{\mathcal J}(\bfU,\bfV) = J(\bfU,\bfV) + \lambda \xi_2(U_1,U_2,V_2) + \mu_1[\tau_1(U_1,V_1)-h_1] + \mu_2[\tau_2(U_1,U_2,V_2)-h_2+h_1]
$$
where $\lambda, \mu_1, \mu_2$ are Lagrange multipliers.  Differentiation with respect to the four independent variables and some elementary algebra gives
\begin{eqnarray*}
\frac{\partial {\mathcal J}}{\partial U_1} & = & \frac{\mu_1 - \mu_2 -[r(V_1) - \lambda +\mu_1/V_1]U_1}{H(U_1) - r(U_1)} \\
\frac{\partial {\mathcal J}}{\partial U_2} & = & \left[ \rule{0cm}{0.4cm} \lambda U_2 - \mu_2 \right] \left[ \frac{1}{r(U_2)} - \frac{1}{K(U_2)+r(U_2)} \right] \\
\frac{\partial {\mathcal J}}{\partial V_1} & = & \left[ \rule{0cm}{0.4cm} r^{\, \prime}(V_1) - \mu_1/V_1^2 \right] \xi_1(U_1) \\
\frac{\partial {\mathcal J}}{\partial V_2} & = & \left[ \rule{0cm}{0.4cm} \varphi(V_2) - \lambda V_2 + \mu_2 \right] \frac{H(V_2)}{H(V_2) - r(V_2)}. 
\end{eqnarray*}
Setting the partial derivatives equal to zero gives $\mu_1 = \psi(V_1)$, $\lambda = \varphi(V_2)/(V_2-U_2)$, $\mu_2 = U_2 \varphi(V_2)/(V_2-U_2)$ and
\begin{equation}
\label{ospsd}
U_1 = [\psi(V_1) - U_2\varphi(V_2)/(V_2-U_2)]/[\varphi^{\, \prime}(V_1) - \varphi(V_2)/(V_2-U_2)] = U_s^{\dag}(U_2,V_1,V_2).
\end{equation}
The formula for $U_1$ is valid provided $U_2 \in [\psi(V_2)/\varphi^{\, \prime}(V_2), V_2)$.  For $U_2 = \psi(V_2)/\varphi^{\, \prime}(V_2)$ the formula reduces to (\ref{osps}).  The strategy is defined by the three independent variables $U_2$, $V_1$ and $V_2$.  The values are fixed by the times $h_1$ and $h_2$ and the distance $x_2-x_1$. 

\subsection{Cost gradient with respect to journey time for journeys with intermediate time constraints}
\label{s:rccitc}

We outline the derivation of the key formul{\ae} for the partial rate of change of journey cost with respect to the prescribed signal times.  We consider Case $1$ and Case $2$ from Appendix~\ref{s:itc} and use the same notation.

{\bf \boldmath Case $1$:}  We assume the optimal strategy takes the form described in Case $1$ of Appendix~\ref{s:itc}.  The cost of the strategy is
\begin{equation}
\label{jc1}
J = \int_0^{V_2} H(v)dx_a(v) + r(V_1) \xi_1(U_1) + r(V_2) \xi_2(U_1,U_2, V_2)
\end{equation}
where $V_1$ and $V_2$ are independent variables, $U_1 = U_s(V_1,V_2)$ is the optimal speed at $x_1$ given by (\ref{osps}) and $U_2 = U_b(V_2)$ is the optimal braking speed given by (\ref{obs}).  The time constraints are
\begin{equation}
\label{t1c1}
\int_0^{U_1} dt_a(v) + \xi_1(U_1)/V_1 = h_1
\end{equation}
and
\begin{equation}
\label{t2c1}
\int_{U_1}^{V_2} dt_a(v) + \int_{U_2}^{V_2} (-1)dt_c(v) +\int_0^{U_2} (-1)dt_b(v) + \xi_2(V_1,V_2)/V_2 = h_2 - h_1.
\end{equation}
Differentiation of (\ref{t1c1}) with respect to $h_1$ and $h_2$ and rearrangement gives
\begin{equation}
\label{dt1c1dh1}
\xi_1 \frac{\partial V_1}{\partial h_1} = -V_1^2 + V_1(V_1-U_1)t_a^{\, \prime}(U_1) \frac{\partial U_1}{\partial h_1}
\end{equation}
and
\begin{equation}
\label{dt1c1dh2}
\xi_1 \frac{\partial V_1}{\partial h_2} = - V_1(V_1-U_1)t_a^{\, \prime}(U_1) \frac{\partial U_1}{\partial h_2}. 
\end{equation}
Differentiation of (\ref{t2c1}) with respect to $h_1$ and $h_2$ and rearrangement gives
\begin{equation}
\label{dt2c1dh1}
\xi_2 \frac{\partial V_2}{\partial h_1} = V_2^2 - V_2(V_2-U_1)t_a^{\, \prime}(U_1) \frac{\partial U_1}{\partial h_1}
\end{equation}
and
\begin{equation}
\label{dt2c1dh2}
\xi_2 \frac{\partial V_2}{\partial h_2} = -V_2^2 - V_2(V_2-U_1)t_a^{\, \prime}(U_1) \frac{\partial U_1}{\partial h_2} + V_2 (V_2 - U_2)[t_c^{\, \prime}(U_2) - t_b^{\, \prime}(U_2)] \frac{\partial U_2}{\partial h_2}.
\end{equation}
Differentiation of (\ref{jc1}) with respect to $h_1$ and $h_2$ gives
\begin{equation}
\label{djc1dh1}
\frac{\partial J}{\partial h_1} = r^{\, \prime}(V_1)\xi_1 \frac{\partial V_1}{\partial h_1}  - r(V_1) x_a^{\, \prime}(U_1) \frac{\partial U_1}{\partial h_1} + r^{\, \prime}(V_2)\xi_2 \frac{\partial V_2}{\partial h_1} + r(V_2) x_a^{\, \prime}(U_1) \frac{\partial U_1}{\partial h_1} 
\end{equation}
and
\begin{eqnarray}
\label{djc1dh2}
\frac{\partial J}{\partial h_2} & = & r^{\, \prime}(V_1)\xi_1 \frac{\partial V_1}{\partial h_2}  - r(V_1) x_a^{\, \prime}(U_1) \frac{\partial U_1}{\partial h_2} + r^{\, \prime}(V_2)\xi_2 \frac{\partial V_2}{\partial h_2} \nonumber \\
& & \hspace{4cm} + r(V_2) x_a^{\, \prime}(U_1) \frac{\partial U_1}{\partial h_2} - r(V_2)[x_c^{\, \prime}(U_2) - x_b^{\, \prime}(U_2)] \frac{\partial U_2}{\partial h_2}.
\end{eqnarray}
It follows from (\ref{dt1c1dh1}) and (\ref{dt2c1dh1}) that (\ref{djc1dh1}) can be rewritten as
\begin{eqnarray}
\label{djc1dh1a}
\frac{\partial J}{\partial h_1} & = & \psi(V_2) - \psi(V_1) + \left[ \rule{0cm}{0.4cm} \psi(V_1) - \varphi^{\, \prime}(V_1)U_1  - \psi(V_2) + \varphi^{\, \prime}(V_2)U_1  \right] t_a^{\, \prime}(U_1) \frac{\partial U_1}{\partial h_1} \nonumber \\
& = & \psi(V_2) - \psi(V_1). 
\end{eqnarray}
It follows from (\ref{dt1c1dh2}) and (\ref{dt2c1dh2}) that (\ref{djc1dh2}) can be rewritten as
\begin{eqnarray}
\label{djc1dh2a}
\frac{\partial J}{\partial h_2} & = & -\psi(V_2) + \left[ \rule{0cm}{0.4cm} \psi(V_1) - \varphi^{\, \prime}(V_1) U_1 + \varphi^{\, \prime}(V_2)U_1 - \psi(V_2) \right] t_a^{\, \prime}(U_1) \frac{\partial U_1}{\partial h_2} \nonumber \\
& & \hspace{5cm} + \left[ \rule{0cm}{0.4cm} \psi(V_2) - \varphi^{\, \prime}(V_2) U_2 \right] [t_c^{\, \prime}(U_2) - t_b^{\, \prime}(U_2)] \frac{\partial U_2}{\partial h_2} \nonumber \\
& = & - \psi(V_2).
\end{eqnarray}

{\bf \boldmath Case $2$:}  We assume the optimal strategy takes the form described in Case~$2$ of Appendix~\ref{s:itc}.  The cost of the strategy is
\begin{equation}
\label{jc2}
J = \int_0^{V_2} H(v)dx_a(v) + r(V_1) \xi_1(U_1)
\end{equation}
where $U_2$, $V_1$ and $V_2$ are independent variables and $U_1 = U_s^{\dag}(U_2,V_1,V_2)$ given by (\ref{ospsd}) is the optimal signal location speed at $x_1$.  There is one distance constraint
\begin{equation}
\label{x2c2}
x_2 - x_1 - \int_{U_1}^{V_2} dx_a(v) - \int_{U_2}^{V_2} (-1)dx_c(v) - \int_0^{U_2} (-1)dx_b(v) = 0
\end{equation}
and two time constraints
\begin{equation}
\label{t1c2}
\int_0^{U_1} dt_a(v) + \xi_1(U_1)/V_1 = h_1
\end{equation}
and
\begin{equation}
\label{t2c2}
\int_{U_1}^{V_2} dt_a(v) + \int_{U_2}^{V_2} (-1)dt_c(v) +\int_0^{U_2} (-1)dt_b(v) = h_2 - h_1.
\end{equation}
Differentiation of (\ref{x2c2}) with respect to $h_1$ and $h_2$ and using the relationship $dx(v) = vdt(v)$ gives
\begin{equation}
\label{dx2c2dh1}
 V_2[t_a^{\, \prime}(V_2) - t_c^{\, \prime}(V_2)] \frac{\partial V_2}{\partial h_1} - U_1t_a^{\, \prime}(U_1) \frac{\partial U_1}{\partial h_1} + U_2[t_c^{\, \prime}(U_2) - t_b^{\, \prime}(U_2)] \frac{\partial U_2}{\partial h_1} = 0,
\end{equation}
and
\begin{equation}
\label{dx2c2dh2}
V_2[t_a^{\, \prime}(V_2) - t_c^{\, \prime}(V_2)] \frac{\partial V_2}{\partial h_2} - U_1t_a^{\, \prime}(U_1) \frac{\partial U_1}{\partial h_2} + U_2[t_c^{\, \prime}(U_2) - t_b^{\, \prime}(U_2)] \frac{\partial U_2}{\partial h_2} = 0.
\end{equation}
Differentiation of (\ref{t1c2}) with respect to $h_1$ and $h_2$ and some elementary algebra gives
\begin{equation}
\label{dt1c2dh1}
\xi_1 \frac{\partial V_1}{\partial h_1} = -V_1^2 + V_1(V_1-U_1)t_a^{\, \prime}(U_1)\frac{\partial U_1}{\partial h_1}
\end{equation}
and
\begin{equation}
\label{dt1c2dh2}
\xi_1 \frac{\partial V_1}{\partial h_2} = V_1(V_1-U_1)t_a^{\, \prime}(U_1)\frac{\partial U_1}{\partial h_2}.
\end{equation}
Differentiation of (\ref{t2c2}) with respect to $h_1$ and $h_2$ gives
\begin{equation}
\label{dt2c2dh1}
[ t_a^{\, \prime}(V_2) - t_c^{\, \prime}(V_2)] \frac{\partial V_2}{\partial h_1} - t_a^{\, \prime}(U_1) \frac{\partial U_1}{\partial h_1} +[t_c^{\, \prime}(U_2) - t_b^{\, \prime}(U_1)] \frac{\partial U_2}{\partial h_1} = -1
\end{equation}
and
\begin{equation}
\label{dt2c2dh2}
[ t_a^{\, \prime}(V_2) - t_c^{\, \prime}(V_2)] \frac{\partial V_2}{\partial h_2} - t_a^{\, \prime}(U_1) \frac{\partial U_1}{\partial h_2} +[t_c^{\, \prime}(U_2) - t_b^{\, \prime}(U_1)] \frac{\partial U_2}{\partial h_2} = 1.
\end{equation}
Differentiation of (\ref{jc2}) with respect to $h_1$ and $h_2$ gives
\begin{equation}
\label{djc2dh1}
\frac{\partial J}{\partial h_1} = H(V_2)x_a^{\, \prime}(V_2) \frac{\partial V_2}{\partial h_1} + r^{\, \prime}(V_1) \xi_1  \frac{\partial V_1}{\partial h_1} - r(V_1)x_a^{\, \prime}(U_1) \frac{\partial U_1}{\partial h_1}
\end{equation}
and
\begin{equation}
\label{djc2dh2}
\frac{\partial J}{\partial h_2} = H(V_2)x_a^{\, \prime}(V_2) \frac{\partial V_2}{\partial h_2} + r^{\, \prime}(V_1)  \frac{\partial V_1}{\partial h_2} - r(V_1)x_a^{\, \prime}(U_1) \frac{\partial U_1}{\partial h_2}.
\end{equation}
We can combine (\ref{dx2c2dh1}) and (\ref{dt2c2dh1}) to eliminate the terms in $\partial U_2/\partial h_1$ and obtain
$$
(V_2 - U_2)[t_a^{\, \prime}(V_2) - t_c^{\, \prime}(V_2)] \frac{\partial V_2}{\partial h_1} + (U_2-U_1)t_a^{\, \prime}(U_1) \frac{\partial U_1}{\partial h_1} = U_2.
$$
Substitution of the basic formul{\ae} $t_a^{\, \prime}(v) = 1/(H(v) - r(v))$, $t_c^{\, \prime}(v) = -1/r(v)$ and $t_b^{\, \prime}(v) = (-1)/(K(v) + r(v))$ and rearrangement now gives
\begin{equation}
\label{sub1}
H(V_2)x_a^{\, \prime}(V_2) \frac{\partial V_2}{\partial h_1} = \frac{\varphi(V_2)U_2}{V_2-U_2} + \frac{\varphi(V_2)(U_1-U_2)}{V_2-U_2}  t_a^{\, \prime}(U_1) \frac{\partial U_1}{\partial h_1}.
\end{equation}
Now (\ref{dt1c2dh1}) and (\ref{sub1}) can be used to expand (\ref{djc2dh1}) and show that 
\begin{eqnarray}
\label{djc2dh1a}
\frac{\partial J}{\partial h_1} & = & \frac{\varphi(V_2)U_2}{V_2-U_2} - \psi(V_1) + \left\{ \left[ \frac{\varphi(V_2)}{V_2-U_2} - \varphi^{\, \prime}(V_1) \right] U_1 + \psi(V_1) - \frac{\varphi(V_2)U_2}{V_2-U_2} \right\} t_a^{\, \prime}(U_1) \frac{\partial U_1}{\partial h_1} \nonumber \\
& = & \frac{\varphi(V_2)U_2}{V_2-U_2} - \psi(V_1).
\end{eqnarray}
We can combine (\ref{dx2c2dh2}) and (\ref{dt2c2dh2}) to eliminate the terms in $\partial U_2/\partial h_2$ and obtain
$$
(V_2 - U_2)[t_a^{\, \prime}(V_2) - t_c^{\, \prime}(V_2)] \frac{\partial V_2}{\partial h_2} + (U_2-U_1)t_a^{\, \prime}(U_1) \frac{\partial U_1}{\partial h_2} = -U_2.
$$
Substitution of the basic formul{\ae} $t_a^{\, \prime}(v) = 1/(H(v) - r(v))$, $t_c^{\, \prime}(v) = -1/r(v)$ and $t_b^{\, \prime}(v) = (-1)/(K(v) + r(v))$ and rearrangement now gives
\begin{equation}
\label{sub2}
H(V_2)x_a^{\, \prime}(V_2) \frac{\partial V_2}{\partial h_2} = \frac{- \varphi(V_2)U_2}{V_2-U_2} + \frac{\varphi(V_2)(U_1-U_2)}{V_2-U_2}  t_a^{\, \prime}(U_1) \frac{\partial U_1}{\partial h_2}.
\end{equation}
Now (\ref{dt1c2dh2}) and (\ref{sub2}) can be used to expand (\ref{djc2dh2}) and show that 
\begin{eqnarray}
\label{djc2dh2a}
\frac{\partial J}{\partial h_2} & = & \frac{- \varphi(V_2)U_2}{V_2-U_2} + \left\{ \left[ \frac{\varphi(V_2)}{V_2-U_2} - \varphi^{\, \prime}(V_1) \right] U_1 + \psi(V_1) - \frac{\varphi(V_2)U_2}{V_2-U_2} \right\} t_a^{\, \prime}(U_1) \frac{\partial U_1}{\partial h_2} \nonumber \\
& = & \frac{- \varphi(V_2)U_2}{V_2-U_2}.
\end{eqnarray}

\subsection{The Donsker invariance principle}
\label{s:dip}

The Donsker invariance principle can be seen as a specific generalization of the central limit theorem.  Let $(\Omega, \Sigma, \mu)$ be a probability space and suppose that $\{\xi_k\}_{k \in {\mathbb N}}$ is a sequence of independent and identically distributed random variables on $(\Omega, \Sigma,\mu)$ with mean ${\mathbb E}[ \xi_k ] = 0$ and variance ${\mathbb E}[ \xi_k^2] = 1$.  Let $S_m = \sum_{k=1}^m \xi_k$ for all $m \in {\mathbb N}$ and define a random continuous function $X_n:C[0,1] \times \Omega \rightarrow {\mathbb R}$ by the formula
$$
X_n(t,\omega) = \left(1/\sqrt{n}\, \right) \mbox{$\sum_{1 \leq k \leq \lfloor nt \rfloor}$} \xi_k(\omega) + (nt - \lfloor nt \rfloor) \left(1/\sqrt{n}\, \right) \xi_{ \lfloor nt \rfloor +1}(\omega).
$$
The function $X_n(t,\omega)$ is a piecewise linear approximation to the Wiener distribution.  When $t = 1$ we have
$$
X_n(1,\omega) = \left(1/\sqrt{n}\, \right) \mbox{$\sum_{1 \leq k \leq n}$} \xi_k(\omega) \rightarrow z(\omega) \sim {\mathcal N}(0,1)
$$
as $n \rightarrow \infty$ by the central limit theorem.  When $t = p/q$ where $p, q \in {\mathbb N}$ and $0 \leq p < q$ we have
$$
X_{nq}(p/q, \omega) = \left(1/\sqrt{nq}\, \right) \mbox{$\sum_{1 \leq k \leq np}$} \xi_k(\omega) \rightarrow \sqrt{p/q} \cdot z(\omega) \sim {\mathcal N} \left(0, p/q\, \right)
$$
as $n \rightarrow \infty$, once again by the central limit theorem.  The Donsker invariance principle \cite{don1, dud1} extends this convergence\textemdash in distribution\textemdash to all $t \in [0,1]$ with $X_n(t, \omega) \rightarrow \sqrt{t} \cdot z(\omega) \sim {\mathcal N}(0, t)$ as $n \rightarrow \infty$.  The simpler step function $S_n(t,\omega)$ can be used in place of the continuous function $X_n(t,\omega)$ in practice because $\|X_n(t,\omega) - S_n(t,\omega)\|_{\infty} \rightarrow 0$ as $n \rightarrow \infty$.  The usual notational convention is to suppress the dependence on $ \omega$.  Thus, for a scaled Wiener process $\epsilon(t)$ with scale parameter $\theta$ we will normally write $\theta X_n(t) \rightarrow \epsilon(t) = \theta W(t) = \theta \sqrt{t} \cdot z \sim {\mathcal N}(0,\theta^2 t)$.

\end{document}